\documentclass[12pt,leqno]{amsart}
\usepackage{graphicx}
\usepackage{tgtermes}

\usepackage{indentfirst,csquotes}

\usepackage[utf8]{inputenc}
\usepackage{amsmath}
\usepackage{amsfonts}
\usepackage{amssymb}
\usepackage[T1]{fontenc}
\usepackage[english]{babel}
\usepackage{fancyhdr}
\usepackage{indentfirst}
\usepackage{graphicx}
\usepackage{newlfont}
\usepackage{amssymb}
\usepackage{latexsym}
\usepackage{amsthm}
\usepackage{mathtools}
\usepackage{subfig}
\usepackage{floatflt}
\usepackage{floatrow}
\usepackage{color}
\usepackage{textcomp}
\usepackage{framed}
\usepackage{verbatim}
\usepackage{hyperref}
\usepackage{booktabs}
\usepackage{caption}
\usepackage[mathscr]{euscript}
\usepackage{braket}
\usepackage{bbm}
\usepackage{mathrsfs}
\newcommand*{\transp}[2][-3mu]{\ensuremath{\mskip1mu\prescript{\smash{\mathrm t\mkern#1}}{}{\mathstrut#2}}}%

\DeclarePairedDelimiter{\abs}{\lvert}{\rvert}
\DeclarePairedDelimiter{\norma}{\lVert}{\rVert}

\theoremstyle{plain}
\newtheorem{theorem}{Theorem}[section]     
\newtheorem{proposition}[theorem]{Proposition}    
\newtheorem{corollary}[theorem]{Corollary}       
\newtheorem{lemma}[theorem]{Lemma}            
\theoremstyle{definition}
\newtheorem{definition}[theorem]{Definition}	
\newtheorem{example}[theorem]{Example}     
\newtheorem{remark}[theorem]{Remark} 
 
\numberwithin{equation}{section}

\begin{document} 
\title[] {Propagation of Shubin-Sobolev singularities of Weyl-quantizations of complex quadratic forms}
\author[]{Marcello Malagutti, Alberto Parmeggiani and Davide Tramontana}

\address{Marcello Malagutti
\newline \indent University College London, Department of Mathematics
\newline \indent 25 Gordon St., London, WC1H 0AY, United Kingdom}
\email {m.malagutti@ucl.ac.uk}

\address{Alberto Parmeggiani
\newline \indent Dipartimento di Matematica, Universit\`a di Bologna
\newline \indent Piazza di Porta San Donato 5, 40126, Bologna, Italy}
\email{alberto.parmeggiani@unibo.it}

\address{Davide Tramontana
\newline \indent Dipartimento di Matematica, Universit\`a di Bologna
\newline \indent Piazza di Porta San Donato 5, 40126, Bologna, Italy}
\email{davide.tramontana4@unibo.it}

\thanks{The authors are members of the Research Group GNAMPA of INdAM}
\thanks{The first author acknowledges the support by  EPSRC Early Career Fellowship: EP/V001760/1.}
\thanks{The second author was partially supported by the Italian Ministry of University and Research, under PRIN2022 (Scorrimento) “Anomalies
in partial differential equations and applications”, 2022HCLAZ8\_002, J53C24002560006.}

\begin{abstract}
The aim of this work is to develop the H\"ormander microlocal theory in the isotropic framework and use the 
results we obtain to study the propagation of singularities for an evolution problem, with generator given by a Weyl-quantization of a complex quadratic form on the phase space.  
\end{abstract} 
\maketitle
\tableofcontents                        
\rhead[\fancyplain{}{\bfseries\thepage}]{\fancyplain{}{\bfseries\thepage}}
\lhead[\fancyplain{}{\bfseries\thepage}]{\fancyplain{}{\bfseries\thepage}}

\section{Introduction}            
This paper is motivated by two main goals: a generalization of H\"ormander's microlocal results of \cite{HQ} in the isotropic framework and the consequent study of the propagation of singularities
for an evolution equation of the form 
\begin{equation}\label{problemintro}
\begin{cases}
\partial_t u+Au=0, \quad \text{in} \quad \mathbb{R}^+ \times \mathbb{R}^n, \\
u(0,\cdot)=u_0 \in \mathscr{S}'(\mathbb{R}^{n}),
\end{cases}
\end{equation}
with $A=\mathrm{Op}^{\mathrm{w}}(a)\in \Psi_{\mathrm{iso}}^2(\mathbb{R}^n)$ the Weyl-quantization (see Section \ref{s:PrelimIsoCalc}) of a complex quadratic form on the phase space
\[
a:\; \mathbb{R}^{2n} \longrightarrow \mathbb{C}, \quad
X \longmapsto a(X)=\langle X,QX \rangle,
\]
defined
by the symmetric complex matrix
$Q \in \mathsf{M}_{2n}(\mathbb{C})$ with $\mathrm{Re} \, Q \geq 0$. The existence of a solution operators $\lbrace e^{-tA} \rbrace_{t \geq 0}$ for such a problem may be found in \cite{HS} 
and we will discuss this matter in more detail in Section \ref{sec.pos}.

The importance of studying such isotropic theory comes from the fact that it allows us to study a broader class of meaningful problems.
Indeed, every evolution problem
in the whole of $\mathbb{R}^n$ for which the generator is given by an operator whose symbol is a quadratic form with positive real part - such as the Heat equation and the Schr\"odinger 
equation - are included in our theory. We will also study the particular case of more general operators tailored to a complex combination of \textit{quantum Harmonic Oscillators}, such as
\[
A=-\Delta_{x'}+\abs{x'}^2+i(-\Delta_{x''}+\abs{x''}^2),
\]
where $x'\in \mathbb{R}^{n_1}$ and $x''\in \mathbb{R}^{n_2}$, with $n_1+n_2=n$ (note that if $n_2=n$ we recover the standard Harmonic Oscillator). 
This operator has a symbol that falls within the scope of our study (see Theorem \ref{theoHO} and Example \ref{exaHO}) and, as we expect, since the Heat equation has a regularizing behavior, 
we have that the singularities of the solution of the problem at time $t>0$ live in the second group of variables (in other words, there is a heat diffusion with respect to some variables, 
and a quantum diffusion in some others.) 

In order to extend H\"ormander's results we will recall the isotropic calculus introduced by Shubin (see \cite{SH}), that is the calculus with respect to the \textit{isotropic} H\"ormander metric 
\[
g_{0,X}=\frac{\abs{dX}^2}{\langle X \rangle^2}, \quad X=(x,\xi)\in \mathbb{R}^{2n},
\]
and we will provide the main facts needed in this framework, based on the notion of isotropic wave front set introduced in \cite{HQ} (see also \cite{DGW}).

We study the "stratification" of $WF_{\mathrm{iso}}$ given by the $s$-wave-front sets $WF^s_{\mathrm{iso}}$,
that are a measure of microlocal $B^s$-regularity, where $B^s$ is the Shubin-Sobolev space on $\mathbb{R}^n$ (see Definition \ref{defBs} below). 
The sets $WF_{\mathrm{iso}}^s$ allow for a more precise information about the regularity of the distribution that we are dealing with. Indeed, one of the reasons to study such stratification arises 
from the fact that once the results are achieved in the stratified setting, it is possible to reconstruct the corresponding properties in the "Schwartz regularity" setting (see Proposition \ref{WF&WFs}).

To obtain such results we need to use different global metrics in addition to the isotropic one, that is the \textit{iso,iso} metric on $\mathbb{R}^{2(n+m)}$ and the $SG$ 
metric on $\mathbb{R}^{2n}$ (see Definition \ref{defSG_isoiso} and  Definition \ref{defSG}).
 Furthermore, we will examine the relation with the \textit{Time-Frequency} calculus introduced in \cite{G}, and developed by many authors (see e.g. \cite{NR}, \cite{RW}, 
\cite{SW}; see also \cite{CM} for an SG-situation). In that context, we will show that the tools we have introduced are related to a complementary perspective that, in the end, 
gives an essentially equivalent information, however through a more pseudodifferential (and flexible) viewpoint, that in many cases is better adapted to phase-space geometry.
 
We finally use our results to study the propagation of singularities, aiming to relate the Shubin-Sobolev singularity of the solution of the problem \eqref{problemintro} 
at time $t$ to that of the initial datum.

The propagation of singularities problem \eqref{problemintro} was already studied from a time-frequency point of view using, as we do, the symplectic geometry of the quadratic form.
To explain briefly the result on propagation of singularities we begin by mentioning that, thanks to the symplectic classification of quadratic form due to H\"ormander \cite{HS}, 
we are able to give the first rough inclusion 
(see \cite{HQ}, \cite{HS})
\[
\mathrm{WF}_{\mathrm{iso}}(e^{-tA}u_0) \subseteq e^{2t\mathrm{Im}\,F} \mathrm{WF}_{\mathrm{iso}}(u_0),
\]  
where $F$ is the complex matrix that satisfies the equation (with $\sigma$ the standard complex symplectic form) 
\[
\sigma(X,FX)=a(X),
\]
that is called the \textit{Hamilton map} associated with $a$. In what follows we will write $a=a_R+ia_I,$ where $a_R=\mathrm{Re}\,a$ and $a_I=\mathrm{Im}\,a.$
The goal is then to improve the inclusion above.

First of all, by using the geometry of the singular space $S$, introduced by M. Hitrik and K. Pravda-Starov in \cite{HPS}, the authors of \cite{PRW} proved that, for all $t>0$
\begin{equation}\label{incPRW}
WF_{\mathrm{iso}}(e^{-tA}u_0) \subseteq (e^{tH_{a_I}}(WF_{\mathrm{iso}}(u_0) \cap S))\cap S,
\end{equation}
that essentially says that the $\mathscr{S}$-singularities of the solution at time $t$ are localized inside the singular space and propagate along the bicharacteristic curves 
of the imaginary part $\mathrm{Im} \, a$ of the symbol.
In the case of $WF_{\mathrm{iso}}^s$, i.e. for studying the 
\textit{Shubin-Sobolev} singularities, the inclusion \eqref{incPRW} was already studied by using again the time-frequency tools in \cite{W} in which one has a different loss of
derivatives. 
In recent years, the inclusion has also been studied for an initial condition in the Gelfand-Shilov class and in the anisotropic context (see \cite{CW}, \cite{CRW}, \cite{W1}, \cite{Wh}). 

\noindent We will also obtain this inclusion from our complementary point of view by using the isotropic microlocal theory developed in Section \ref{s:IsoMicroSing} and in Section \ref{s:MainProp}.

In particular, in order to study the singularities of a tensor product of two tempered distributions we have to keep into account the global regularity of both of them and this produces 
a first loss of derivatives (the $s_\ast$ in Theorem \ref{WFmainprope}) due to the interaction of the single tempered distributions $u$ and $v$ in the tensor product 
$u\otimes v\in \mathscr{S}'(\mathbb{R}^n_x\times\mathbb{R}^m_y)$. Moreover, when we consider the "pullback" of a temperate distribution - that satisfies suitable conditions - 
through an injective map from $\mathbb{R}^m$ to $\mathbb{R}^n$ we encounter an additional loss of derivatives, that depends on the codimension $n-m$, since it basically corresponds 
to taking the trace of the distribution. In fact, we reduce matters to the case in which the pullback is done via an immersion $L:\mathbb{R}^m\to\mathbb{R}^m\times \{0_{n-m}\}$, 
and after making sense of the notion of the pullback via an injective map of distributions in a certain class, we show that on sufficiently regular distributions the pullback is 
the trace on $\{x_{m+1}=\dots=x_n=0\}$ and the Trace Theorem leads to a loss of regularity of $(n-m)/2$.
\\ 
We next describe our main result on the Schwartz kernels that plays a crucial role in our work. Let $\mathscr{K}:\mathscr{S}(\mathbb{R}^n)\rightarrow \mathscr{S}'(\mathbb{R}^m)$,
with $n \geq m$, be a linear operator with Schwartz kernel $K\in \mathscr{S}'(\mathbb{R}^{n+m})$. Then, if $K\in \mathscr{S}(\mathbb{R}^{n+m})$ and $u \in \mathscr{S}(\mathbb{R}^{n})$, we may write 
\[
\mathscr{K}u(x)=\int L^\ast (K \otimes u)(x,y)dy,
\]
where $K\otimes u\in \mathscr{S}(\mathbb{R}^{m+n+n})$ is the function 
\[
(x,y,z) \mapsto K(x,y)u(z),
\]
and $L$ is the map $(x,y)\mapsto (x,y,y)$.
Therefore, after extending this expression to those $K\in \mathscr{S}'(\mathbb{R}^{m+n})$ and $u \in \mathscr{S}'(\mathbb{R}^n)$ that satisfy condition \eqref{condWF2} below 
(see Section \ref{s:MainProp} for the notation), we put together the previous results into the following theorem. 

\begin{theorem}
\label{WFmainprope}
Let $\mathscr{K}:\mathscr{S}(\mathbb{R}^n)\rightarrow \mathscr{S}'(\mathbb{R}^m)$,
with $n \geq m$,
be a linear operator with Schwartz kernel $K \in B^{-r_1}(\mathbb{R}^{m+n})$ 
for $r_1\geq 0$.
Then, if $u \in \mathscr{S}(\mathbb{R}^n)$
and $\mu>n$ is fixed, we have
that for all $s_1 \in \mathbb{R}$ 
\begin{equation}\label{WFmainpropeeq1}
WF_{\mathrm{iso}}^{s_1-\mu}(\mathscr{K}u)\subseteq WF_{\mathrm{iso},X}^{s_1}(K).
\end{equation}
 
\noindent In addition, the definition of the map $u \mapsto \mathscr{K}u$ 
can be extended by continuity to those $u \in \mathscr{S}'(\mathbb{R}^n)$ 
for which
\begin{equation}\label{condWF2}
WF_{\mathrm{iso}}(u) \cap WF_{\mathrm{iso},Y}(K)=\emptyset.
\end{equation}
Finally, if $u\in B^{-r_2}(\mathbb{R}^n)$, for some $r_2 \geq 0$, satisfies
\eqref{condWF2} 
and $s_1,s_2 \in \mathbb{R}$ are such that 
$s_\ast:=\min\lbrace s_1-r_2, s_2-r_1\rbrace\leq s_1+s_2$, 
then for any fixed $\mu>n$ we have  
\begin{equation}\label{WFmainpropeeq2}
WF_{\mathrm{iso}}^{s_\ast-\mu}(\mathscr{K}u)\subseteq WF_{\mathrm{iso},X}^{s_1}(K)\cup WF_{\mathrm{iso}}^{s_1}(K)' \circ WF_{\mathrm{iso}}^{s_2}(u).
\end{equation}
\end{theorem}
\vspace{2mm}

By using the geometric structure of the kernel $K_{e^{-2itF}}$ of the propagator $\mathscr{K}_{e^{-2itF}}=e^{-tA}$ (see \cite{HS}) - and of the associated singular space $S$ - 
we will specialize Theorem \ref{WFmainprope} to improve the inclusion \eqref{incPRW} in the stratified setting giving in addition a lower bound on the loss of derivatives of 
$4n+\varepsilon$, with $\varepsilon$ as small as we wish. This will be our main result concerning the propagation of singularities, stated as follows.  
 
\begin{theorem}\label{maintheo}
 Let $t>0$ be fixed and consider the operator $\mathscr{K}_{e^{-2itF}}:\mathscr{S}(\mathbb{R}^n)\rightarrow \mathscr{S}'(\mathbb{R}^n)$.
 Then, for any given $\varepsilon>0$ (small), the Schwartz kernel $K_{e^{-2it'F}}$ of $\mathscr{K}_{e^{-2it'F}}$ belongs to 
$B^{-(n+\varepsilon)}(\mathbb{R}^{2n})$, uniformly in $t' \in (0,t]$.\\
 Moreover, on defining $r_0=\inf\{r\geq 0;\,\,K_{e^{-2it'F}}\in B^{-r}(\mathbb{R}^{2n}),\,\,\forall t'\in (0,t]\}$, one has that 
 if $u_0 \in \mathscr{S}'(\mathbb{R}^n)$ and $\mu>2(n+r_0)$, then for all $s \in \mathbb{R}$ 
 \begin{equation}\label{incmaintheo}
 WF_{\mathrm{iso}}^{s-\mu}(e^{-tA}u_0) \subseteq (e^{tH_{a_I}}(WF_{\mathrm{iso}}^{s}(u_0) \cap S))\cap S.
 \end{equation}
\end{theorem}
\vspace{3mm}

We conclude this introduction by giving the plan of the paper. In Section \ref{s:PrelimIsoCalc} we recall the isotropic calculus. Namely, we define the class of symbols 
whose quantization leads to the class of the isotropic pseudodifferential operators ($\psi$dos) that we will use in the following sections. \\
The aim of Section \ref{s:IsoMicroSing} is to define the isotropic wave front sets and relate them with the action of isotropic $\psi$dos on tempered distributions. 
We will give the notion of elliptic and characteristic set and
 show the basic microlocal properties in this isotropic framework, that will be fundamental to study the microlocal regularity. 

Section \ref{s:TFCalc} is devoted to the relation between the isotropic pseudodifferential calculus and the time-frequency analysis. In fact, we recall the definition 
and some basic properties of the Short-Time Fourier Transform (STFT), and use it to introduce the class of anti-Wick quantized operators, 
that will be useful to define modulation spaces and the Gabor wave front set. The class of anti-Wick quantized operators leads also to the fact that this wave front set 
coincides with the isotropic one of Section \ref{s:IsoMicroSing}. 

In Section \ref{s:MainProp} 
 we extend H\"ormander's results concerning the operations for the wave front set, in terms of the $WF_{\mathrm{iso}}^s$
(see \cite{HoV3} for the $C^\infty$ case and \cite{HQ} for the $\mathscr{S}$ case) needed to prove Theorem \ref{WFmainprope}.

Finally, the goal of Section \ref{sec.pos} is to give the results on the propagation of singularities. Following \cite{HS}, we recall the notion of Gaussian distribution 
through which we will provide a geometric characterization for the kernel of the propagator of \eqref{problemintro}. This will lead to the propagation results and, in particular, 
to the proof of Theorem \ref{maintheo}.

An Appendix where we collect (and prove) some useful results, concludes the paper.

\section{Preliminaries of isotropic calculus} \label{s:PrelimIsoCalc}
In this section we recall some basic notions of microlocal analysis and pseudodifferential calculus in the isotropic context. Throughout, we will work on 
$T^\ast \mathbb{R}^n\cong \mathbb{R}^{2n}=\mathbb{R}^n\times\mathbb{R}^n$. Hence, we write $X=(x,\xi) \in \mathbb{R}^{2n}$ for the points in the phase-space and 
we write $\dot{\mathbb{R}}^{2n}=\mathbb{R}^{2n}\setminus \lbrace 0 \rbrace$. 
Recall also that $\langle \cdot \rangle=(1+\abs{\cdot}^2)^{1/2}$, $D_{x_j}=-i\partial_{x_j}$ and $(\cdot,\cdot)_0$ is the $L^2$-product on $\mathbb{R}^n$.
Furthermore we denote by $\mathbb{N}_{0}=\mathbb{N}\cup \lbrace 0 \rbrace$, by $\lbrace 0_n \rbrace$ the origin on $\mathbb{R}^n$
and, finally, by $\langle \cdot, \cdot \rangle$ the \textit{non-Hermitian} product on $\mathbb{C}^{2n}$ defined by $\langle X,Y \rangle:=\sum_iX_iY_i$ for $X,Y \in \mathbb{C}^{2n}$. 
We say that a set $\Gamma\subseteq \dot{\mathbb{R}}^{2n}$ is \textit{conic} if it is invariant under multiplication by positive real numbers. Hence a cone 
(i.e. a conic set) $\Gamma \subseteq \dot{\mathbb{R}}^{2n}$ is uniquely determined by its intersection with the unit sphere $\mathbb{S}^{2n-1}$. Therefore $\Gamma$ is open, resp. closed, 
in $\dot{\mathbb{R}}^{2n}$ if $\Gamma$ is relatively open, resp. closed, in $\dot{\mathbb{R}}^{2n}$ (equivalently, $\Gamma\cap\mathbb{S}^{2n-1}$ is open, resp. closed, in $\mathbb{S}^{2n-1}$).

Next, we briefly recall the isotropic calculus introduced by Shubin (in the context of the Weyl-H\"ormander calculus). We refer to \cite{HoV3} and \cite{SH} (see also \cite{Le}). 

\begin{definition}
Let $a \in C^\infty(\mathbb{R}^{2n})$ and $m \in \mathbb{R}$. We say that \textit{$a$ is a Shubin symbol of order $m$}, and we write $a \in S_{\mathrm{iso}}^m(\mathbb{R}^{n})$, 
if for all $\alpha \in \mathbb{N}_0^{2n}$  there exists $C_{\alpha}>0$ such that
\begin{equation}
\label{Shubinsymbol}
\abs*{\partial_X^{\alpha}a(X)}\leq C_{\alpha}\braket{X}^{m-\abs*{\alpha}}, \quad X \in \mathbb{R}^{2n}.
\end{equation}
\end{definition}

Shubin symbols of order $m$ form a Fr\'echet space where, if $a \in S_{\mathrm{iso}}^m(\mathbb{R}^n)$, the semi-norms are given by
\begin{equation}\label{semi-norms}
\abs{a}_k^{(m)}:=\max_{\abs{\alpha}\leq k}\sup_{X \in \mathbb{R}^{2n}}\abs{\partial_X^\alpha a(X) }\langle X \rangle^{-(m-\abs{\alpha})}, \quad k \in \mathbb{N}_0.
\end{equation}
We denote
\[
S_{\mathrm{iso}}^{\infty}(\mathbb{R}^{n})=\bigcup_{m \in \mathbb{R}} S^m_{\mathrm{iso}}(\mathbb{R}^{n}).
\]
Note also that
\[
\bigcap_{m \in \mathbb{R}}S_{\mathrm{iso}}^m(\mathbb{R}^{n})=\mathscr{S}(\mathbb{R}^{2n}).
\]

If $a \in S^m_{\mathrm{iso}}(\mathbb{R}^n)$, we can associate
to it a \textit{pseudodifferential operator}
using the so-called \textit{quantization formulas}.
The different quantizations we will use belong to the one-parameter 
family of quantizations known as $t$-quantizations, defined as follows.  

\begin{definition}
Let $t \in [0,1]$ and $a \in S^m_\mathrm{iso}(\mathbb{R}^n)$. The \textit{$t$-quantized pseudodifferential operator} associated with $a$ is defined as
\[
\mathrm{Op}_t(a)u(x)= (2\pi)^{-n}\int_{\mathbb{R}^{2n}}e^{i\langle x-y, \xi \rangle} a((1-t)x+ty,\xi)u(y)dyd\xi, \quad u \in \mathscr{S}(\mathbb{R}^n).
\]
\end{definition}

It is easy to verify that $\mathrm{Op}_t(a)$ 
is a linear operator that acts continuosly on 
$\mathscr{S}(\mathbb{R}^n)$ and extends by duality
to a linear continuous operator on $\mathscr{S}'(\mathbb{R}^n)$.

Different quantizations generally lead to different operators 
and the most important ones are the case
for which $t=0$ (left-quantization) 
or $t=1/2$ (Weyl-quantization).
In these cases, if $a \in S_{\mathrm{iso}}^m(\mathbb{R}^n)$, 
we define by $\mathrm{Op}(a)=\mathrm{Op}_0(a)$
\textit{the left-quantized pseudodifferential operator} associated
with it, given by
\[
\mathrm{Op}(a)u(x)= (2\pi)^{-n}\int_{\mathbb{R}^{2n}}e^{i\langle x-y, \xi \rangle} a(x,\xi)u(y)dyd\xi, \quad u \in \mathscr{S}(\mathbb{R}^n)
\]
and by $\mathrm{Op}^{\mathrm{w}}(a)=\mathrm{Op}_{1/2}(a)$ the 
\textit{Weyl-quantized pseudodifferential operator}, given by 
\[
\mathrm{Op}^{\mathrm{w}}(a)u(x)= (2\pi)^{-n}\int_{\mathbb{R}^{2n}}e^{i\langle x-y, \xi \rangle} a((x+y)/2,\xi)u(y)dyd\xi, \quad u \in \mathscr{S}(\mathbb{R}^n).
\]
Given $a \in S_{\mathrm{iso}}^m(\mathbb{R}^n)$ and $u \in \mathscr{S}'(\mathbb{R}^n)$ we sometimes use also the notation 
\[
a(x,D)u:=\mathrm{Op}(a)u, \quad a^\mathrm{w}(x,D)u:=\mathrm{Op}^\mathrm{w}(a)u.
\]
In what follows, by the convenient properties
of the Weyl Calculus, such as the adjoint property (see Theorem \ref{adjA*})
and the metaplectic invariance (see \eqref{Metainv}),
we focus on the case of Weyl-quantized 
pseudodifferential operators.
Therefore, we define the set of pseudodifferential operators
that we will deal with as
\[
\Psi_{\mathrm{iso}}^m(\mathbb{R}^n):=\Bigl\lbrace A:\mathscr{S}'(\mathbb{R}^n) \rightarrow \mathscr{S}'(\mathbb{R}^n); 
\ \exists \, a \in S_{\mathrm{iso}}^m(\mathbb{R}^n), \; A=\mathrm{Op}^{\mathrm{w}}(a) \Bigr\rbrace.
\]
Nevertheless, given $A \in \mathrm{Op}^\mathrm{w}(a)\in \Psi_{\mathrm{iso}}^m(\mathbb{R}^n)$,
it is well known (see for instance \cite{SH}, Theorem 23.2)
that for all $t \in [0,1]$
there exists $a_t \in S_{\mathrm{iso}}^m(\mathbb{R}^n)$
such that 
\begin{equation}\label{wt}
\mathrm{Op}^\mathrm{w}(a)=\mathrm{Op}_t(a_t),
\end{equation}
and hence, the set of pseudodifferential operators
could be equivalently described as
\[
\Psi_{\mathrm{iso}}^m(\mathbb{R}^n)=\lbrace A:\mathscr{S}'(\mathbb{R}^n) \rightarrow \mathscr{S}'(\mathbb{R}^n); \ \exists \, a \in S_{\mathrm{iso}}^m(\mathbb{R}^n), \; A=\mathrm{Op}(a) \rbrace,
\]
or more generally as
\[
\Psi_{\mathrm{iso}}^m(\mathbb{R}^n)=\lbrace A:\mathscr{S}'(\mathbb{R}^n) \rightarrow \mathscr{S}'(\mathbb{R}^n); \ \exists \, a \in S_{\mathrm{iso}}^m(\mathbb{R}^n), \; A=\mathrm{Op}_t(a) \rbrace.
\]
If $A=\mathrm{Op}^{\mathrm{w}}(a)=\mathrm{Op}(a_L) \in \Psi_{\mathrm{iso}}^m(\mathbb{R}^n)$
we call $a,a_L \in S_{\mathrm{iso}}^m(\mathbb{R}^n)$ 
\textit{the Weyl symbol} and \textit{the left symbol} of $A$,
respectively.

\begin{definition}
If $a \in \mathscr{S}(\mathbb{R}^{2n})$, then $\mathrm{Op}^\mathrm{w}(a):\mathscr{S}'(\mathbb{R}^n)\rightarrow \mathscr{S}(\mathbb{R}^n)$ 
is called \textit{a smoothing operator} and we write $A \in \Psi_\mathrm{iso}^{-\infty}(\mathbb{R}^n)$.
\end{definition}

\begin{definition}\label{DefClassic}
Let $(a_j)_{j\geq 1}$ be a sequence of symbols such that $a_j \in S_{\mathrm{iso}}^{m_j}(\mathbb{R}^{n})$, $m_j \searrow -\infty$ as $j \rightarrow + \infty$ 
and $a \in S_{\mathrm{iso}}^{m_1}(\mathbb{R}^{n})$. We write
\[
a \sim \sum_{j \geq 1}a_j,
\] 
if for any integer ${r \geq 2}$
\[
a-\sum_{j=1}^{r-1} a_j \in S_{\mathrm{iso}}^{m_r}(\mathbb{R}^{n}),
\]
where $m_r=\max_{j \geq r} m_j$. This is called \textit{an asymptotic expansion}.
\end{definition}

\begin{lemma}
\label{constofa}
Let $(a_j)_{j\geq 1}$ be a sequence of symbols such that $a_j \in S_{\mathrm{iso}}^{m_j}(\mathbb{R}^{n})$ and $m_j \searrow -\infty$ as $j \rightarrow + \infty$. 
Then there exists a function $a \in S_{\mathrm{iso}}^{m_1}(\mathbb{R}^{n})$ such that 
\[
a \sim \sum_{j \geq 1}a_j.
\]
Moreover, if another function $a'$ has the same property, then $a-a' \in \mathscr{S}(\mathbb{R}^{2n})$.
\end{lemma}

\begin{theorem}\label{compAB}
Let $a \in S_{\mathrm{iso}}^{m_1}(\mathbb{R}^{n})$ and $b \in S_{\mathrm{iso}}^{m_2}(\mathbb{R}^{n})$. Then the composition $\mathrm{Op}^{\mathrm{w}}(a)\mathrm{Op}^{\mathrm{w}}(b)$ is well defined and
\[
\mathrm{Op}^{\mathrm{w}}(a)\mathrm{Op}^{\mathrm{w}}(b)=\mathrm{Op}^{\mathrm{w}}(a \sharp b),
\]
where the symbol 
\[
(a \sharp b)(X) =\Bigl(e^{i\sigma(D_X;D_Y)/2}a(X)b(Y)\Bigr)\Bigl|_{Y=X} \in S_{\mathrm{iso}}^{m_1+m_2}(\mathbb{R}^{n})
\]
is called the Weyl-product and the map
\[
\sharp : S_{\mathrm{iso}}^{m_1}(\mathbb{R}^{n}) \times S_{\mathrm{iso}}^{m_2}(\mathbb{R}^{n}) \longrightarrow S_{\mathrm{iso}}^{m_1+m_2}(\mathbb{R}^{n}), 
\] 
is continuous. Moreover,
\[
(a \sharp b) (X) \sim \sum_{j=0}^{+\infty}\frac{(i/2)^j}{j!}\Bigl(\sigma(D_X;D_Y)\Bigr)^ja(X)b(Y)\Bigl|_{Y=X},
\]
where for $X=(x,\xi), Y=(y,\eta) \in \mathbb{R}^{2n}$,
\[
\sigma(D_X;D_Y):=\sigma(D_x,D_\xi;D_y,D_\eta)=\langle D_\xi,D_y \rangle - \langle D_x,D_\eta \rangle.
\] 
\end{theorem}
\begin{theorem}\label{adjA*}
Let $A \in \Psi_\mathrm{iso}^m(\mathbb{R}^n)$. Then the adjoint operator $A^\ast \in \Psi_\mathrm{iso}^m(\mathbb{R}^n)$ and $A^\ast=\mathrm{Op}^\mathrm{w}(a^\ast)$, where $a^\ast=\bar{a}$.  
\end{theorem}
 
We conclude this section by giving the notion of \textit{Hamilton vector field} 
and, in case of $a$ is a quadratic form, the relation 
with the so called \textit{Hamilton map}. 

Let $\sigma$ be the standard symplectic form on $\mathbb{R}^{2n}$: in coordinates $(x_1,...,x_n,\xi_1,...,\xi_n)$ it is given by $\sigma=\sum_j d\xi_j \wedge dx_j$.

\begin{definition}\label{def.HV}
\textit{The Hamilton vector field} of $a\in C^\infty(\mathbb{R}^{2n})$ 
is defined as the unique vector field $H_{a}\in T\mathbb{R}^{2n}$ such that
\[
d a(V)=\sigma(V,H_{a}), \quad \forall V \in T\mathbb{R}^{2n}.
\]
\end{definition} 

Therefore in coordinates the Hamilton vector field of $a$ takes the form
\[
H_{a}=\sum_{j=1}^n \Bigl(\partial_{\xi_j} a\cdot \partial_{x_j}-\partial_{x_j} a\cdot \partial_{\xi_j}\Bigr).
\]

Finally, we give the notion of \textit{Hamilton map} - that is the linearized version near a critical point of the Hamilton vector field (see \cite{HoV3}, Section 21.5) - 
in case of $a$ is a complex quadratic form on the phase space. More precisely, let
\[
a:\; \mathbb{R}^{2n} \longrightarrow \mathbb{C}, \quad
X \longmapsto a(X)=\langle X,QX \rangle,
\]
be a \textit{complex}
quadratic form defined
by the symmetric matrix
$Q \in \mathsf{M}_{2n}(\mathbb{C})$.
We define \textit{the Hamilton map} associated with $a=a(X)$ 
as the matrix $F$ defined by the equation
\begin{equation}\label{Hamiltonmap}
\sigma(X,FX)=a(X),
\end{equation}
where $\sigma$ is the symplectic form of $\mathbb{R}^{2n}$ extended to $\mathbb{C}^{2n}$ defined by 
\[
\sigma((x,\xi),(y,\eta))=\langle y,\xi \rangle- \langle x,\eta \rangle. 
\]
It is useful to observe that the matrix $F$ can be expressed also as $F=JQ$, where 
\[
J=
\begin{pmatrix}
0 & I_n \\
-I_n & 0
\end{pmatrix}.
\]
Therefore, if $a$ is a \textit{real-valued} quadratic form, one has
\begin{equation}
\label{lin.HVF}
H_{a}(x,\xi)=
\begin{pmatrix}
\partial_\xi a(x,\xi) 
\vspace{2mm} \\
-\partial_x a(x,\xi)
\end{pmatrix}
=2F
\begin{pmatrix}
x \\
\xi
\end{pmatrix}, \quad (x,\xi) \in \mathbb{R}^{2n}.
\end{equation}
\section{Isotropic microsingularities} \label{s:IsoMicroSing}
The aim of this section is to give the necessary tools to study the isotropic microsingularities on the phase space (see also \cite{H} and \cite{DGW}). 
To do that we first define a characteristic set for the isotropic class of symbols.

\begin{definition}\label{defprincsymb}
The operator $A=\mathrm{Op}^{\mathrm{w}}(a) \in \Psi_{\mathrm{iso}}^m(\mathbb{R}^{n})$ is said to be \textit{isotropic elliptic at $X_0\in \dot{\mathbb{R}}^{2n}$} 
if there exist an open cone $\Gamma_{X_0}\subseteq \dot{\mathbb{R}}^{2n}$ that contains $X_0$ and some constants $c,C>0$, such that 
\[
\abs*{a(X)}\geq c\abs*{X}^m, \quad \forall X\in \Gamma_{X_0}, \ \abs*{X}\geq C.
\]
If $A$ is isotropic elliptic at $X_0$, for every $X_0 \in \dot{\mathbb{R}}^{2n}$ it is said to be \textit{isotropic elliptic}. If $A$ is not isotropic elliptic at $X_0$, 
it is said to be \textit{isotropic characteristic at $X_0$}.
\end{definition}

\begin{definition}
We define the set of the \textit{isotropic elliptic points} for $A$ as the set 
\[
\mathrm{Ell}_{\mathrm{iso}}(A)=\lbrace X_0 \in \dot{\mathbb{R}}^{2n}; \ A \ \mathrm{is} \ \mathrm{elliptic} \ \mathrm{at} \ X_0 \rbrace,
\]
and the set of \textit{isotropic characteristic points} as the set
\[
\mathrm{Char}_{\mathrm{iso}}(A)=\dot{\mathbb{R}}^{2n} \setminus \mathrm{Ell}_{\mathrm{iso}}(A).
\]
With $A=\mathrm{Op}^\mathrm{w}(a)$, we say that $a$ is \textit{non-characteristic} at $X_0 \in \mathbb{R}^{2n}$ if $X_0 \notin \mathrm{Char}_{\mathrm{iso}}(A)$.

From now on, since we will work in the isotropic setting, we will omit the term isotropic for such a points, when there is no possibility of confusion.   
\end{definition}

\begin{remark}
If $A \in \Psi_\mathrm{iso}^m(\mathbb{R}^n)$, $\mathrm{Ell}_\mathrm{iso}(A)$ is an open set and moreover for $A \in \Psi_\mathrm{iso}^m(\mathbb{R}^n)$ and $B \in \Psi_\mathrm{iso}^{m'}(\mathbb{R}^n)$, 
\[
\mathrm{Ell}_\mathrm{iso}(A\circ B)=\mathrm{Ell}_\mathrm{iso}(A)\cap \mathrm{Ell}_\mathrm{iso}(B).
\] 
\end{remark}

\begin{definition}\label{defessconesupp}
Let $a \in S_{\mathrm{iso}}^m(\mathbb{R}^n)$. We say that $X_0 \notin \mathrm{essconesupp}_\mathrm{iso} \, a\subseteq \dot{\mathbb{R}}^{2n}$, 
if there exists a conic neighborhood $\Gamma_{X_0}$ of $X_0$ such that, for all $\alpha\in \mathbb{N}_0^{2n}, \, N \in \mathbb{N}$ there exists $C_{\alpha,N}>0$ such that 
\begin{equation}\label{essconesupp}
\abs*{\partial_X^\alpha a(X)}\leq C_{\alpha,N} \braket{X}^{-N}, \quad \forall X \in \Gamma_{X_0}, \ \abs*{X}\geq 1. 
\end{equation}
The set $\mathrm{essconesupp}_\mathrm{iso} \, a$ is called \textit{the isotropic essential conic support of $a$}. 
\end{definition}

\begin{remark}
Note that, by definition, $\mathrm{essconesupp}_\mathrm{iso}\, a$ is a closed conic subset of $\dot{\mathbb{R}}^n$. 
\end{remark}

Moreover, note that this definition differs from that of conic support in \cite{HQ} that we recall next.

\begin{definition}
Let $a \in \mathscr{S}'(\mathbb{R}^{2n})$. \textit{The isotropic conic support of $a$}, denoted by $\mathrm{conesupp}_\mathrm{iso} \, a\subset\dot{\mathbb{R}}^{2n}$, 
is the set of all $X_0 \in \dot{\mathbb{R}}^{2n}$ 
such that, for all $\Gamma_{X_0}\subseteq \dot{\mathbb{R}}^{2n}$ conic neighborhood of $X_0$, the closure in $\mathbb{R}^{2n}$
\begin{center}
$\overline{\mathrm{supp} \, a \cap \Gamma_{X_0}}$ is not compact. 
\end{center}
\end{definition}

\begin{remark}
Note that if $a \in S_{\mathrm{iso}}^m(\mathbb{R}^n)$ then
\[
\mathrm{essconesupp}_\mathrm{iso} \, a \subseteq \mathrm{conesupp}_\mathrm{iso} \, a.
\]
In fact, if $X_0 \in \dot{\mathbb{R}}^{2n}$ is such that $\overline{\mathrm{supp} \, a \cap \Gamma_{X_0}}$ is compact, for some conic neighborhood $\Gamma_{X_0}$, then $a$ satisfies \eqref{essconesupp} there. 
\end{remark}

The notion of essential conic support of a symbol leads immediately to the definition of operator wave front set on this isotropic context. 

\begin{definition}
Let $A=\mathrm{Op}^\mathrm{w}(a)\in \Psi_{\mathrm{iso}}^m(\mathbb{R}^n)$. The \textit{isotropic wave front set} of $A$ is defined as 
\[
WF'(A):=\mathrm{essconesupp}_\mathrm{iso} \,a\subseteq \dot{\mathbb{R}}^{2n},
\]
which is a closed conic set of $\dot{\mathbb{R}}^{2n}$.
\end{definition}

Note first that the operator wave front set enjoys the following properties. 

\begin{proposition}\label{opWFprop}
Let $A\in \Psi_{\mathrm{iso}}^m(\mathbb{R}^n)$ and $B \in \Psi_{\mathrm{iso}}^{m'}(\mathbb{R}^n)$. Then\\ 
(i)
\[
WF'(A\circ B)\subseteq WF'(A) \cap WF'(B);
\]
(ii)
\[
WF'(A+B)\subseteq WF'(A) \cup WF'(B);
\]
(iii)
\[
WF'(A^\ast)=WF'(A).
\]
(iv)
\[
WF'(A)=\emptyset \quad \iff A \in \Psi^{-\infty}_{\mathrm{iso}}(\mathbb{R}^n).
\]
\end{proposition}

\begin{proof}
The properties follow directly from the formulas for the composition and the adjoint (Theorem \ref{compAB} and Theorem \ref{adjA*}) and from the definition of isotropic essential conic support.
\end{proof}

\begin{proposition}\label{usefulopWF}
Let $U\subseteq\dot{\mathbb{R}}^{2n}$ be an open cone and $K \subseteq U$ be a closed cone. Then, there exists $A=\mathrm{Op}^{\mathrm{w}}(a)\in \Psi^0_\mathrm{iso}(\mathbb{R}^n)$ such that 
$WF'(A)\subseteq U$ and $WF'(I-A)\cap K=\emptyset$. In particular, one actually has $\mathrm{supp}\,a\subset U$.
\end{proposition}
\begin{proof}
Let $\Gamma_1$ be a closed cone with $K\subset\Gamma_1\subset U.$ Take $\tilde{a}\in C^\infty(\dot{\mathbb{R}}^{2n})$ homogeneous of degree $0$
such that $\tilde{a}\bigl|_{\Gamma_1}\equiv 1$ and $\mathrm{supp}\,\tilde{a}\bigl|_{\mathbb{S}^{2n-1}}\subset U\cap\mathbb{S}^{2n-1}.$ Let then $\chi\in C^\infty(\mathbb{R}^{2n})$
be such that $0\leq\chi\leq 1$, $\chi\equiv 0$ if $|X|\leq 1/2$ and $\chi\equiv 1$ if $|X|\geq 1.$ Put $a=\chi\tilde{a}$ and $A=\mathrm{Op}^{\mathrm{w}}(a)$.
\end{proof}

To study the microsingularities of a distribution, one of the most important tools is the notion of microlocal parametrix. 

\begin{proposition}\label{micropmx}
Let $A \in \Psi_{\mathrm{iso}}^m(\mathbb{R}^n)$ and let $K \subseteq \dot{\mathbb{R}}^{2n}$ be closed such that
$\mathrm{pr}(K)\subsetneq\mathrm{Ell}_{\mathrm{iso}}(A)\cap\mathbb{S}^{2n-1}$, where $\mathrm{pr}\colon\dot{\mathbb{R}}^{2n}\longrightarrow\mathbb{S}^{2n-1}$ is
the projection. Then, there exists an operator 
$E \in \Psi_{\mathrm{iso}}^{-m}(\mathbb{R}^n)$, called microlocal parametrix of $A$ on $K$, such that  
\begin{equation}
K\cap WF'(AE-I)=\emptyset, \quad K \cap WF'(EA-I)=\emptyset.
\end{equation}
\end{proposition}

\begin{proof}
By Proposition \ref{usefulopWF}, applied when considering the closed cone $\{tX;\,\,X\in\mathrm{pr}(K),\,\,t>0\}\subsetneq\mathrm{Ell}_{\mathrm{iso}}(A)$, 
there exists $B_1\in \Psi_\mathrm{iso}^0(\mathbb{R}^n)$ such that 
\[
WF'(I-B_1)\cap K=\emptyset, \quad WF'(B_1)\subseteq \mathrm{Ell}_\mathrm{iso}(A)
\]
and by repeating this argument (since $WF'(B_1)$ is a closed conic set) there exists $B_2 \in \Psi_\mathrm{iso}^0(\mathbb{R}^n)$ such that 
\[
WF'(I-B_2)\cap WF'(B_1)=\emptyset, \quad WF'(B_2)\subseteq \mathrm{Ell}_\mathrm{iso}(A).
\]
Write $B_1=\mathrm{Op}^{\mathrm{w}}(b_1)$, $B_2=\mathrm{Op}^{\mathrm{w}}(b_2)$ and $A=\mathrm{Op}^{\mathrm{w}}(a)$. From Proposition \ref{usefulopWF} 
we have $\mathrm{supp}\,b_2\subsetneq\mathrm{Ell}_{\mathrm{iso}}(A)$.
We now choose closed cones $\Gamma_1\subset\Gamma_2\subset\mathrm{Ell}_{\mathrm{iso}}(A)$ such that $K\subset\Gamma_1.$ Then there exist $c,C>0$ such that
$|a(X)|\geq c|X|^m$ for all $X\in\Gamma_2$ with $|X|\geq C.$ Choose $\chi\in C^\infty(\mathbb{R}^{2n})$ with $0\leq\chi\leq 1$ such that $\mathrm{supp}\,\chi\subset(\Gamma_2
\cap\{|X|\geq 2C\})$ and $\chi\bigl|_{\Gamma_1\cap\{|X|\geq 3C\}}\equiv 1$. Put  

\[
e_0:=b_2\chi/a, \quad E_0 :=\mathrm{Op}^\mathrm{w}(e_0)\in \Psi_\mathrm{iso}^{-m}(\mathbb{R}^n),\quad B_2^\chi=\mathrm{Op}^{\mathrm{w}}(b_2\chi)\in\Psi^0_{\mathrm{iso}}(\mathbb{R}^n),
\]
and write, with $R=I-B_2^\chi+\Psi^{-1}_{\mathrm{iso}}\in\Psi^0_{\mathrm{iso}}$, 
\begin{equation}\label{eqpara}
AE_0=I-R=I-B_1R-(I-B_1)R.
\end{equation}
Next note that $WF'((I-B_1)R) \cap K=\emptyset$, and since $B_1R=\mathrm{Op}^{\mathrm{w}}(b_1\sharp r)$, we have 
$$b_1\sharp r=b_1(1-b_2\chi)+b_1\sharp S^{-1}_{\mathrm{iso}}+S^{-1}_{\mathrm{iso}}\sharp S^{-1}_{\mathrm{iso}}=0+S^{-1}_{\mathrm{iso}}.$$
Therefore $B_1R \in \Psi_\mathrm{iso}^{-1}(\mathbb{R}^n)$.

\noindent At this point, we use a Neumann series argument to invert microlocally $I-B_1R$. We choose
\[
E' \sim \sum_{j=1}^\infty (B_1R)^j \in \Psi^{-1}(\mathbb{R}^n).
\]
Since $(I-B_1R)(I+\sum_{j=1}^N(B_1R)^j)=I-(B_1R)^{N+1}$, for all $N\in\mathbb{N}$, we get
\begin{equation}\label{eqpara2}
(I-B_1R)(I+E')=I+S, \quad S \in \Psi_\mathrm{iso}^{-\infty}(\mathbb{R}^n).
\end{equation}
Hence, we put 
\[
E:=E_0(I+E') \in \Psi_\mathrm{iso}^{-m}(\mathbb{R}^n).
\]
In conclusion, by \eqref{eqpara} and \eqref{eqpara2}
\[
\begin{split}
AE&=AE_0(I+E')\\
&=\Bigr(I-B_1R-(I-B_1)R\Bigl)(I+E')\\
&=I+S-((I-B_1)R)(I+E'),
\end{split}
\]
with $WF'(((I-B_1)R)(I+E')-S)\cap K=\emptyset$. We have thus proved that there exists $E \in \Psi_\mathrm{iso}^{-m}(\mathbb{R}^n)$ 
such that $WF'(AE-I)\cap K=\emptyset$. Similarly, we construct $E_1 \in \Psi_{\mathrm{iso}}^{-m}(\mathbb{R}^n)$ such that $WF'(E_1A-I)\cap K=\emptyset$ and
\[
E=(E_1A+R_1)E=E_1(AE)+R_1'=E_1+R_1'+R_2,
\] 
with $R_1,R_1',R_2$ having operator wave front sets all disjoint from $K$. This proves that any right microlocal parametrix is also a left microlocal parametrix and viceversa. 
\end{proof}

In order to study the isotropic regularity we recall the definition of Shubin-Sobolev space and to do that we refer to \cite{H}, \cite{P1}.

Since $\Lambda^s$ is globally elliptic, we may find $E_s \in \Psi_{\mathrm{iso}}^{-s}(\mathbb{R}^n)$ and $R \in \Psi_{\mathrm{iso}}^{-\infty}(\mathbb{R}^n)$ such that $E_s\Lambda^s=I+R$. Hence, we define the Shubin-Sobolev spaces as follows.

\begin{definition}\label{defBs}
Let $s \in \mathbb{R}$ and let $p \in \mathbb{N}$, $p\geq s$. We define the \textit{Shubin-Sobolev space of order $s \in \mathbb{R}$} as 
\[
B^s(\mathbb{R}^n)=\lbrace u \in \mathscr{S}'(\mathbb{R}^n); \ \Lambda^s u \in L^2(\mathbb{R}^n) \rbrace,
\]
endowed with the inner product
\[
(u,v)_{B^s}:=(\Lambda^s u,\Lambda^s v)_0 + \sum_{\abs{\alpha}+\abs{\beta}\leq p}(x^\alpha D_x^\beta Ru,x^\alpha D_x^\beta Rv)_0, \quad u,v \in B^s(\mathbb{R}^n),
\]
and norm 
\[
\norma*{u}_{B^s}=(u,u)_{B^s}^{1/2}, \quad u \in B^s(\mathbb{R}^n).
\]
\end{definition}
It can be shown that a different choice of $p$ gives rise to an equivalent norm.

These spaces have the property that 
$$\mathscr{S}(\mathbb{R}^n)\hookrightarrow B^s(\mathbb{R}^n)\quad\text{\rm with dense image,}$$
and
\[
B^{-s}(\mathbb{R}^n)=(B^s(\mathbb{R}^n))^\ast,\,\,\,s\in\mathbb{R}.
\]
Moreover, for $s \in \mathbb{N}_0$,
\[
B^s(\mathbb{R}^n)=\lbrace u \in L^2(\mathbb{R}^n); \quad \sum_{\abs*{\alpha}+\abs*{\beta} \leq s}\norma*{x^\alpha \partial_x^\beta u}^2_0 <+\infty \rbrace
\]
and
\[
\norma*{u}_{B^s}^2 \approx \sum_{\abs*{\alpha}+\abs*{\beta} \leq s}\norma*{x^\alpha \partial_x^\beta u}^2_0,
\]
that is, the two norms are equivalent. (We say that $A\approx B$, for $A,B>0$ if there exists $C>0$ such that $C^{-1}A\leq B\leq CA$.)
 
\begin{theorem}\label{boundtheo}
Let $A \in \Psi_\mathrm{iso}^m(\mathbb{R}^n)$. Then, for all $s \in \mathbb{R}$,  
\[
A:B^s(\mathbb{R}^n)\rightarrow B^{s-m}(\mathbb{R}^n)
\]
is bounded.
\end{theorem}

Using those spaces we want to introduce an isotropic wave front set as follows (see \cite{DGW}). 

\begin{definition}\label{def.WFiso}
Let $u \in \mathscr{S}'(\mathbb{R}^n)$. \textit{The Isotropic Wave Front Set} of $u \in \mathscr{S}'(\mathbb{R}^n)$ is defined as
\[
WF_{\mathrm{iso}}(u) = 
\bigcap_{\substack{A \in \Psi^0_{\mathrm{iso}}(\mathbb{R}^n) \\ Au \in \mathscr{S}(\mathbb{R}^n)}}
\mathrm{Char}_{\mathrm{iso}}(A)
\]
and \textit{the Isotropic Wave Front Set of $u$ of order $s \in \mathbb{R}$} is defined as
\[
WF^s_{\mathrm{iso}}(u) = 
\bigcap_{\substack{A \in \Psi^0_{\mathrm{iso}}(\mathbb{R}^n) \\ Au \in B^s(\mathbb{R}^n)}}
\mathrm{Char}_{\mathrm{iso}}(A).
\]
\end{definition}

\begin{remark}
Note that if $s'<s$, since $B^{s}(\mathbb{R}^n) \subseteq B^{s'}(\mathbb{R}^n)$, one has $WF_{\mathrm{iso}}^{s'} (u)\subseteq WF_{\mathrm{iso}}^s(u)$.
\end{remark}

\begin{remark}\label{sminusl}
By multiplying by an elliptic operator that does not change the characteristic set, one easily sees that the isotropic wave front set of order $s\in \mathbb{R}$ can be equivalently defined as 
\[
WF_{\mathrm{iso}}^s (u)= \bigcap_{\substack{A \in \Psi^s_{\mathrm{iso}}(\mathbb{R}^n) \\ Au \in L^2(\mathbb{R}^n)}}\mathrm{Char}_{\mathrm{iso}}(A),
\]
ore more generally, for $\ell \in \mathbb{R}$,
\[
WF_{\mathrm{iso}}^s (u)= \bigcap_{\substack{A \in \Psi^{s-\ell}_{\mathrm{iso}}(\mathbb{R}^n) \\ Au \in B^\ell(\mathbb{R}^n)}}\mathrm{Char}_{\mathrm{iso}}(A).
\]
\end{remark}

\begin{remark}
Note that the isotropic wave front set
of $u \in \mathscr{S}'(\mathbb{R}^n)$
can be rephrased
by saying that $X_0 \notin WF_{\mathrm{iso}}(u)$
if there exists $A=\mathrm{Op}^\mathrm{w}(a)$,
with $a \in S^0_{\mathrm{iso}}(\mathbb{R}^n)$, 
such that $X_0 \notin \mathrm{Char}(A)$ and $Au \in \mathscr{S}(\mathbb{R}^n)$. \\
Note that $A$ may also be written as $A=\mathrm{Op}(a_L)$, where $a_L$ is the left-symbol computed from $a$ which has the same principal part.
Hence, considering $A$ written as a Weyl-quantizazion or a left-quantization gives the same property (as a matter of fact, the characteristic set
depends only on the principal part of the symbol of $A$ and hence it does not see the kind of quantization we are using).
  
\noindent The same obviously holds for the $s$-isotropic
wave front set $WF^s_{\mathrm{iso}}(u)$.
\end{remark}

\begin{proposition}\label{propregu}
Let $u \in \mathscr{S}'(\mathbb{R}^n)$. Then for $s \in \mathbb{R}$,
\begin{equation}\label{regs}
WF_\mathrm{iso}^s(u)=\emptyset \iff u \in B^s(\mathbb{R}^n)
\end{equation}
and
\begin{equation}\label{reg}
WF_{\mathrm{iso}}(u)= \emptyset \iff u \in \mathscr{S}(\mathbb{R}^n).
\end{equation}
\end{proposition}

\begin{proof}
Let us prove \eqref{regs}. The proof of \eqref{reg} is completely similar.

\noindent If $u \in B^s(\mathbb{R}^n)$, for every $X_0 \in \dot{\mathbb{R}}^{2n}$ we may find an operator $A \in \Psi^0(\mathbb{R}^n)$ such that $X_0 \notin \mathrm{Char}_\mathrm{iso}(A)$ and $Au \in B^s(\mathbb{R}^n)$, by Theorem \ref{boundtheo}.
Then $X_0 \notin WF_\mathrm{iso}^s(u)$. 

\noindent Conversely, if $WF_\mathrm{iso}^s(u)=\emptyset$, for every $X\in \mathbb{S}^{2n-1}$ there exists $A_{X} \in \Psi^0_\mathrm{iso}(\mathbb{R}^n)$ and $\Gamma_{X}$ open cone such that $A_{X}$ is elliptic on $\Gamma_{X}$ and $A_Xu \in B^s(\mathbb{R}^n)$. By compactness of $\mathbb{S}^{2n-1}$, there exists $X_1,\dots,X_N \in \mathbb{S}^{2n-1}$ and open cones $\Gamma_{X_j}$ 
(conic neighborhoods of the $X_j$) such that
 \[
\dot{\mathbb{R}}^{2n}=\bigcup_{j=1}^N \Gamma_{X_j}.
\]
Then, if we set 
\[
A:=\sum_{j=1}^N A_{X_j}^\ast A_{X_j} \in \Psi_\mathrm{iso}^0(\mathbb{R}^n),
\]
we have that $A$ is an elliptic operator on $\dot{\mathbb{R}}^{2n}$ such that $Au \in B^s(\mathbb{R}^n)$. Hence, we may find a global parametrix 
$E\in \Psi_\mathrm{iso}^0(\mathbb{R}^n)$ for $A$. Therefore
\[
u=EAu+Ru,
\]
where $R \in \Psi_\mathrm{iso}^{-\infty}(\mathbb{R}^n)$, so that we can conclude that $u \in B^s(\mathbb{R}^n)$. 
\end{proof}

We now study, for a $\psi$do $A$ and a temperate distribution $u$, the relation among the isotropic microsingularities of $Au$,
those of $u$ and those of the operator wave front set of $A$. 
We first prove that the isotropic pseudodifferential operators are \textit{microlocal}. 

\begin{theorem}\label{isomicro}
Let $A \in \Psi_\mathrm{iso}^m(\mathbb{R}^n)$ and $u \in \mathscr{S}'(\mathbb{R}^n)$. Then, for all $s \in \mathbb{R}$
\begin{equation}\label{microlocals}
WF_\mathrm{iso}^{s-m}(Au) \subseteq WF'(A)\cap WF_\mathrm{iso}^s(u)
\end{equation}
and 
\begin{equation}\label{microlocal}
WF_\mathrm{iso}(Au) \subseteq WF'(A)\cap WF_\mathrm{iso}(u).
\end{equation}
\end{theorem}
\begin{proof}
We prove \eqref{microlocals}. Let us consider first $X_0\in\dot{\mathbb{R}}^{2n}$ such that $X_0 \notin WF'(A)$. Thanks to Proposition \ref{usefulopWF} 
we may construct $B \in \Psi_\mathrm{iso}^0(\mathbb{R}^n)$ such that $B$ is elliptic at $X_0$ and $WF'(A)\cap WF'(B)=\emptyset$. 
Therefore, by Proposition \ref{opWFprop}$-(i)$, $WF'(B\circ A)=\emptyset$ and so by Proposition \ref{opWFprop}$-(iv)$
\[
B(Au)=(BA)u \in \mathscr{S}(\mathbb{R}^n)\subseteq B^{s-m}(\mathbb{R}^n).
\]
Thus, since $X_0 \notin \mathrm{Char}_\mathrm{iso}(B)$ we have $X_0 \notin WF_{\mathrm{iso}}^{s-m}(Au)$. 

\noindent We consider now $X_0 \notin WF^s_\mathrm{iso}(u)$. By definition there exists $C \in \Psi_\mathrm{iso}^0(\mathbb{R}^n)$ elliptic at $X_0$ such that 
$Cu \in B^s(\mathbb{R}^n)$ and again by Proposition \ref{usefulopWF}, we can construct $B \in \Psi^0_\mathrm{iso}(\mathbb{R}^n)$, elliptic at $X_0$, 
such that $WF'(B)\subseteq \mathrm{Ell}_\mathrm{iso}(C)$. By Proposition \ref{micropmx} there exists $E$ microlocal 
parametrix of $C$ on $WF'(B)$, that is, setting $R=EC-I \in \Psi_\mathrm{iso}^0(\mathbb{R}^n)$ we have
\begin{equation}\label{eqRB}
WF'(R)\cap WF'(B)=\emptyset.
\end{equation}
Then,
\[
B(Au)=BA(EC-R)u=BAE(Cu)-BARu.
\]
Now, since $Cu \in B^s(\mathbb{R}^n)$, $A \in \Psi_\mathrm{iso}^m(\mathbb{R}^n)$, $B,E\in \Psi_\mathrm{iso}^0(\mathbb{R}^n)$, by Theorem \ref{boundtheo} we have that $BAE(Cu) \in B^{s-m}(\mathbb{R}^n)$ and by \eqref{eqRB} and Proposition \ref{opWFprop}-$(i)$, $BARu \in \mathscr{S}(\mathbb{R}^n)\subseteq B^{s-m}(\mathbb{R}^n)$. In conclusion 
\[
B(Au) \in B^{s-m}(\mathbb{R}^n)
\]
and $X_0 \notin \mathrm{Char}_{\mathrm{iso}}(B)$. Thus, also in this case $X_0 \notin WF^{s-m}_\mathrm{iso}(Au)$.

\noindent The proof of \eqref{microlocal} is analogous, requiring Schwartz regularity instead of $B^s$-regularity.
\end{proof}

\begin{corollary}\label{corAu}
Let $s,m \in \mathbb{R}$, $u \in \mathscr{S}'(\mathbb{R}^n)$ and $A \in \Psi_\mathrm{iso}^m(\mathbb{R}^n)$ such that $WF'(A)\cap WF_\mathrm{iso}^s(u)=\emptyset$. Then $Au \in B^{s-m}(\mathbb{R}^n)$. 
\end{corollary}
\begin{proof}
By Theorem \ref{isomicro} we have $WF_\mathrm{iso}^{s-m}(Au)=\emptyset$ and so, by Proposition \ref{propregu}, $Au \in B^{s-m}(\mathbb{R}^n)$.
\end{proof}

We next prove the \textit{microlocal elliptic regularity} in the present isotropic context. 

\begin{theorem}
Let $A\in \Psi_\mathrm{iso}^m(\mathbb{R}^n)$ and $u \in \mathscr{S}'(\mathbb{R}^n)$. Then for all $s \in \mathbb{R}$
\begin{equation}\label{microells}
WF_\mathrm{iso}^s(u)\subseteq WF_{\mathrm{iso}}^{s-m}(Au) \cup \mathrm{Char}_\mathrm{iso}(A),
\end{equation}
and 
\begin{equation}\label{microell}
WF_\mathrm{iso}(u)\subseteq WF_\mathrm{iso}(Au) \cup \mathrm{Char}_\mathrm{iso}(A).
\end{equation}
In particular, if $A$ is elliptic 
\begin{equation}\label{ellpres}
WF_\mathrm{iso}^s(u)=WF^{s-m}_\mathrm{iso}(Au) \quad \mathit{and} \quad WF_\mathrm{iso}(u)=WF_\mathrm{iso}(Au).
\end{equation}
\end{theorem}
\begin{proof}
Suppose $X_0 \notin WF_\mathrm{iso}^{s-m}(Au)$ and $X_0 \in \mathrm{Ell}_\mathrm{iso}(A)$. Since $X_0 \notin WF_\mathrm{iso}^{s-m}(Au)$ by definition there exists 
$B \in \Psi_\mathrm{iso}^{s-m}(\mathbb{R}^n)$ such that $X_0 \in \mathrm{Ell}_\mathrm{iso}(B)$ and $BAu \in L^2(\mathbb{R}^n)$. In conclusion 
$B\circ A \in \Psi_\mathrm{iso}^s(\mathbb{R}^n)$, $X_0 \in \mathrm{Ell}_\mathrm{iso}(B\circ A)=\mathrm{Ell}_\mathrm{iso}(A)\cap \mathrm{Ell}_\mathrm{iso}(B)$ and 
$BA u \in L^2(\mathbb{R}^n)$. Hence $X_0 \notin WF_\mathrm{iso}^{s}(u)$. Once again the proof of \eqref{microell} is the same.
Finally, \eqref{ellpres} follows from \eqref{microells}, \eqref{microell} and Theorem \ref{isomicro}.
\end{proof}

We conclude this section by showing that, as for the usual wave front set, the following result holds.

\begin{proposition}\label{WF&WFs}
Let $u \in \mathscr{S}'(\mathbb{R}^n)$. Then 
\[
WF_{\mathrm{iso}}(u)=\overline{\bigcup_{s \in \mathbb{R}}WF_{\mathrm{iso}}^s(u)}.
\]
\end{proposition} 
\begin{proof}
The inclusion '$\supseteq$' follows from the fact that $WF_\mathrm{iso}^s(u) \subseteq WF_\mathrm{iso}(u)$ for all $s \in \mathbb{R}$ combined with the fact that $WF_\mathrm{iso}(u)$ is a closed set.

\noindent We now prove '$\subseteq$'. If $X_0 \notin \overline{\bigcup_{s \in \mathbb{R}}WF_{\mathrm{iso}}^s(u)}$, there exists an 
open conic set $\Gamma_{X_0}$ containing $X_0$, such that 
\[
\Gamma_{X_0} \cap\overline{\bigcup_{s \in \mathbb{R}}WF_{\mathrm{iso}}^s(u)}=\emptyset.
\]
Let now $a \in S^0_\mathrm{iso}(\mathbb{R}^n)$ such that $\mathrm{supp} \, a \subseteq \Gamma_{X_0}$ and $X_0 \in \mathrm{Ell}_\mathrm{iso}(A)$, with $A=\mathrm{Op}^\mathrm{w}(a)$. 
By Theorem \ref{isomicro} and Theorem \ref{propregu}, since $WF'(A)\cap WF_\mathrm{iso}^s(u)=\emptyset$, for all $s \in \mathbb{R}$, we have  
\[
Au \in \bigcap_{s \in \mathbb{R}}B^s(\mathbb{R}^n)=\mathscr{S}(\mathbb{R}^n).
\] 
Therefore $X_0 \notin WF_{\mathrm{iso}}(u)$. 
\end{proof}
\begin{remark}
    Note that Proposition \ref{WF&WFs} provides a stratification of  $WF_\mathrm{iso}$ in terms of $WF^s_\mathrm{iso}$. 
    This, in particular, leads to the fact that the results for $WF_{\mathrm{iso}}$ can be described starting from the refined version given for $WF_{\mathrm{iso}}^s$. 
For instance, \textit{the microlocal property} \eqref{microlocal} and \textit{the microlocal elliptic regularity} \eqref{microell} follow from \eqref{microlocals} and \eqref{microells}, respectively. 
\end{remark}

\section{Relation with the time-frequency calculus} \label{s:TFCalc}
We next wish to establish a relation between this isotropic pseudodifferential calculus and the \textit{time-frequency} calculus developed by Gr\"ochenig \cite{G} 
(see also \cite{SW}), of which we briefly recall the basic notions below. 

\begin{definition}\label{def.STFT}
Let $u \in \mathscr{S}'(\mathbb{R}^n)$ and $\psi \in \mathscr{S}(\mathbb{R}^n)$. 
We define \textit{the short-time Fourier transform (STFT)} of $u$ with respect to the window function $\psi$, as the function
\[
\mathbb{R}^{2n} \ni (x,\xi)=X \mapsto V_\psi u (X)=(u, \Pi(X)\psi) \in \mathbb{C},
\] 
where $\Pi(X)=M_\xi T_x$ is the time-frequency shift, composition of the translation operator $T_x \psi(\cdot)=\psi(\cdot-x)$ 
and the modulation operator $M_\xi\psi(\cdot)=e^{i\langle \cdot, \xi\rangle}\psi(\cdot)$. 
\end{definition}

\begin{remark}
For the Short-Time Fourier Transform of $u\in \mathscr{S}'(\mathbb{R}^n)$ we have that
$V_{\psi}u \in C^\infty(\mathbb{R}^{2n})$ and there exists $\mathrm{N} \in \mathbb{N}_0$ such that 
\[
\abs{V_{\psi}u(X)} \lesssim \braket{X}^N, \ X \in \mathbb{R}^{2n}.
\]
\end{remark}

Moreover, if $u$ and $\psi$ are in $L^2(\mathbb{R}^n)$, similarly to the Fourier transform we have the following orthogonality relation (cf. \cite{G}, Theorem 3.2.1 and Corollary 3.2.2).

\begin{theorem}
Let $u_1,u_2,\psi_1,\psi_2 \in L^2(\mathbb{R}^n)$, then $V_{\psi_j}u_j \in L^2(\mathbb{R}^{2n})$ for $j=1,2,$ and 
\[
( V_{\psi_1}u_1,V_{\psi_2}u_2 )_{L^2(\mathbb{R}^{2n})}=(2\pi)^n( u_1,u_2 )_{L^2(\mathbb{R}^n)} \overline{( \psi_1,\psi_2 )}_{L^2(\mathbb{R}^n)}.
\] 
In particular, if $u,\psi \in L^2(\mathbb{R}^n)$, then 
\begin{equation}
\norma*{V_{\psi}u}_{L^2(\mathbb{R}^{2n})} = (2\pi)^{n/2}\norma*{u}_{L^2(\mathbb{R}^n)} \norma*{\psi}_{L^2(\mathbb{R}^n)}.
\label{eqV-isometry}\end{equation}

So, if $\norma*{\psi}_{L^2(\mathbb{R}^n)}=1$,

\[
\norma*{V_{\psi}u}_{L^2(\mathbb{R}^{2n})} = (2\pi)^{n/2}\norma*{u}_{L^2(\mathbb{R}^n)}, \quad \mathrm{for} \ \mathrm{all} \ u \in L^2(\mathbb{R}^n),
\]
thus in this case the STFT is a conformal isometry from $L^2(\mathbb{R}^n)$ into $L^2(\mathbb{R}^{2n})$.
\end{theorem}

\begin{remark}
\label{rmkInvSTFT}
Note that if we consider $F \in L^2(\mathbb{R}^{2n})$, then the conjugate-linear functional
\[
T^F_\psi h= \int_{\mathbb{R}^{2n}}F(X)\overline{(h,\Pi(X)\psi)}_{L^2(\mathbb{R}^n)} \ dX
\]
is a bounded functional on $L^2(\mathbb{R}^n)$.

\noindent In fact, by Cauchy-Schwartz
\[
\abs*{T^F_\psi h}\leq \norma*{F}_{L^2(\mathbb{R}^{2n})}\norma*{V_\psi h}_{L^2(\mathbb{R}^{2n})}=
(2\pi)^{n/2}\norma*{F}_{L^2(\mathbb{R}^{2n})}\norma*{\psi}_{L^2(\mathbb{R}^n)}\norma*{h}_{L^2(\mathbb{R}^n)}.
\]
Therefore, by the Riesz Representation Theorem, there exists $f \in L^2(\mathbb{R}^n)$ such that
\[
f=\int_{\mathbb{R}^{2n}}F(X)\Pi(X)\psi \ dX.
\]
\end{remark}

This result leads to the following inversion formula.

\begin{corollary}
Suppose that $\psi, \gamma \in L^2(\mathbb{R}^n)$ and $(\psi,\gamma)_{L^2(\mathbb{R}^n)}\neq 0$. Then for all $u \in L^2(\mathbb{R}^n)$ \[
u=\frac{1}{(\psi,\gamma)_{L^2(\mathbb{R}^n)}}\int_{\mathbb{R}^{2n}}V_{\psi}u(X)\Pi(X)\gamma \ dX.
\]
\end{corollary}

Let now $T_\psi$ be the linear operator defined by
\[
T_\psi F=\int_{\mathbb{R}^{2n}}F(X)\Pi(X)\psi \ dX, \quad F \in L^2(\mathbb{R}^{2n}).
\]
By Remark \ref{rmkInvSTFT}, $T_\psi$ is a bounded
operator from $L^2(\mathbb{R}^{2n})$ to $L^2(\mathbb{R}^n)$, and for $F \in L^2(\mathbb{R}^{2n}), \ h \in L^2(\mathbb{R}^n)$ we have
\[
\begin{split}
(T_\psi F,h)_{L^2(\mathbb{R}^n)}&= \int_{\mathbb{R}^{2n}}F(X)(\Pi(X)\psi,h)_{L^2(\mathbb{R}^n)} \ dX \\
&=(F,V_{\psi}h)_{L^2(\mathbb{R}^{2n})}=(V_\psi^\ast F,h)_{L^2(\mathbb{R}^{2n})}. 
\end{split}
\]
Thus, $V_\psi^\ast F=T_\psi F$ and when $|\!|\psi|\!|_{L^2}=1$ the inversion formula is given by
\[
V_{\psi}^\ast V_\psi=(2\pi)^nI.
\]
The main tool that is used to reach the correspondence between this time-frequency calculus and the isotropic calculus, 
is the so called Weyl-Wick lemma, and in order to introduce this tool we introduce a new class of operators, also known as localization operators, that are related to pseudodifferential operators.

\begin{definition}
For $a \in \mathscr{S}'(\mathbb{R}^{2n})$, we define the corresponding localization operator called \textit{anti-Wick quantized operators} $A_a$, 
weakly by its action on $f,g \in \mathscr{S}(\mathbb{R}^n)$, that is,
\begin{equation}
(A_af,g)_0=(2\pi)^{-n}(aV_{\psi}f,V_{\psi}g)_0=(2\pi)^{-n}(V_{\psi}^*aV_{\psi}f,g)_0,
\label{eqAa}\end{equation}
that is, $A_a=(2\pi)^{-n}V_{\psi}^*aV_{\psi}$.
\end{definition}

This operator can be written as a Weyl-pseudodifferential operator using the following result known as Weyl-Wick lemma (cf. \cite{NR} Proposition 1.7.9, Theorem 1.7.10). 

\begin{lemma}\label{Weylwicktheorem}
Let $a \in \mathscr{S}'(\mathbb{R}^{2n})$. There exists $b \in \mathscr{S}'(\mathbb{R}^{2n})$ such that $A_a=\mathrm{Op}^{\mathrm{w}}(b)$ where 
\[
b=\pi^{-n}e^{-\abs*{\cdot}^2} \ast a.
\]
If $a \in S_{\mathrm{iso}}^m(\mathbb{R}^{n})$ then $b \in S_{\mathrm{iso}}^m(\mathbb{R}^{n})$, where
\[
b \sim \sum_{\alpha \in \mathbb{N}_0^{2n}}c_\alpha \partial^\alpha a,
\]
with $c_0=1$. Moreover, denoting $A=\mathrm{Op}^{\mathrm{w}}(a)$ and $B=\mathrm{Op}^{\mathrm{w}}(b)$ we have 
$\mathrm{Char}_{\mathrm{iso}}(B) =\mathrm{Char}_{\mathrm{iso}}(A)$ and $\mathrm{essconesupp}_\mathrm{iso} \, b \subseteq \mathrm{essconesupp}_\mathrm{iso} \, a$.
\end{lemma}

Now, we want to use the Short-Time Fourier Transform to define the Gabor wave front set.

\begin{definition}
\label{isoSTFT}
We say that 
$0 \neq X_0 \not\in WF_{G}(u)$ if there exists an open cone $\Gamma_{X_0} \subseteq \dot{\mathbb{R}}^{2n}$ that contains
 $X_0$ such that 
\[
\sup_{X \in \Gamma_{X_0}} \braket{X}^N \abs*{V_\psi u (X)} <+\infty, \quad \forall N\geq 0.
\]
The set $WF_{G}(u)\subseteq\dot{\mathbb{R}}^{2n}$ is called the \textit{Gabor wave front set}. 
\end{definition}

\begin{remark}
This definition is independent of the choice of window 
function $\psi$ (cf. \cite{RW}, Corollary 3.3).
\end{remark}

Now, note that it is possible to characterize the Shubin-Sobolev wave front sets in terms of the 
Short-Time Fourier Transform as follows. 

\begin{remark}
Note that by the properties of elliptic operator the Shubin-Sobolev spaces of order $s \in \mathbb{R}$ 
previously defined can be equivalently defined as 
\[
B^s(\mathbb{R}^n)=\lbrace u \in \mathscr{S}'(\mathbb{R}^n); \ Au \in L^2(\mathbb{R}^n) \rbrace,
\]
for any given elliptic operator $A \in \Psi_{\mathrm{iso}}^s(\mathbb{R}^{n})$.
\end{remark}

One may give an equivalent definition of Shubin-Sobolev spaces in terms of localization operator and modulation spaces. 

\begin{proposition}
Let $s \in \mathbb{R}$ and $\psi\in\mathscr{S}(\mathbb{R}^n)$. All the following three spaces coincide with $B^s(\mathbb{R}^n)$: \\
$(i)$ 
\[
\begin{split}
B_A^s(\mathbb{R}^n):=&\lbrace u \in \mathscr{S}'(\mathbb{R}^n); \ Au \in L^2(\mathbb{R}^n)\rbrace,  \quad 
A \in \Psi_{\mathrm{iso}}^s(\mathbb{R}^{n}) \ \mathit{elliptic},\\
& |\!|u|\!|_{B^s_A}:=|\!|Au|\!|_0+\sum_{\abs{\alpha}+\abs{\beta}\leq p} \norma*{x^\alpha D_x^\beta R' u}_0;
\end{split}
\]
with $p \in \mathbb{N}$, $p \geq s$ and $E'A=I+R'$ with $E' \in \Psi_{\mathrm{iso}}^{-s}(\mathbb{R}^n)$ and $R' \in \Psi_{\mathrm{iso}}^{-\infty}(\mathbb{R}^n)$.\\
$(ii)$ 
\[
B^s_{A_a}(\mathbb{R}^n):=\lbrace u \in \mathscr{S}'(\mathbb{R}^n); \ A_{\braket{\cdot}^s}u \in L^2(\mathbb{R}^{2n}) \rbrace,
\quad |\!|u|\!|_{B^s_{A_a}}:=|\!|A_{\langle\cdot\rangle^s}u|\!|_0;
\]
$(iii)$
\[
B^s_{ST}(\mathbb{R}^n):=\lbrace u \in \mathscr{S}'(\mathbb{R}^n); \ \braket{\cdot}^sV_\psi u \in L^2(\mathbb{R}^{2n}) \rbrace,
\quad |\!|u|\!|_{B^s_{ST}}=|\!|\langle\cdot\rangle^sV_\psi u|\!|_0. 
\]
\end{proposition} 
\begin{proof}
As an overall observation, we notice that the duals of the above spaces with order $s$ are the same with order $-s$.\\
Next, we prove $iii)\Leftrightarrow ii)$. Let first $u \in B^s_{ST}(\mathbb{R}^n)$. Since
$$|\!|V^*_\psi F|\!|_{L^2(\mathbb{R}^n_y)}\leq|\!|\psi|\!|_{L^2(\mathbb{R}^n_y)}|\!|F|\!|_{L^2(\mathbb{R}^{2n}_{x,\xi})},\quad\forall F\in  L^2(\mathbb{R}^{2n}_{x,\xi}),$$
when $\braket{\cdot}^sV_\psi u \in L^2(\mathbb{R}^{2n})$ then by (\ref{eqAa})  we immediately have 
\[
\norma*{A_{\braket{\cdot}^s}u}_{L^2}\leq(2\pi)^{-n/2}\norma*{\psi}_{L^2}\norma*{\braket{\cdot}^sV_\psi u}_{L^2}\leq C \norma*{\braket{\cdot}^sV_\psi u}_{L^2}, \quad \forall s \in \mathbb{R}.
\]
This implies 
\[
B_{ST}^s (\mathbb{R}^n)\subseteq B_{A_a}^s(\mathbb{R}^n), \quad \forall s \in \mathbb{R},
\]
and, by duality
\[
B_{A_a}^s(\mathbb{R}^n)=(B_{A_a}^{-s}(\mathbb{R}^n))'\subseteq (B_{ST}^{-s}(\mathbb{R}^n))'=B_{ST}^s(\mathbb{R}^n), \quad \forall s \in \mathbb{R}.
\]
We next prove $ii) \Leftrightarrow i)$. Let $ u \in B_{A_a}^s(\mathbb{R}^n)$. By Lemma \ref{Weylwicktheorem} there exists $b \in S_{\mathrm{iso}}^s(\mathbb{R}^n)$ 
such that $B:=\mathrm{Op}^{\mathrm{w}}(b)=A_{\braket{\cdot}^s}$ and $\mathrm{Char}_{\mathrm{iso}}(B)=\mathrm{Char}_{\mathrm{iso}}(\mathrm{Op}^{\mathrm{w}}(\braket{\cdot}^s))$. 
Then $B$ is elliptic and by hypothesis
\[
B u=A_{\braket{\cdot}^s}u \in L^2(\mathbb{R}^n),
\]
so that 
\[
B_{A_a}^s(\mathbb{R}^n) \subseteq B^s(\mathbb{R}^n), \quad \forall s \in \mathbb{R}.
\] 
On the other hand, if $u \in B^s(\mathbb{R}^n)$, we define $a=\braket{\cdot}^s$ and 
\[
b=\pi^{-n}e^{-\abs*{\cdot}^2} \ast \braket{\cdot}^s.
\]
Again by the Weyl-Wick correspondence $B=\mathrm{Op}^{\mathrm{w}}(b)$ is elliptic, hence by hypothesis 
\[
Bu \in L^2(\mathbb{R}^n).
\]
By construction $\mathrm{Op}^{\mathrm{w}}(b)=A_{\braket{\cdot}^s}$, whence we also have
\[
B^s(\mathbb{R}^n) \subseteq B_{A_a}^s(\mathbb{R}^n) \quad \forall s \in \mathbb{R}.
\]
\end{proof}

We are in a position to give the definition of Sobolev-Gabor wave front set.

\begin{definition}
If $u \in \mathscr{S}'(\mathbb{R}^n)$ and $\psi \in \mathscr{S}(\mathbb{R}^n)\setminus \lbrace 0 \rbrace$ then 
the \textit{Gabor wave front set $WF_{G}^s(u)$ of order $ s \in \mathbb{R}$} is defined by saying that 
$X_0 \not\in WF_{G}^s(u)$ if there exists an open cone $\Gamma=\Gamma_{X_0} \subseteq \dot{\mathbb{R}}^{2n}$ containing 
$X_0$ such that 
\[
\chi_\Gamma\braket{\cdot}^{s}V_{\psi}u \in L^2(\mathbb{R}^{2n}),
\]
where here $\chi_\Gamma$ denotes the characteristic function of the cone $\Gamma$.
\end{definition}

Note that for $s'<s$ we have $WF_G^{s'}(u) \subseteq WF_G^{s}(u)$ and also in this case the definition of $WF_G^s(u)$ 
does not depend on the window function (cf. \cite{SW}, Proposition 3.2). Finally, the following property is still true 
(cf. \cite{SW}, Proposition 3.6).

\begin{proposition}\label{WFG&WFGs}
If $u \in \mathscr{S}'(\mathbb{R}^n)$ then
\[
WF_G(u)=\overline{\bigcup_{s\in \mathbb{R}}WF_G^s(u)}\subseteq \dot{\mathbb{R}}^{2n}.
\]
\end{proposition}

Our next goal is to put in relation this Gabor wave front set and the isotropic wave front set introduced before. 

\begin{theorem}\label{theoisoG}
Let $u \in \mathscr{S}'(\mathbb{R}^n)$ and $s\in \mathbb{R}$, then
\[
WF_{\mathrm{iso}}^s(u)=WF_G^s(u).
\]
\end{theorem}
\begin{proof}
Let $X_0 \not\in WF_G^s(u)$. We want to show that $X_0 \not\in WF_{\mathrm{iso}}^s(u)$. 
By hypothesis there exists an open cone $\Gamma=\Gamma_{X_0}$ that contains $X_0$, such that 
\[
\chi_{\Gamma}\braket{\cdot}^sV_{\psi}u \in L^2(\mathbb{R}^{2n}),
\]
where, once again, $\chi_{\Gamma}$ is the characteristic function of $\Gamma$. Then, since 
$V_{\psi}^*:L^2(\mathbb{R}^{2n}) \rightarrow L^2(\mathbb{R}^n)$, we have
\[
(V_{\psi}^*\chi_{\Gamma}\braket{\cdot}^sV_{\psi})u \in L^2(\mathbb{R}^{n})
\]
and $A_{\chi_{\Gamma}\braket{\cdot}^s}=(2\pi)^{-n}(V_{\psi}^*\chi_{\Gamma}\braket{\cdot}^sV_{\psi})$.
By the Weyl-Wick correspondence 
we therefore have that there exists $b:=(\pi)^{-n}e^{-\abs*{\cdot}^2} \ast (\chi_{\Gamma}\braket{\cdot}^s)\in S_{\mathrm{iso}}^s(\mathbb{R}^n)$ such that
\[
B:=\mathrm{Op}^{\mathrm{w}}(b)=A_{\chi_{\Gamma}\braket{\cdot}^s} \quad \text{and} \quad \mathrm{Char}_{\mathrm{iso}}(B)=\mathrm{Char}_{\mathrm{iso}}(\chi_{\Gamma}\braket{\cdot}^s).
\]
Hence, we have $X_0 \not \in \mathrm{Char}_{\mathrm{iso}}(B)$. Moreover, since $\mathrm{Op}^{\mathrm{w}}(b)u \in L^2(\mathbb{R}^n)$, we conclude that $X_0 \notin WF_{\mathrm{iso}}^s(u)$. 
On the other hand, if $X_0 \notin WF_{\mathrm{iso}}^s(u)$ there exists $b \in S_{\mathrm{iso}}^s(\mathbb{R}^n)$ such that, 
with $B=\mathrm{Op}^{\mathrm{w}}(b)$,  
\[
Bu \in L^2(\mathbb{R}^n) \quad \text{and} \quad X_0 \notin \mathrm{Char}_{\mathrm{iso}}(B).
\]
Now we construct $a$ in the following way. We define 
\[
\begin{split}
a_s&:=b \ ; \\
a_{s-1}&:=-\sum_{\abs*{\alpha}=1}c_\alpha \partial^{\alpha}a_s \ ; \\
a_{s-k}&:=-\sum_{1\leq\abs*{\alpha} \leq k}c_{\alpha}\partial^{\alpha}a_{s-k+\abs*{\alpha}}, \quad \text{for all} \ k \in \mathbb{N}_0.
\end{split}
\]
where the coefficient $c_\alpha$ are given in \cite{NR}, Proposition 1.7.9 and Theorem 1.7.10. 
So, by Lemma \ref{constofa} there exists $a \in S_{\mathrm{iso}}^s(\mathbb{R}^n)$ such that 
\[
a \sim \sum_{j \geq 0} a_{s-j},
\]
and so, using again the Weyl-Wick connection, we may find $\tilde{b} \in S_{\mathrm{iso}}^s(\mathbb{R}^n)$ such that,
with $\tilde{B}:=\mathrm{Op}^{\mathrm{w}}(\tilde b)$, 
\[
\tilde{B}=A_a
\]
and
\[
\tilde{b} \sim \sum_{\alpha}c_\alpha \partial^\alpha a.
\]
Observe that $B-\tilde{B}=S$ where $S:\mathscr{S}'(\mathbb{R}^n) \rightarrow \mathscr{S}(\mathbb{R}^n)$ is continuous,
so that 
\[
A_a=\mathrm{Op}^{\mathrm{w}}(b)+S.
\]
We have proved in this way that $A_au=(2\pi)^{-n}(V_\psi^{\ast}aV_{\psi}) u \in L^2(\mathbb{R}^n)$ and then $aV_{\psi}u \in L^2(\mathbb{R}^{2n})$. 
Finally, $X_0\notin\mathrm{Char}_{\mathrm{iso}}(\mathrm{Op}^{\mathrm{w}}(a))=\mathrm{Char}_{\mathrm{iso}}(B)$.
Therefore, there exists an open cone $\Gamma=\Gamma_{X_0}$ containing $X_0$ such that 
$\Gamma \, \cap \, \mathrm{Char}_{\mathrm{iso}}(\mathrm{Op}^{\mathrm{w}}(a))= \emptyset$.
Then, since $a(X)$ is elliptic on such a cone, we have that 
\[
\chi_{\Gamma}\langle \cdot\rangle^{s}V_\psi u \in L^2(\mathbb{R}^{2n}).
\]
This proves that $X_0 \notin WF_G(u)$.
\end{proof}

\begin{remark}
As a consequence by Proposition \ref{WF&WFs} and Proposition \ref{WFG&WFGs} one has
\[
WF_G(u)=WF_{\mathrm{iso}}(u), \quad \text{for every} \ u \in \mathscr{S}'(\mathbb{R}^n).
\]
\end{remark}

\section{Main properties of the $s$-isotropic wave front set} \label{s:MainProp}
In this section we get back to isotropic wavefront sets, with the aim of giving 
a new interpretation of H\"ormander's results concerning the isotropic wave front set of \cite{HQ}, 
in terms of the $s$-isotropic wave front set. 
Through these results, in particular we obtain
a relation between the isotropic wave front set
of the kernel of the operator 
$\mathscr{K}_{e^{-2itF}}$ in \cite{HQ},
of an initial datum $u \in \mathscr{S}'(\mathbb{R}^n)$. 

First of all it is known that if $\chi$ is a linear
symplectic map in $\mathbb{R}^{2n}$, 
there exists a unitary operator $U_\chi$ from $L^2(\mathbb{R}^n)$ to $L^2(\mathbb{R}^n)$, associated with $\chi$, 
and determined up to a constant factor (cf. \cite{HoV3}, Theorem 18.5.9), such that  
\begin{equation}
\label{Metainv}
U_\chi^{-1}\mathrm{Op}^{\mathrm{w}}(a)U_\chi=\mathrm{Op}^{\mathrm{w}}(a \circ \chi),
\end{equation}
where $a \in S_\mathrm{iso}^m(\mathbb{R}^n)$ for some $m \in \mathbb{R}$. 
Thanks to this result the following proposition easily follows. 

\begin{proposition}\label{WFmetainv}
Let $u \in \mathscr{S}'(\mathbb{R}^n)$, $\chi$ be a symplectic transformation in $\mathbb{R}^{2n}$ and $U_\chi$ 
be a unitary operator on $L^2$ associated with $\chi$ as above. If $X \in \dot{\mathbb{R}}^{2n}$, one has
\begin{equation}\label{eq.WFmetainv}
X \in WF_{\mathrm{iso}}(u) \Leftrightarrow \chi(X) \in WF_{\mathrm{iso}}(U_\chi u),
\end{equation}
and more precisely, for all $s \in \mathbb{R}$ 
\begin{equation}\label{eq.WFmetainvs}
X \in WF_{\mathrm{iso}}^s(u) \Leftrightarrow \chi(X) \in WF_{\mathrm{iso}}^s(U_\chi u),
\end{equation}
\end{proposition}
\begin{proof}
We first note that \eqref{eq.WFmetainv} follows by \eqref{eq.WFmetainvs}
and by Proposition \ref{WF&WFs}, therefore it is sufficient
to show \eqref{eq.WFmetainvs}. 
To prove the latter we observe that by \eqref{Metainv}
\[
\mathrm{Op}^{\mathrm{w}}(a \circ \chi)U_\chi^{-1}u=U_\chi^{-1}\mathrm{Op}^{\mathrm{w}}(a)u,
\]
so, if $a \in S_{\mathrm{iso}}^s(\mathbb{R}^n)$ and $\mathrm{Op}^{\mathrm{w}}(a)u \in L^2(\mathbb{R}^n)$  we have
$\mathrm{Op}^{\mathrm{w}}(a \circ \chi)U_\chi^{-1}u \in L^2(\mathbb{R}^n)$, since $U_\chi$ is a unitary operator from
$L^2(\mathbb{R}^n)$ to $L^2(\mathbb{R}^n)$. Hence $\chi^{-1}(X) \notin WF_{\mathrm{iso}}^s(U_\chi^{-1}u)$ if and
only if $X\notin WF_{\mathrm{iso}}^s(u)$. 
\end{proof}
As first application of this theorem we have the following fundamental remark,
that roughly speaking says that the ($s$-)isotropic wave front set
is invariant under change of coordinates induced by
an invertible linear map - that we identify
with the associated matrix - as stated below (cf. \cite{HQ}, p. 119).
\begin{remark}\label{invcor}
Let $\kappa:\mathbb{R}^n 
\rightarrow \mathbb{R}^n$ be an invertible 
linear map and let $\chi(x,\xi)=(\kappa^{-1}x,\transp{\kappa} \, \xi)$
be the corresponding
symplectomorphism.

\noindent Therefore, 
\[
U_\chi f(x)=f(\kappa x) \abs*{\mathrm{det} \, \kappa}^{1/2}, \quad f \in \mathscr{S}(\mathbb{R}^n),
\]
so that, if $u \in \mathscr{S}'(\mathbb{R}^n)$, by Proposition \ref{WFmetainv},
for all $s \in \mathbb{R}$, we get
\[
WF_{\mathrm{iso}}^s(\kappa^\ast u)=\kappa^\ast WF_{\mathrm{iso}}^s(u):=\{(x,\transp{\kappa} \, \xi)\in\dot{\mathbb{R}}^{2n};\,(\kappa x,\xi)\in WF_{\mathrm{iso}}^{s}(u)\}.
\]
\end{remark}

We now want to relate the isotropic microsingularities
of two tempered distributions with the isotropic microsingularities of 
their tensor product. 

\begin{proposition}\label{WFtenprod}
Let $ u \in B^{-r_1}(\mathbb{R}^n)$ and $v \in B^{-r_2}(\mathbb{R}^m)$, for some $r_1,r_2 \geq 0$ and let $s_1,s_2\in \mathbb{R}$ be such that $s_*=\min\{s_1-r_2,s_2-r_1\}\leq s_1+s_2$. Then
\[
WF^{s_*}_{\mathrm{iso}}(u \otimes v) \subseteq \Bigl(WF_{\mathrm{iso}}^{s_1}(u) \cup \lbrace (0,0) \rbrace \Bigr) \times 
\Bigl(WF_{\mathrm{iso}}^{s_2}(v) \cup \lbrace (0,0) \rbrace \Bigr).
\]
\end{proposition}
\begin{proof}
We write $X=(x,\xi) \in \mathbb{R}^{2n}$ and $Y=(y,\eta) \in \mathbb{R}^{2m}$. Let $0 \neq X_0 \notin WF_{\mathrm{iso}}^{s_1}(u)$ 
and $0 \neq Y_0 \notin WF_{\mathrm{iso}}^{s_2}(v)$. Then, by definition, there exist 
$a=a(X) \in S_{\mathrm{iso}}^{s_1}(\mathbb{R}^{n})$ non-characteristic at $X_0$ and 
$b=b(Y) \in S_{\mathrm{iso}}^{s_2}(\mathbb{R}^{n})$ non-characteristic at $Y_0$, such that 
$\mathrm{Op}^{\mathrm{w}}(a)u \in L^2(\mathbb{R}^n)$ and $\mathrm{Op}^{\mathrm{w}}(b)v \in L^2(\mathbb{R}^m)$.
We define the symbol $c(x,y,\xi,\eta):=a(X)b(Y)$, then 
\[
\mathrm{Op}^{\mathrm{w}}(c)(u \otimes v) \in L^2(\mathbb{R}^{n+m}).
\]
However, this is a symbol for the split \textit{iso,iso}-metric in $\mathbb{R}^{2n}_X\times\mathbb{R}^{2m}_Y$
\[
g_2=\abs*{dX}^2/(1+\abs*{X}^2)+\abs*{dY}^2/(1+\abs*{Y}^2),
\]
but not for the isotropic metric on $\mathbb{R}^{2(n+m)}_{X,Y}$ 
\[
g_0=(\abs*{dX}^2+\abs*{dY}^2)/(1+\abs*{X}^2+\abs*{Y}^2),
\]
i.e. $c \in S_{\mathrm{iso},\mathrm{iso}}^{s_1,s_2}(\mathbb{R}^{n} \times \mathbb{R}^{m})$ (see Appendix A for the properties
of the calculus. The idea is to have operators which are a Weyl quantization separately in the $X$ and $Y$ variables,
where $(X,Y)\in\mathbb{R}^n\times\mathbb{R}^m$). 
To cope with this problem we introduce a symbol of order $0$ for the metric $g_0$ i.e. 
$\psi=\psi(x,y,\xi,\eta) \in S_{\mathrm{iso}}^0(\mathbb{R}^{n}\times \mathbb{R}^m) \subseteq S^{0,0}_{\mathrm{iso},\mathrm{iso}}(\mathbb{R}^{n} \times \mathbb{R}^{m})$
(see Proposition \ref{propSGisoiso} of the Appendix), such that 
\[
\mathrm{supp} \, \psi \subseteq \lbrace C^{-1}\leq (1+\abs*{X})/(1+\abs*{Y})\leq C \rbrace,
\]
for some $C>0$ fixed. 
In this way we may consider the composition of this symbol with respect to the metric $g_2$ to obtain
\[
\mathrm{Op}^{\mathrm{w}}(\psi)\mathrm{Op}^{\mathrm{w}}(c)=\mathrm{Op}^{\mathrm{w}}(\varphi),
\]
with $\varphi\in S^{s_1,s_2}_{\mathrm{iso},\mathrm{iso}}(\mathbb{R}^n\times\mathbb{R}^m)$.
Now, using the support of $\psi$ we have that the symbol $\varphi$ is actually a symbol for the metric $g_0$. 
Moreover $\mathrm{Op}^{\mathrm{w}}(\varphi)(u \otimes v) \in L^2(\mathbb{R}^{n+m})$. So, choosing $a,b,\psi$ 
such that $\varphi$ is non-characteristic at $(x_0,y_0,\xi_0,\eta_0)$ we get that 
$(x_0,y_0,\xi_0,\eta_0)\notin WF_{\mathrm{iso}}^{s_1+s_2}(u \otimes v)$ and then 
$(x_0,y_0,\xi_0,\eta_0)\notin WF_{\mathrm{iso}}^{s_*}(u \otimes v)$ (since $s_*\leq s_1+s_2$). 
Now, let $Y_0=0_{2m} $ and $X_0 \notin WF_{\mathrm{iso}}^{s_1}(u)$. 
By definition there exists $a=a(X) \in S_{\mathrm{iso}}^{s_1}(\mathbb{R}^n)$ non-characteristic at $X_0$ 
such that $\mathrm{Op}^\mathrm{w}(a)u \in L^2(\mathbb{R}^n)$. We next choose 
$\psi_0=\psi_0(x,y,\xi,\eta) \in S_{\mathrm{iso}}^0(\mathbb{R}^{n}\times \mathbb{R}^m)$ non-characteristic at 
$(x_0,0,\xi_0,0)$ such that 
\[
\mathrm{supp} \ \psi_0 \subseteq \lbrace 1+\abs*{Y}\leq C(1+\abs*{X}) \rbrace
\]
for some $C>0$. Then $\psi_0\in S^{0,0}_{\mathrm{iso},\mathrm{iso}}(\mathbb{R}^n\times
\mathbb{R}^m).$
Define 
$$\tilde{a}(x,y,\xi,\eta):=a(X) \otimes \mathbbm{1}_{Y}\in S^{s_1,0}_{\mathrm{iso},\mathrm{iso}}(
\mathbb{R}^n\times\mathbb{R}^m).$$ 
At this point we have to prove a lemma about the composition 
(with respect to the metric $g_2$) of $\mathrm{Op}^{\mathrm{w}}(\psi_0)$ and $\mathrm{Op}^{\mathrm{w}}(\tilde{a})$.

\begin{lemma}\label{lemmaComposition}
Let $u\in B^{-r_1}(\mathbb{R}^n)$ and $v\in B^{-r_2}(\mathbb{R}^m).$ 
Let $\tilde{a}$ and $\psi_0$ be the two symbols defined above. Then 
\[
\mathrm{Op}^{\mathrm{w}}(\psi_0) \mathrm{Op}^{\mathrm{w}}(\tilde{a})(u \otimes v) \in B^{-r_2}(\mathbb{R}^{n+m})
\]
\end{lemma}
\begin{proof}
Let 
$$\Lambda^{-r_{2}}=\mathrm{Op}^{\mathrm{w}}(\braket{(x,y,\xi,\eta)}^{-r_{2}})\in\Psi^{-r_2}_{\mathrm{iso}}(\mathbb{R}^{n+m}).$$
Then (again by Proposition \ref{propSGisoiso}), 
$$\Lambda^{-r_{2}}\mathrm{Op}^{\mathrm{w}}(\psi_0)\in\Psi^{0,-r_2}_{\mathrm{iso},\mathrm{iso}}(\mathbb{R}^n\times\mathbb{R}^m).$$
We next have to prove that 
$$\Lambda^{-r_{2}}\mathrm{Op}^{\mathrm{w}}(\psi_0)\mathrm{Op}^{\mathrm{w}}(\tilde{a})(u\otimes v)
\in L^2(\mathbb{R}^n\times\mathbb{R}^m).$$
By the hypothesis on $a$ we have that
$$\mathrm{Op}^{\mathrm{w}}(\tilde a)(u\otimes v)=\mathrm{Op}^{\mathrm{w}}(a)u\otimes v=:\tilde u\otimes v,$$
with $\tilde u\in L^2(\mathbb{R}^n).$
Let now 
$B_{1}:=\mathbbm{1}_X\otimes\mathrm{Op}^{\mathrm{w}}(\braket{Y}^{-r_{2}})\in\Psi^{0,-r_2}_{\mathrm{iso},\mathrm{iso}}(\mathbb{R}^n\times\mathbb{R}^m)$. Therefore there exists
$E_1\in\Psi^{0,r_2}_{\mathrm{iso},\mathrm{iso}}(\mathbb{R}^n\times\mathbb{R}^m)$
such that
$$E_1 B_1=I+R,$$
where $R\in\Psi^{0,-\infty}_{\mathrm{iso},\mathrm{iso}}(\mathbb{R}^n\times\mathbb{R}^m)$ is of the form
$\mathbbm{1}_X\otimes(\text{\rm smoothing operator in $y$}).$
Hence
$$E_1B_1(\tilde{u}\otimes v)=\tilde u\otimes v+\tilde{u}\otimes\tilde v,$$
with $\tilde u\in L^2(\mathbb{R}^n),\,\,v\in B^{-r_2}(\mathbb{R}^m),\,\,
\tilde v\in \mathscr{S}(\mathbb{R}^m).$ It is important to notice that
$$\Lambda^{-r_2}\mathrm{Op}^{\mathrm{w}}(\psi_0)E_1\in \Psi^{0,0}_{\mathrm{iso},\mathrm{iso}}(\mathbb{R}^n\times
\mathbb{R}^m).$$
Now we have
$$\Lambda^{-r_2}\mathrm{Op}^{\mathrm{w}}(\psi_0)(\tilde{u}\otimes v)=
\Lambda^{-r_2}\mathrm{Op}^{\mathrm{w}}(\psi_0)E_1B_1(\tilde{u}\otimes v)+r,$$
where $r\in L^2(\mathbb{R}^{n+m})$ and $B_1(\tilde u\otimes v)\in L^2(\mathbb{R}^n\times\mathbb{R}^m).$ Since
$$ \Lambda^{-r_2}\mathrm{Op}^{\mathrm{w}}(\psi_0)E_1\in
\Psi^{0,0}_{\mathrm{iso},\mathrm{iso}}(\mathbb{R}^n\times\mathbb{R}^m)$$
this concludes the proof of the lemma.
\end{proof}


\noindent By Lemma \ref{lemmaComposition}, we have thus found a symbol $\varphi=\varphi(x,y,\xi,\eta)$ defined by the relation 
\[
\mathrm{Op}^{\mathrm{w}}(\varphi)=\Lambda^{-r_2}\mathrm{Op}^{\mathrm{w}}(\psi_0)\mathrm{Op}^{\mathrm{w}}(\tilde{a})
\]
such that $\mathrm{Op}^{\mathrm{w}}(\varphi)(u \otimes v) \in L^2(\mathbb{R}^{n+m})$, $\varphi$ is non-characteristic at $(x_0,0,\xi_0,0)$ and 
$\varphi \in S_{\mathrm{iso}}^{s_1-r_2}(\mathbb{R}^n \times \mathbb{R}^m)$. The latter is due to the fact that $\varphi\in S^{s_1,-r_2}_{\mathrm{iso,iso}}(\mathbb{R}^n\times\mathbb{R}^m)$
and that on $\mathrm{supp}\,\psi_0$ one has $\langle X\rangle\leq\langle(X,Y)\rangle\leq C\langle X\rangle$.
Hence $(x_0,0,\xi_0,0) \notin WF_{\mathrm{iso}}^{s_1-r_2}(u \otimes v)$ and consequently 
$(x_0,0,\xi_0,0) \notin WF_{\mathrm{iso}}^{s_*}(u \otimes v)$. \\
The case $X_0=0_{2n}$ and $0\not=Y_0 \notin WF_{\mathrm{iso}}^{s_2}(v)$ is analogous. 
For the case $X_0\in WF_{\mathrm{iso}}^{s_1}(u)$ and $0\not=Y_0\notin WF_{\mathrm{iso}}^{s_2}(v)$ or 
$0\not=X_0\notin WF_{\mathrm{iso}}^{s_1}(u)$ and $Y_0\in WF_{\mathrm{iso}}^{s_2}(v)$ it is sufficient to choose a $\psi$ having the support on a set of type 
$\lbrace C^{-1}\leq (1+\abs*{X})/(1+\abs*{Y})\leq C \rbrace$ for some $C>0$ fixed, as in the first part of the proof and repeat the argument of the second part.
This concludes the proof of the proposition.
\end{proof}

We next want to examine the pull-back of a temperate distribution
$u \in \mathscr{S}'(\mathbb{R}^n)$
through a surjective linear map and then the case in which
the map is injective. The latter case is more complicated
and will require some extra assumptions. 

\begin{proposition}
Let $0\neq u\in\mathscr{S}'(\mathbb{R}^{n})$. Let $L:\mathbb{R}^{N}\rightarrow\mathbb{R}^{n}$
be a linear surjective map. Then, for all fixed $\mu>(N-n)/2$ and for all $s \in \mathbb{R}$
\begin{equation}\label{eq:Surjective_inclusion_1}
\begin{split}
WF_{\mathrm{iso}}^{s-\mu}(L^{\ast}u) &\subseteq L^\ast WF_{\mathrm{iso}}^{s}(u)  \cup(\mathrm{Ker}\,L\times\{0_{N}\}) \\ 
 &= \{(x,\transp{L} \,\xi')\in\dot{\mathbb{R}}^{2N};\,(L x,\xi')\in WF_{\mathrm{iso}}^{s}(u)\} \\
 & \quad \cup(\mathrm{Ker}\,L\times\{0_{N}\}).
\end{split}
\end{equation}

\end{proposition}
\begin{proof}
By Remark \ref{invcor} may assume, possibly by changing coordinates, that $L$ is the projection
\[
L:\mathbb{R}^{n}\times\mathbb{R}^{N-n}\ni(x',x'')\mapsto x'\in\mathbb{R}^{n}.
\]
The idea is to reduce the proof to showing that
\begin{equation}
\label{eq:Surjective_equivalence_1}
(x_{0}',x_{0}'',\xi'_{0},0)\in WF_{\mathrm{iso}}^{s-\mu}(L^\ast u),\,(x'_{0},\xi'_{0})\neq0\,\implies\,(x'_{0},\xi'_{0})\in WF_{\mathrm{iso}}^{s}(u),
\end{equation}
since we will show that this conditions implies
(\ref{eq:Surjective_inclusion_1}).
We first assume that this condition holds 
and we prove that it implies 
(\ref{eq:Surjective_inclusion_1}). 

\noindent To prove that (\ref{eq:Surjective_inclusion_1}) is satisfied, 
it is sufficient to notice that, since 
$\Delta_{x''}L^\ast u=0$,
by definition of the isotropic wave front set we have  
\[
WF_{\mathrm{iso}}^{s-\mu}(L^\ast u)\subseteq \mathbb{R}^{N}\times(\mathbb{R}^{n}\times\{0_{N-n}\})
\]
and then by \eqref{eq:Surjective_equivalence_1}
\[
\begin{split}
WF_{\mathrm{iso}}^{s-\mu}(L^\ast u)&=WF_{\mathrm{iso}}^{s-\mu}(L^\ast u)\cap\Bigl(\mathbb{R}^{N}\times(\mathbb{R}^{n}\times\{0_{N-n}\})\Bigl)\\
&\subseteq  \{(x,\transp{L}\xi')\in\dot{\mathbb{R}}^{2N};\,(Lx,\xi')\in WF_{\mathrm{iso}}^{s}(u)\}\cup(\mathrm{Ker}\,L\times\{0_{N}\}),
\end {split}
\]
(the second set in the above inclusion is needed 
to keep into account the case $(x',\xi')=(0,0)$).

\noindent In order to complete the proof, it remains to prove that
(\ref{eq:Surjective_equivalence_1}) holds
and to do that
we show that for all $\mu>(N-n)/2$
\begin{equation}
0\neq (x'_{0},\xi'_{0})\notin WF_{\mathrm{iso}}^{s}(u) \Rightarrow (x_{0}',x_{0}'',\xi_{0}',0)\notin WF_{\mathrm{iso}}^{s-\mu}(L^\ast u), \quad \forall x_{0}''\in\mathbb{R}^{N-n}.\label{eq:WF_F*}
\end{equation}
We assume at first that $x'_{0}\neq0$ (we will treat the case $\xi'_{0}\neq0$
later). By hypothesis 
there exists $C=\mathrm{Op}^{\mathrm{w}}(c)\in \Psi_{\mathrm{iso}}^{s}(\mathbb{R}^{n})$ 
elliptic 
at $(x'_{0},\xi'_{0})$ such that $Cu\in L^{2}(\mathbb{R}^{n})$. We then pick $\psi\in C_{c}^{\infty}(\mathbb{R}_{+};[0,1])$, define 
$\tilde\psi(x,\xi)=\psi((1+|x''|^2+|\xi''|^2)/(1+|x'|^2+|\xi'|^2))$ and consider (with $\sharp$ in the iso,iso classes)
\[
\tilde{c}(x',x'',\xi',\xi''):=\tilde{\psi}(x,\xi)\sharp \, c(x',\xi').
\]
We have that 
$\tilde\psi\in S^0_{\mathrm{iso}}(\mathbb{R}^N)$ and that
$\tilde{c}\in S_{\mathrm{iso}}^{s}(\mathbb{R}^{N})$ (by virtue of the $\mathrm{supp} \, \tilde{\psi}$) and is elliptic at $(x_{0}',x_{0}'',\xi_{0}',0)$ for all $x_0'' \in \mathbb{R}^{N-n}$,
when $\psi$ is picked with a suitably small support and $\tilde{\psi}$ is non-characteristic at $(x_0',x_0'',\xi_0',0)$, by virtue of the ellipticity of $c$ at $(x_{0}',\xi_{0}')$.
Moreover, since $L^\ast u=u\otimes\mathbbm{1}_{x''}$, denoting
$\tilde{C}:=\mathrm{Op}^{\mathrm{w}}(\tilde{c})$, we get
\[
\tilde{C}L^\ast u=\mathrm{Op}^{\mathrm{w}}(\tilde\psi)\mathrm{Op}^{\mathrm{w}}(c)(u\otimes\mathbbm{1}_{x''}) 
=\mathrm{Op}^{\mathrm{w}}(\tilde\psi)(Cu\otimes\mathbbm{1}_{x''}).
\]
Now, with $\Lambda^{-\mu}=\mathrm{Op}^{\mathrm{w}}(\langle (x',x'',\xi',\xi'') \rangle^{-\mu})$, we have $\Lambda^{-\mu}\mathrm{Op}^{\mathrm{w}}(\tilde\psi)\in\Psi_{\mathrm{iso}}^{-\mu}(\mathbb{R}^{N})$
and hence (by Proposition \ref{propSGisoiso}) also in $\Psi_{\mathrm{iso},\mathrm{iso}}^{0,-\mu}(\mathbb{R}^{n}\times\mathbb{R}^{N-n})$.

Since $\mathbbm{1}_{x''} \in B^{-\mu}(\mathbb{R}^{N-n})$ for all $\mu>(N-n)/2$ in view of Proposition \ref{prop.1in-mu}), 
we have that (see the proof of Lemma \ref{lemmaComposition})
\[
\Lambda^{-\mu}\tilde{C}L^\ast u \in L^2(\mathbb{R}^{N}).
\]
Therefore, as $\Lambda^{-\mu}\tilde{C}\in \Psi^{s-\mu}_{\mathrm{iso}}(\mathbb{R}^N)$ is elliptic at $(x_{0}',x_{0}'',\xi_{0}',0)$, we have 
\eqref{eq:WF_F*} in the case $x_{0}'\neq0$.

\noindent Now, we consider the case $\xi'_{0}\neq0$.
By Proposition \ref{WFmetainv} used with $U_{\chi}$ the normalized Fourier
transform $\mathscr{F}_{x'\to\xi'}$ in the $x'$ variable (and, hence,
$\chi:(x',\xi')\mapsto(\xi',-x')$), we prove (in the same way as when $x'_{0}\neq0$),
\[
(\xi_{0}',x_{0}'',-x'_{0},0)\in WF_{\mathrm{iso}}^{s-\mu}(L^\ast \mathscr{F}_{x'\to\xi'}u)\implies\,(\xi'_{0},-x'_{0})\in WF_{\mathrm{iso}}^{s}(\mathscr{F}_{x'\to\xi'}u),
\]
and since $\mathscr{F}_{x'\to\xi'}L^\ast u=L^\ast \mathscr{F}_{x'\to\xi'}u$ we obtain
\[
(\xi_{0}',x_{0}'',-x'_{0},0)\in WF_{\mathrm{iso}}^{s'-\mu}(\mathscr{F}_{x'\to\xi'}L^\ast u)\implies\,(\xi'_{0},-x'_{0})\in WF_{\mathrm{iso}}^{s'}(\mathscr{F}_{x'\to\xi'}u).
\]
Hence, again by Proposition \ref{WFmetainv}, 
we have \eqref{eq:WF_F*} in the case $\xi'_{0}\neq0$. 

\noindent Therefore, we have proved \eqref{eq:WF_F*} 
and this concludes the proof of the proposition. 
\end{proof}

We next examine the case in which the map is \textit{injective}.
To do that for $\Gamma \subseteq \dot{\mathbb{R}}^n$,
a closed cone,
we introduce the following space 
\[
\mathscr{S}'_{\Gamma}(\mathbb{R}^n)=\lbrace u \in \mathscr{S}'(\mathbb{R}^n); \ WF_{\mathrm{iso}}(u) \subseteq  \Gamma \rbrace
\]  
and we endow it with a notion of convergence by saying that: \\
\textit{A sequence $(u_j)_{j \in \mathbb{N}}$ converges to $0$ in $\mathscr{S}'_{\Gamma}(\mathbb{R}^n)$ if $u_j \rightarrow 0$ in $\mathscr{S}'(\mathbb{R}^n)$ for $j\rightarrow +\infty$, and for every $a \in S_{\mathrm{iso}}^0(\mathbb{R}^{n})$ with $\Gamma \, \cap \, \mathrm{essconesupp}_{\mathrm{iso}} \, a= \emptyset$ one has $\mathrm{Op}^{\mathrm{w}}(a)u_j \rightarrow 0$ in $\mathscr{S}(\mathbb{R}^n)$ for $j \rightarrow +\infty$}.
\begin{remark}\label{rmk.BsGamma}
    Note that $\mathscr{S}(\mathbb{R}^n)$ is a dense subspace of $\mathscr{S}'_\Gamma(\mathbb{R}^n)$. In fact, let $u \in \mathscr{S}'(\mathbb{R}^n)$ and $\phi \in C_c^\infty(\mathbb{R}^{2n})$, with $\phi(0)=1$. Setting $\phi_\varepsilon(X):=\phi(\varepsilon X)$, by the standard rules of calculus 
\[
\mathrm{Op}^{\mathrm{w}}(\phi_{\varepsilon})u\rightarrow u \quad \text{in} \ \mathscr{S}_{\Gamma}'(\mathbb{R}^n) \ \text{as} \ \varepsilon \rightarrow 0. 
\]
Moreover, if $u \in B^s(\mathbb{R}^n)$, for all $\delta>0$, we also have 
\begin{equation}\label{convBmu-delta}
\mathrm{Op}^{\mathrm{w}}(\phi_{\varepsilon})u\rightarrow u \quad \text{in} \ B^{s-\delta}(\mathbb{R}^n) \ \text{as} \ \varepsilon \rightarrow 0. 
\end{equation}
Indeed, by Proposition \ref{propBs-delta}, there exist $k \in \mathbb{N}_0$ and a constant $C>0$ such that
\begin{equation}\label{convBmu-delta1}
\norma*{\mathrm{Op}^{\mathrm{w}}(\phi_\varepsilon-1)u}_{B^{s-\delta}}\leq C \abs{\phi_\varepsilon-1}_k^{(\delta)}\norma*{u}_{B^s}.
\end{equation}
Therefore, since by Proposition $\ref{propmoll}$
\begin{equation}\label{convmu-delta2}
\abs{\phi_\varepsilon-1}^{(\delta)}_k\rightarrow 0 \quad \text{as} \ \varepsilon \rightarrow 0,
\end{equation}
we have that \eqref{convBmu-delta} follows from \eqref{convBmu-delta1}. 
\end{remark}

\begin{proposition}
\label{WFpull}
Let $L:\mathbb{R}^m \rightarrow \mathbb{R}^n$, with $n \geq m$ be an injective map
and let $\Gamma\subseteq \dot{\mathbb{R}}^{2n}$ be a closed cone such that
\begin{equation}
\label{GammacapF}
\Gamma \cap \lbrace(0,\xi)\in \dot{\mathbb{R}}^{2n}; \ \transp{L} \, \xi=0 \rbrace = \emptyset.
\end{equation}
Then, the pullback through $L$, defined on $\mathscr{S}(\mathbb{R}^n)$ as $(L^*u)(x')=u(Lx')$, $x'\in\mathbb{R}^m$,
can be uniquely extended as a sequentially continuous linear map 
\[
L^\ast:\mathscr{S}'_{\Gamma}(\mathbb{R}^n)\rightarrow \mathscr{S}'(\mathbb{R}^m).
\]
In addition, if $u \in \mathscr{S}'_{\Gamma}(\mathbb{R}^n)$ and $\mu>(n-m)/2$ (with $x'=(x_1,\dots,x_m)$), for all $s \in \mathbb{R}$
\begin{equation}\label{Pullback_incl}
WF_{\mathrm{iso}}^{s-\mu}(L^\ast u) \subseteq L^\ast WF_{\mathrm{iso}}^s(u)=\lbrace (x',\transp{L}\,\xi) \in \dot{\mathbb{R}}^{2m}; \ (Lx',\xi) \in WF_{\mathrm{iso}}^s(u)\rbrace.
\end{equation} 
\end{proposition}
\begin{proof}
Note first that if $n=m$ then $L:\mathbb{R}^n \rightarrow \mathbb{R}^n$ is invertible
and then the result follows from Remark \ref{invcor}. 
Moreover, since $L$ is an injective map, and, again by Remark \ref{invcor}, the result is invariant under change of coordinates,
we may restrict ourselves to the case for which $L$ is an immersion of $\mathbb{R}^m_{x'}$
into $\mathbb{R}^{m}_{x'}\times \mathbb{R}^{n-m}_{x''}$,
i.e. 
\[
L:\mathbb{R}^m \rightarrow \mathbb{R}^n, \quad  L:\mathbb{R}^m \ni x'\mapsto(x',x''=0) \in \mathbb{R}^n.
\] 
Hence, hypothesis \eqref{GammacapF} becomes
\[
\lbrace (x=0,\xi); \ \xi_j=0, \ j\leq m \rbrace \cap \Gamma=\emptyset,
\]
and, consequently, there exists $C>0$, such that, denoting $\xi'$ the dual variable of $x'$, one has
\begin{equation}
\label{Gammaprope1}
\abs*{x}+\abs*{\xi}\leq C(\abs*{x}+\abs*{\xi'}), \quad \forall (x,\xi)\in \Gamma.
\end{equation}
Now, we choose a symbol
$a=a(x,\xi) \in S_{\mathrm{iso}}^0(\mathbb{R}^{n})$
such that\\
$(i)$ $\Gamma \cap \mathrm{essconesupp}_\mathrm{iso} \, a=\emptyset$;\\
$(ii)$ with $b:=1-a \in S_{\mathrm{iso}}^0(\mathbb{R}^{n})$, one has $\mathrm{supp} \, b\subseteq \Gamma_1\cup\{0\}$  
and 
\begin{equation}
\label{Gammaprope2}
\abs*{x}+\abs*{\xi}\leq C_1(\abs*{x}+\abs*{\xi'}), \quad \forall (x,\xi)\in \Gamma_1,
\end{equation}
where $\Gamma_1$ is a closed cone that contains $\Gamma$ and $C_1>C>0$.\\
One may construct $a$ as follows. Take $\chi\in C^\infty(\mathbb{R}^{2n}_X)$ the usual excision function, that is $0\leq\chi\leq 1$, $\chi\equiv 0$ in $|X|\leq 1/2$
and $\chi\equiv 1$ in $|X|\geq 1.$ If $\omega\in C^\infty(\mathbb{S}^{2n-1})$ such that $0\leq\omega\leq 1$, $\omega\equiv 1$ on $\Gamma\cap\mathbb{S}^{2n-1}$ and 
$\mathrm{supp}\,\omega\subset\Gamma_1\cap\mathbb{S}^{2n-1}$. Define finally $a(X)=1-\chi(X)\omega(X/|X|).$

To extend $L^*$ as a sequentially continuous map from $\mathscr{S}'_\Gamma(\mathbb{R}^n)$ into $\mathscr{S}'(\mathbb{R}^m)$ we proceed as follows. The uniqueness will be a 
consequence of the sequential continuity and of the density of $\mathscr{S}(\mathbb{R}^n)$ in $\mathscr{S}'_\Gamma(\mathbb{R}^n)$.

Write
\begin{equation}
\label{usplit}
u=\mathrm{Op}^{\mathrm{w}}(a)u+\mathrm{Op}^{\mathrm{w}}(1-a)u=\mathrm{Op}^{\mathrm{w}}(a)u+\mathrm{Op}^{\mathrm{w}}(b)u,
\end{equation}
and consider the two terms separately. The plan is to define $L^*u=L^*\mathrm{Op}^{\mathrm{w}}(a)u+L^*\mathrm{Op}^{\mathrm{w}}(b)u.$ 

\noindent Let us consider first the term $\mathrm{Op}^{\mathrm{w}}(a)u$. In this case, since $WF_{\mathrm{iso}}(u)\subseteq \Gamma$ and 
\[
\mathrm{essconesupp}_\mathrm{iso} \, a\cap \Gamma=\emptyset,
\]
by Theorem \ref{isomicro} we have
\[
\mathscr{S}'_{\Gamma}(\mathbb{R}^n)\ni u \mapsto \mathrm{Op}^{\mathrm{w}}(a) u \in \mathscr{S}(\mathbb{R}^n),
\]
which is hence sequentially continuous by definition of $\mathscr{S}'_\Gamma(\mathbb{R}^n).$ 
Thus, $L^\ast \mathrm{Op}^{\mathrm{w}}(a)u\in \mathscr{S}(\mathbb{R}^m)$ and
$L^\ast \mathrm{Op}^{\mathrm{w}}(a)\colon\mathscr{S}_{\Gamma}'(\mathbb{R}^n)\rightarrow\mathscr{S}'(\mathbb{R}^m)$ is sequentially continuous.

\noindent We next investigate the second term $\mathrm{Op}^{\mathrm{w}}(b)u$, the most delicate of the two. Note first that, by Remark \ref{rmksum},
it is possible to write $u$ as a finite sum 
\begin{equation}
u=\sum_{\abs*{\alpha}+\abs*{\beta}\leq N} D^{\alpha}(x^{\beta}u_{\alpha \beta}),\label{u_sum}
\end{equation}
for some $N \in \mathbb{N}_0$,
where $D^\gamma u_{\alpha \beta} \in L^2(\mathbb{R}^n)$ for $\abs*{\gamma}\leq n$
and hence, in particular,
by the Sobolev embedding theorem,
$u_{\alpha \beta}$ is continuous. Furthermore, we note that 
if $(u_j)_{j \in \mathbb{N}} \subseteq \mathscr{S}'(\mathbb{R}^n)$,
$u_j \rightarrow 0$ in $\mathscr{S}'_{\Gamma}(\mathbb{R}^n)$ as $j \rightarrow +\infty$ for all $\alpha,\beta$ in the sum,
we may choose corresponding sequences $(u_{\alpha\beta,j})_{j}$ converging to $0$ in $L^2(\mathbb{R}^n)$.
This will shortly be useful.

\noindent By \eqref{u_sum} we may write
\[
\mathrm{Op}^{\mathrm{w}}(b)u=\sum_{\abs*{\alpha}+\abs*{\beta}\leq N}\mathrm{Op}^{\mathrm{w}}(b)D^{\alpha}(x^{\beta}u_{\alpha \beta})
\]
and then focus our attention on a fixed term of the finite sum on the RHS. To do that we
define 
\[
B_1:=\mathrm{Op}^{\mathrm{w}}((1+\abs*{x}^2+\abs*{\xi'}^2)^{\ell} \, \sharp \, \psi)
\] 
with $\ell \in \mathbb{N}$ sufficiently large to be determined and where $\psi\in C^\infty(\mathbb{R}^{2n})$, $0\leq\psi\leq 1$, 
is such that $\psi\equiv 1$ on $\mathrm{supp}\,b$ and has support in a larger open cone $\Gamma_2$ (with $\Gamma_2\cap\mathbb{S}^{2n-1}$ relatively compact)
for which \eqref{Gammaprope2} is still true (by choosing a constant $C_2>C_1$). 
It suffices to choose $\tilde\omega\in C^\infty(\mathbb{S}^{2n-1}),$ with $0\leq\tilde\omega\leq 1$ and $\tilde\omega\omega=\omega$,
$\mathrm{supp}\,\tilde\omega\subset\Gamma_2\cap\mathbb{S}^{2n-1}$, and put $\psi(X)=\chi(X)\tilde\omega(X/|X|).$

In this way $B_1\in \Psi_{\mathrm{iso}}^{2\ell}(\mathbb{R}^n)$. We have $I=B_1E+R$, where $E\in\Psi^{-2\ell}_{\mathrm{iso}}(\mathbb{R}^n)$ is a microlocal
parametrix of $B_1$ on $\Gamma_1$ and $WF'(R)\cap \Gamma_1=\emptyset$ (cf. Proposition \ref{micropmx}). Therefore,
\[
\begin{split}
\mathrm{Op}^{\mathrm{w}}(b)D^{\alpha}x^{\beta}u_{\alpha \beta}&=B_1(E\mathrm{Op}^{\mathrm{w}}(b)D^{\alpha}x^{\beta})u_{\alpha \beta}+r_{\alpha\beta} \\
&=B_1\mathrm{Op}^{\mathrm{w}}(b_{\alpha\beta})u_{\alpha\beta}+r_{\alpha\beta},
\end{split}
\]
where, for $\ell$ sufficiently large, 
we have that $b_{\alpha \beta}$ is a symbol of order $\leq 0$ and $r_{\alpha\beta} \in \mathscr{S}(\mathbb{R}^n)$ 
(note that $WF'(R)\cap WF'(\mathrm{Op}^{\mathrm{w}}(b))=\emptyset$).
Therefore, by using the Leibniz rule 
\begin{equation}\label{sumopb}
\mathrm{Op}^{\mathrm{w}}(b)u=\sum_{\abs{\alpha}+\abs{\beta}\leq N}\sum_{\abs*{\gamma}+\abs*{\delta} \leq \tilde{N}} D_{x'}^{\gamma}\Bigl(x^{\delta}\mathrm{Op}^{\mathrm{w}}(b_{\gamma \delta})u_{\alpha \beta}\Bigr)+r,
\end{equation}
for some $\tilde{N}=\tilde{N}(\alpha,\beta) \in \mathbb{N}_0$, 
$r:=\sum_{\abs{\alpha}+\abs{\beta}\leq N}r_{\alpha\beta}\in \mathscr{S}(\mathbb{R}^n)$,
$b_{\gamma \delta}$ of order $0$ and with
$D_{x'}^{\gamma}:=D_{x_1}^{\gamma_1} \dots D_{x_m}^{\gamma_m}$ 
for all $\gamma \in \mathbb{N}_0^m, \delta \in \mathbb{N}_0^n$ in the sum. At this point, by \eqref{u_sum}, we have that
\[
D^{\varepsilon}\mathrm{Op}^{\mathrm{w}}(b_{\gamma \delta})u_{\alpha \beta} \in L^2(\mathbb{R}^n)
\]
for $\abs*{\varepsilon} \leq n$, for all $\abs{\alpha}+\abs{\beta} \leq N$ and $ \abs{\gamma}+\abs{\delta}\leq \tilde{N}$. 
Hence $\mathrm{Op}^{\mathrm{w}}(b_{\gamma \delta})u_{\alpha \beta} \in C^0(\mathbb{R}^n)\cap L^\infty(\mathbb{R}^n)$. In this way, 
for a term $D_{x'}^\gamma x^{\delta}\mathrm{Op}^{\mathrm{w}}(b_{\gamma \delta})u_{\alpha \beta}$ it is possible to define $L^\ast$ as
\[
L^\ast D_{x'}^\gamma x^{\delta}\mathrm{Op}^{\mathrm{w}}(b_{\gamma \delta})u_{\alpha \beta}:=D_{x'}^\gamma L^\ast  x^{\delta}\mathrm{Op}^{\mathrm{w}}(b_{\gamma \delta})u_{\alpha \beta},
\]
because $L^\ast$ commutes with $D_{x'}^\gamma$ on sufficiently regular functions, as
$D_{x'}^\gamma$ operates only in the $x'$ variable. Moreover, since $x^\delta$ is a multiplier and 
\[
\mathrm{Op}^{\mathrm{w}}(b_{\gamma \delta})u_{\alpha \beta}\subseteq L^\infty(\mathbb{R}^n) \subseteq \mathscr{S}'(\mathbb{R}^n),
\]
one has 
\[
L^\ast D_{x'}^\gamma x^{\delta}\mathrm{Op}^{\mathrm{w}}(b_{\gamma \delta})u_{\alpha \beta}\in \mathscr{S}'(\mathbb{R}^m).
\]
By using this for every term in the sum \eqref{sumopb} we at last may define $L^\ast\mathrm{Op}^{\mathrm{w}}(b)u$ as an element of $\mathscr{S}'(\mathbb{R}^m)$. 
Finally, as mentioned earlier,  if we take a sequence $u_j \rightarrow 0$ in $\mathscr{S}'(\mathbb{R}^n)$ as $j \rightarrow +\infty$,
then for each term of the sum $\eqref{sumopb}$ we may construct a sequence of continuous functions such that the sequence of restrictions (that is, $L^*$)
converges to $0$ in $L^2(\mathbb{R}^m)$. This implies that the restriction of the sequence $(u_j)_j$ converges to $0$ in $\mathscr{S}'(\mathbb{R}^m)$.  
Thus the operator $L^\ast\mathrm{Op}^{\mathrm{w}}(b)\colon\mathscr{S}'_{\Gamma}(\mathbb{R}^n)\rightarrow\mathscr{S}'(\mathbb{R}^m)$ 
is sequentially continuous, whence

\[
L^\ast:\mathscr{S}'_{\Gamma}(\mathbb{R}^n)\rightarrow \mathscr{S}'(\mathbb{R}^m)
\]
is a sequentially continuous linear map. 

\noindent Now we prove (\ref{Pullback_incl}).
Namely, we prove that if $0\neq (x_0',\xi_0') \notin L^\ast  WF_{\mathrm{iso}}^s(u)$
then $(x_0',\xi_0') \notin WF_{\mathrm{iso}}^{s-\mu}(L^\ast u)$, for $\mu>(n-m)/2$.
More explicitly, we have to prove that if $\mu>(n-m)/2$ one has 
\[ 
(x_0',0,\xi_0',\xi'') \notin WF_{\mathrm{iso}}^s(u), \ \forall \xi'' \in \mathbb{R}^{n-m} \Rightarrow (x_0',\xi_0') \notin WF_{\mathrm{iso}}^{s-\mu}(L^\ast u).
\] 
By hypothesis we can find a conic neighborhood $V$ of $(x_0',0,\xi_0')$ such that 
\begin{equation}\label{WFu&V}
WF_{\mathrm{iso}}^s(u)\cap (V \times \mathbb{R}^{n-m})=\emptyset.
\end{equation}
We then define 
\[
c'(x,\xi'):=\phi(x,\xi')\langle (x,\xi') \rangle^{s-\mu},
\]
where $\phi \in C^\infty(\mathbb{R}_x^n\times \mathbb{R}_{\xi'}^m)$,
$0\leq \phi \leq 1$, with $\mathrm{supp} \, \phi \subseteq V$,
$\phi \equiv 1$ on $V'$ smaller conic neighborhood 
of $(x_0',0,\xi_0')$ and $\mu>(n-m)/2$.
Therefore, we note that $c(x,\xi):=c'(x,\xi') \otimes \mathbbm{1}_{\xi''}$ is a symbol
with respect to the metric
\[
g'=(\abs*{dx}^2+\abs*{d\xi'}^2)/(1+\abs*{x}^2+\abs*{\xi'}^2)
\]
and $\tilde{c}(x',\xi'):=c(x',0,\xi',0) \in S_{\mathrm{iso}}^{s-\mu}(\mathbb{R}^m)$ 
is a symbol, elliptic at $(x_0',\xi_0') \in \mathbb{R}^{2m}$,
with respect to the isotropic metric in $\mathbb{R}^{2m}$.

\noindent Now, as before, we write  
\[
u=\mathrm{Op}^{\mathrm{w}}(a)u+\mathrm{Op}^{\mathrm{w}}(b)u
\]
and examine the two terms separately.

\noindent Let us first look at $\mathrm{Op}^{\mathrm{w}}(a)u$. 
In this case, since $\mathrm{Op}^{\mathrm{w}}(a)u \in \mathscr{S}(\mathbb{R}^n)$
we have that $L^\ast \mathrm{Op}^{\mathrm{w}}(a)u \in \mathscr{S}(\mathbb{R}^m)$ and then 
\[
\mathrm{Op}^{\mathrm{w}}(\tilde{c})L^\ast (\mathrm{Op}^{\mathrm{w}}(a)u) \in \mathscr{S}(\mathbb{R}^m).
\]
We next examine $\mathrm{Op}^{\mathrm{w}}(b)u$. In the first place,
since $b\in S_{\mathrm{iso}}^0(\mathbb{R}^n)$,  
it is a symbol also for the metric $g'$.
Thus we may consider the composition with respect to the metric
$g'$ between $\mathrm{Op}^{\mathrm{w}}(b)$ and
$\mathrm{Op}^{\mathrm{w}}(c)$. By using $\eqref{Gammaprope2}$, 
which is valid on the support of $b$, we may conclude that $c \,\sharp \, b \in S_{\mathrm{iso}}^{s-\mu}(\mathbb{R}^n)$ and since $WF_{\mathrm{iso}}^s(u) \cap \, \mathrm{essconesupp}_\mathrm{iso} \, c\, \sharp \,b=\emptyset$, by Corollary \ref{corAu} we get
\[
\mathrm{Op}^{\mathrm{w}}(c \,\sharp \, b)u=\mathrm{Op}^{\mathrm{w}}(c)\mathrm{Op}^{\mathrm{w}}(b)u \in B^{\mu}(\mathbb{R}^n).
\]
Now, since 
\[
v:=\mathrm{Op}^{\mathrm{w}}(c)\mathrm{Op}^{\mathrm{w}}(b)u \in B^{\mu}(\mathbb{R}^n) \cap \mathscr{S}'_\Gamma(\mathbb{R}^n)
\]
 by Remark \ref{rmk.BsGamma}, we have that there exists a sequence $(v_j)_{j\geq 1} \subseteq \mathscr{S}(\mathbb{R}^n)$ such that 
\[
v_j \rightarrow v \quad \text{in} \ \mathscr{S}_\Gamma'(\mathbb{R}^n) \ \text{as} \ j\rightarrow +\infty,
\]
and for all $\delta>0$
\[
v_j \rightarrow v \quad \text{in} \ B^{\mu-\delta}(\mathbb{R}^n) \ \text{as} \ j\rightarrow +\infty.
\]
We then choose $\delta>0$ sufficiently small such that $\mu-\delta>(n-m)/2$ and denote by $\tau$ the trace operator 
\[
\tau : \mathscr{S}(\mathbb{R}^n)\rightarrow \mathscr{S}(\mathbb{R}^m), \quad  \mathscr{S}(\mathbb{R}^n) \ni f \mapsto \tau f(x'):=f(x',0)
\]
which extends to a bounded operator 
\[
\tau :H^{s}(\mathbb{R}^n)\rightarrow L^2(\mathbb{R}^m), \quad \text{if} \ \ s>(n-m)/2.
\]
In this way, since by Proposition \ref{prop.BsHs} we have $B^{\mu-\delta}(\mathbb{R}^n)\subseteq H^{\mu-\delta}(\mathbb{R}^n)$, we obtain that $(\tau v_j)_{j\geq 1}$ is a Cauchy sequence in $L^2(\mathbb{R}^m)$ and hence there exists $w \in L^2(\mathbb{R}^m)$ such that 
\[
\tau v_j \rightarrow w=:\tau v \quad \text{in} \ L^2(\mathbb{R}^m) \ \text{as} \ j\rightarrow +\infty \ 
\]
(note also that, as a consequence, $\tau v_j \rightarrow \tau v \ \text{in} \  \mathscr{S}'(\mathbb{R}^m) \ \text{as} \ j\rightarrow +\infty$).

\noindent Now, since $L^\ast v_j=\tau v_j$ for all $j \geq 1$ and since the map $L^\ast:\mathscr{S}_\Gamma '(\mathbb{R}^n)\rightarrow \mathscr{S}'(\mathbb{R}^m)$ is sequentially continuous, as $j \rightarrow + \infty$ we have 
\[
L^\ast v_j \rightarrow L^\ast  v=\tau v \in L^2(\mathbb{R}^m),
\]
and so
\[
L^\ast\mathrm{Op}^{\mathrm{w}}(c)\mathrm{Op}^{\mathrm{w}}(b)u \in L^2(\mathbb{R}^m).
\]
Therefore, we have that
$\tilde{C}:=\mathrm{Op}^{\mathrm{w}}(\tilde{c})\in \Psi_{\mathrm{iso}}^{s-\mu}(\mathbb{R}^m)$
is elliptic at $(x_0',\xi_0')$ and 
\[
\mathrm{Op}^{\mathrm{w}}(\tilde{c})L^\ast u=\mathrm{Op}^{\mathrm{w}}(\tilde{c})L^\ast(\mathrm{Op}^{\mathrm{w}}(a)u)+\mathrm{Op}^{\mathrm{w}}(\tilde{c})L^\ast(\mathrm{Op}^{\mathrm{w}}(b)u) \in L^2(\mathbb{R}^m),
\]
that is $(x_0',\xi_0') \notin WF_{\mathrm{iso}}^{s-\mu}(L^\ast u)$ and so we have proved the proposition.
\end{proof}
Next, as a consequence of the previous proposition, given a distribution $u\in \mathscr{S}'(\mathbb{R}^n)$,
if it satisfies condition \eqref{condrestr} below, 
essentially by integrating in the "regular variables",
we are able to define a distribution $v \in \mathscr{S}'(\mathbb{R}^m)$,
with $m \leq n$, the singularities of which are related to those of $u$. The result can be stated as follows.

\begin{proposition}\label{WFint}
Let $m,n \in \mathbb{N}$, with $n\geq m$, let $\mu>(n-m)/2$, 
and let $u \in \mathscr{S}'(\mathbb{R}^n)$ be such that 
\begin{equation}\label{condrestr}
WF_{\mathrm{iso}}(u) \cap \lbrace(0,x'',0,0); \, x'' \in \mathbb{R}^{n-m} \rbrace = \emptyset.
\end{equation}
Then 
\[
v(\cdot):=\int u(\cdot,x'')dx'' \in \mathscr{S}'(\mathbb{R}^m)
\]
can be defined by continuity,
i.e. through $L^\ast$,
as in Proposition \ref{WFpull} 
and for all $s \in \mathbb{R}$
\[
WF_{\mathrm{iso}}^{s-\mu}(v) \subseteq \lbrace (x',\xi ') \in \dot{\mathbb{R}}^{2m}; \ (x',x'',\xi',0) \in WF_{\mathrm{iso}}^{s}(u), \ \text{for some} \ x'' \in \mathbb{R}^{n-m} \rbrace.
\]
\end{proposition}
\begin{proof}
Note first that, if $f \in \mathscr{S}(\mathbb{R}^n)$ denoting by
$\hat{f}_2$ the Fourier Transform of $f$ with respect to
the $x''$ variables 
and by $L:\mathbb{R}^m \rightarrow \mathbb{R}^n$ 
the injection of $\mathbb{R}^m$ into $\mathbb{R}^m \times \mathbb{R}^{n-m}$,
(i.e $Lx'=(x',0)$)
we may write 
\[
\int f(x',x'')dx''=(L^\ast \hat{f}_2)(x').
\]
Therefore, the idea, essentially, is to use such an expression 
and the continuity of the map $L^\ast$
given in Proposition \ref{WFpull}, to give a characterization of the
$s$-isotropic wave front set of a distribution
$u \in \mathscr{S}'(\mathbb{R}^n)$ satisfying \eqref{condrestr}, as follows. 

\noindent Let $u \in \mathscr{S}'(\mathbb{R}^n)$ that satisfies \eqref{condrestr}  
and let $\hat{u}_2$ and $L$ be given as before. 
Recall that, by Proposition \ref{WFmetainv} (since for $\chi:(x_j,\xi_j) \mapsto (\xi_j,-x_j)$
the associated metaplectic operator 
is $U_\chi=(2\pi)^{-1/2}\mathscr{F}_{x_j \rightarrow \xi_j}$), 
we have
\begin{equation}
\label{MetaFourier}
WF_{\mathrm{iso}}(\hat{u}_2)=\lbrace (x',x'',\xi',\xi'') \in \dot{\mathbb{R}}^{2n}; \ (x',\xi'',\xi',-x'')\in WF_{\mathrm{iso}}(u) \rbrace ,
\end{equation}
and
\begin{equation}
\label{MetaFouriers}
WF_{\mathrm{iso}}^s(\hat{u}_2)=\lbrace (x',x'',\xi',\xi'') \in \dot{\mathbb{R}}^{2n}; \ (x',\xi'',\xi',-x'')\in WF_{\mathrm{iso}}^s(u) \rbrace.
\end{equation}
The two inclusions give rise to the following two facts. 
By \eqref{MetaFourier} we have that 
if $WF_{\mathrm{iso}}(u) \cap \lbrace(0,x'',0,0) \rbrace = \emptyset$, then 
\[
WF_{\mathrm{iso}}(\hat{u}_2)\cap \lbrace (0,\xi)\in \dot{\mathbb{R}}^{2n}; \ \transp{L} \, \xi=\xi'=0 \rbrace = \emptyset,
\]
and consequently
we may apply Proposition \ref{WFpull} to $\hat{u}_2$.
Moreover, by \eqref{MetaFouriers}
\begin{equation}\label{MetaFouriersF}
L^\ast WF_{\mathrm{iso}}^s(\hat{u}_2)\subseteq \lbrace (x',\xi') \in \dot{\mathbb{R}}^{2m}; \ (x',\xi'',\xi',0)\in WF_{\mathrm{iso}}^s(u) \ \text{for some} \ \xi'' \in \mathbb{R}^{n-m}  \rbrace.
\end{equation}
We now take
a sequence $(u_j)_{j\geq 1}\subseteq \mathscr{S}(\mathbb{R}^n)$
such that $u_j\rightarrow u$ as
$j \rightarrow +\infty$ in $\mathscr{S}'(\mathbb{R}^n)$.
In this way we obtain $L^\ast \hat{u}_{j,2} \rightarrow L^\ast  \hat{u}_2=:v$
as $j \rightarrow + \infty$, and again by Proposition \ref{WFpull}, if $\mu>(n-m)/2$,
\begin{equation}\label{Fasthatu2}
WF_{\mathrm{iso}}^{s-\mu}(L^\ast\hat{u}_2)\subseteq L^\ast  WF_{\mathrm{iso}}^s(\hat{u}_2).
\end{equation}
Hence, by \eqref{Fasthatu2} and \eqref{MetaFouriersF}, we obtain
\[
WF_{\mathrm{iso}}^{s-\mu}(L^\ast\hat{u}_2)\subseteq \lbrace (x',\xi') \in \dot{\mathbb{R}}^{2m}; \ (x',\xi'',\xi',0)\in WF_{\mathrm{iso}}^s(u) \ \text{for some} \ \xi'' \in \mathbb{R}^{n-m}\rbrace,
\]
that we may rewrite as 
\[
WF_{\mathrm{iso}}^{s-\mu}(v) \subseteq \lbrace (x',\xi') \in \dot{\mathbb{R}}^{2m}; \ (x',x'',\xi',0) \in WF_{\mathrm{iso}}^{s}(u) \ \text{for some} \ x'' \in \mathbb{R}^{n-m} \rbrace.
\]
\end{proof}

By putting together these auxiliary results
we now give the proof of Theorem \ref{WFmainprope}, but first we provide some useful notations. 

In what follows, given two sets $X_1,X_2$ and $A \subseteq X_1 \times X_2$, $B\subseteq X_2$,
we write  
\[
A\circ B=\lbrace x \in X_1; \ \exists \, y \in B \ \text{s.t.} \ (x,y) \in A \rbrace \subseteq X_1.
\]
Furthermore, if $K \in \mathscr{S}'(\mathbb{R}^{m+n})$ we denote by 
\[
WF_{\mathrm{iso},X}(K):=\lbrace X=(x,\xi) \in \dot{\mathbb{R}}^{2m}; \ (x,0,\xi,0) \in WF_{\mathrm{iso}}(K) \rbrace,
\]
by 
\[
WF_{\mathrm{iso},Y}(K):=\lbrace Y=(y,\eta) \in \dot{\mathbb{R}}^{2n}; (0,y,0,-\eta) \in WF_{\mathrm{iso}}(K)\rbrace, 
\]
and
\[
WF_{\mathrm{iso}}(K)':= \lbrace (x,y,\xi,\eta)\in\dot{\mathbb{R}}^{2(n+m)}; \ (x,y,\xi,-\eta) \in WF_{\mathrm{iso}}(K) \rbrace. 
\]
The definitions for the $s$-isotropic wave front set
are the same considering the $s$-isotropic wave front set on the RHS.

\begin{proof}[Proof of Theorem \ref{WFmainprope}]
We will see that \eqref{WFmainpropeeq1} is a
particular case of \eqref{WFmainpropeeq2} and hence we 
prove \eqref{WFmainpropeeq2}.
Note first that, if $K \in \mathscr{S}(\mathbb{R}^{m+n})$ 
and $u \in \mathscr{S}(\mathbb{R}^{n})$, we may write 
\[
\mathscr{K}u(x)=\int L^\ast (K \otimes u)(x,y)dy,
\]
where $K\otimes u\in \mathscr{S}(\mathbb{R}^{m+n+n})$ is the function 
\[
(x,y,z) \mapsto K(x,y)u(z),
\]
and $L$ is the map $(x,y)\mapsto (x,y,y)$.
Therefore, the idea is to extend this expression by continuity
as we already did in Proposition \ref{WFint}.

\noindent We begin by noting that, if $K \in \mathscr{S}'(\mathbb{R}^{n+m})$ 
and $u \in \mathscr{S}'(\mathbb{R}^n)$,
$K \otimes u \in \mathscr{S}'(\mathbb{R}^{m+n+n})$ is still well defined 
and by Proposition \ref{WFtenprod} we get
\begin{equation*}
\begin{split}
WF_{\mathrm{iso}}^{s_\ast}(K \otimes u)\subseteq & \Bigl \lbrace (x,y,z,\xi,\eta,\zeta); \ (x,y,\xi,\eta) \in WF_{\mathrm{iso}}^{s_1}(K)\cup \lbrace 0_{2(m+n)} \rbrace, \\
& (z,\zeta) \in WF_{\mathrm{iso}}^{s_2}(u) \cup \lbrace 0_{2n} \rbrace \Bigr\rbrace \setminus \lbrace 0_{2(m+n+n)} \rbrace.
\end{split}
\end{equation*}
Moreover, by H\"ormander (\cite{HQ}, Proposition 2.8), we also have 
\begin{equation}\label{RHSset}
\begin{split}
WF_{\mathrm{iso}}(K \otimes u)\subseteq & \Bigl \lbrace (x,y,z,\xi,\eta,\zeta); \ (x,y,\xi,\eta) \in WF_{\mathrm{iso}}(K)\cup \lbrace 0_{2(m+n)} \rbrace, \\
& (z,\zeta) \in WF_{\mathrm{iso}}(u) \cup \lbrace 0_{2n} \rbrace \Bigr\rbrace \setminus \lbrace 0_{2(m+n+n)} \rbrace.
\end{split}
\end{equation}
Hence, with $\theta=(\xi,\eta,\zeta)$,
one has that the RHS of \eqref{RHSset}
does not intersect the set $\lbrace (0,\theta); \ \transp{L} \, \theta=0 \rbrace$.
\noindent Indeed,
if $(x,y,z,\xi,\eta,\zeta) \in \lbrace (0,\theta); \ \transp{L} \, \theta=0 \rbrace$ 
then $x=0$, $y=0$, $z=0$, $\xi=0$, and $\eta+\zeta=0$.
However, if such a point belongs to the RHS of \eqref{RHSset},
by condition \eqref{condWF2} we have
\[
(0,\zeta) \in \Bigl(WF_{\mathrm{iso,Y}}(\mathrm{K})\cup \lbrace 0_{2n} \rbrace \Bigr)\cap \Bigl(WF_{\mathrm{iso}}(u)\cup \lbrace 0_{2n} \rbrace \Bigr)=\lbrace 0_{2n} \rbrace
\]
and so $\zeta=0$ and $\eta=0$.

\noindent Therefore, we can apply Proposition \ref{WFpull}
to $K\otimes u$ to obtain for $\mu'>((m+n+n)-(m+n))/2=n/2$
\[
\begin{split}
WF_{\mathrm{iso}}^{s_\ast-\mu'}(L^\ast(K \otimes u))  \subseteq & \Bigl\lbrace(x,y,\xi,\eta+\zeta);  \ (x,y,\xi,\eta) \in WF_{\mathrm{iso}}^{s_1}(K) \cup \lbrace 0_{2(m+n)} \rbrace, \\
&(y,\zeta) \in WF_{\mathrm{iso}}^{s_2}(u) \cup \lbrace 0_{2n} \rbrace \Bigr \rbrace \setminus \lbrace 0_{2(m+n)}\rbrace. 
\end{split}
\]
Next, again by H\"ormander (\cite{HQ} Proposition 2.9) we have
\[
\begin{split}
WF_{\mathrm{iso}}(L^\ast(K \otimes u))  \subseteq & \Bigl\lbrace(x,y,\xi,\eta+\zeta);  \ (x,y,\xi,\eta) \in WF_{\mathrm{iso}}(K) \cup \lbrace 0_{2(m+n)} \rbrace, \\
&(y,\zeta) \in WF_{\mathrm{iso}}(u) \cup \lbrace 0_{2n} \rbrace \Bigr \rbrace \setminus \lbrace 0_{2(m+n)} \rbrace
\end{split}
\]
and also in this case the RHS of the inclusion does not intersect the set $\lbrace(0,y,0,0); y\in \mathbb{R}^n \rbrace$. In fact, if $x=\xi=\eta+\zeta=0$ and
\[
(y,\zeta) \in \Bigl(WF_{\mathrm{iso,Y}}(\mathrm{K})\cup \lbrace 0_{2n} \rbrace \Bigr)\cap \Bigl(WF_{\mathrm{iso}}(u)\cup \lbrace 0_{2n} \rbrace \Bigr)=\lbrace 0_{2n} \rbrace
\]
we must have $y=0$, $\zeta=0$ and hence $\eta=0$. Therefore, by applying Proposition \ref{WFint} to $L^\ast(K \otimes u)$, the integration with respect to $y$ is well defined
and for $\mu''>((m+n)-m)/2=n/2$ we have 
\[
\begin{split}
WF_{\mathrm{iso}}^{s_\ast-\mu'-\mu''}(\mathscr{K}u)&\subseteq \Bigl\lbrace(x,\xi)\in \dot{\mathbb{R}}^{2m}; (x,y,\xi,-\eta) \in WF_{\mathrm{iso}}^{s_1}(K) \cup \lbrace 0_{2(m+n)} \rbrace, \\
& \quad \quad (y,\eta) \in WF_{\mathrm{iso}}^{s_2}(u) \cup \lbrace 0_{2n} \rbrace, \ \mathrm{for} \ \mathrm{some} \ (y,\eta) \Bigr\rbrace
\\
&= WF_{\mathrm{iso},X}^{s_1}(K)\cup WF_{\mathrm{iso}}^{s_1}(K)' \circ WF_{\mathrm{iso}}^{s_2}(u),
\end{split}
\]
that is, for $\mu>n$,
\[
WF_{\mathrm{iso}}^{s_\ast-\mu}(\mathscr{K}u) \subseteq WF_{\mathrm{iso},X}^{s_1}(K)\cup WF_{\mathrm{iso}}^{s_1}(K)' \circ WF_{\mathrm{iso}}^{s_2}(u).
\]
Finally, to prove \eqref{WFmainpropeeq1} it is sufficient
to note that if $u \in \mathscr{S}(\mathbb{R}^n)$ 
then $WF_{\mathrm{iso}}^{s_2}(u)=\emptyset$ for every $s_2 \in \mathbb{R}$.
Thus, choosing $s_2$ sufficiently large such that 
\[
s_\ast= \min\lbrace s_1-r_2, s_2-r_1 \rbrace=s_1-r_2=s_1
\]
since $r_2=0$ because $u \in \mathscr{S}(\mathbb{R}^n) \subseteq L^2(\mathbb{R}^n)=B^0(\mathbb{R}^n)$. Therefore, for any fixed $\mu>n$
\[
WF_{\mathrm{iso}}^{s_1-\mu}(\mathscr{K}u)\subseteq WF_{\mathrm{iso},X}^{s_1}(K).
\]
\end{proof}

\section{Propagation of Shubin-Sobolev singularities}\label{sec.pos}
As already mentioned in the Introduction, one of the aims of the work is to study the isotropic propagation of singularities for the problem
\begin{equation}
\label{Problem}
\begin{cases}
\partial_t u+Au=0, \quad \text{in} \ \ \mathbb{R}^+\times \mathbb{R}^n, \\
u(0,\cdot)=u_0 \in \mathscr{S}'(\mathbb{R}^{n}),
\end{cases}
\end{equation}
where $A=\mathrm{Op}^{\mathrm{w}}(a)\in \Psi_{\mathrm{iso}}^2(\mathbb{R}^n)$ and
\[
a:\; \mathbb{R}^{2n} \longrightarrow \mathbb{C}, \quad
X \longmapsto a(X)=\langle X,QX \rangle,
\]
is a \textit{complex} quadratic form defined by the symmetric matrix
$Q \in \mathsf{M}_{2n}(\mathbb{C})$ that satisfies $\mathrm{Re} \, Q \geq 0$ 
(where, once again, for $z,w \in \mathbb{C}^n$, $\langle z, w \rangle:=\sum_jz_jw_j$
is the non-Hermitian inner product). Moreover, recall that the \textit{Hamilton map} $F$ associated with $a$ (see \eqref{Hamiltonmap}) is defined by the relation
\[
\sigma(X,FX)=a(X),
\]
where $\sigma$ denotes the standard complex symplectic form. 

In what follows we assume that the linear maps that we will consider act on 
\[
\mathbb{C}^{2n}=\mathbb{R}^{2n}+i\mathbb{R}^{2n}
\]
and we sometimes use the subscript $\mathbb{R}$ to indicate that we are considering the intersection of a  subspace of $\mathbb{C}^{2n}$ with $\mathbb{R}^{2n}$. 
So, for instance, given $M \in \mathsf{M}_{2n}(\mathbb{C})$, we write
\[
\mathrm{Ker}_\mathbb{R}(M)=\lbrace w \in \mathbb{R}^{2n}; \ Mw=0 \rbrace= \lbrace w \in \mathrm{Ker}(M); \ w \in \mathbb{R}^{2n} \rbrace.
\]
Finally, we denote by $a_R=a_R(X)$ and $a_I=a_I(X)$ the real and the imaginary part of the symbol respectively, i.e. 
\[
a(X)=a_R(X)+ia_I(X).
\]

\begin{remark}
Note that a broader class of significant problems is included in our setting.
For instance, with $a(X)=\abs{\xi}^2$ we obtain the Cauchy problem for \textit{the Heat equation}, with $a(X)=i\abs{\xi}^2$, the Cauchy problem for 
\textit{the (free) Schr\"odinger equation} and finally, with $a(X)=i(\abs{x}^2+\abs{\xi}^2)$, we obtain 
the Cauchy problem for \textit{the Schr\"odinger equation of the (purely imaginary) quantum Harmonic oscillator}.
Other interesting cases will be examined below.
\end{remark}

To study such a problem, it is necessary to recall first that by \cite{HS} (p. 426) we know that there exists a $C^0$-contraction semigroup 
$\lbrace e^{-tA} \rbrace_{t\geq 0}$ on $L^2(\mathbb{R}^n)$ that extends to a semigroup $\lbrace e^{-tA} \rbrace_{t\geq 0}$ 
from $\mathscr{S}'(\mathbb{R}^n)$ to $\mathscr{S}'(\mathbb{R}^n)$.
Then, essentially, for $u_0 \in \mathscr{S}'(\mathbb{R}^n)$, our goal is to study the regularity of $u=e^{-tA}u_0$ that solves problem \eqref{Problem}.

The main goal of this section is to prove Theorem \ref{maintheo}, which states, roughly speaking, that the solution at time $t$,
with respect to the initial datum, has a loss of derivatives that depends on the regularity of the kernel of the propagator and on the dimension $n$ of the space.

To approach such a problem, we first recall the notion of Gaussian distribution introduced in \cite{HS} (Section 5) through which we will characterize the kernel of the propagator.

For $0\neq u \in \mathscr{D}'(\mathbb{R}^n)$ we denote by 
\[
\begin{split}
\mathscr{L}_u=\Bigl\lbrace {\ell \,}^\mathrm{w}(x,D) \; ; & \quad {\ell\,}^\mathrm{w}(x,D)u=0,  \\
& \ell(x,\xi)=\sum_{j=1}^na_j\xi_j+\sum_{j=1}^nb_jx_j, \ a_j,b_j \in \mathbb{C}, \forall j=1,\dots,n \Bigr\rbrace,
\end{split}
\]
and we give the following definition. 
\begin{definition}
A distribution $u\in \mathscr{D}'(\mathbb{R}^n)$ is called \textit{a Gaussian distribution} if every 
$v \in \mathscr{D}'(\mathbb{R}^n)$ such that ${\ell \,}^\mathrm{w}(x,D)v=0$ for all ${\ell \,}^{\mathrm{w}}(x,D) \in \mathscr{L}_u$, is a complex multiple of $u$.
\end{definition}

It is possible to obtain a correspondence between these distributions and the space of complex Lagrangian planes in $\mathbb{C}^{2n}$ (with respect to the 
standard complex symplectic form). To do that define
\[
\lambda_u= \lbrace (x,\xi) \in \mathbb{C}^{2n}; \ \ell(x,\xi)=0, \; \forall \; {\ell \,}^\mathrm{w}(x,D) \in \mathscr{L}_u \rbrace
\]
and recall the following result (cf. \cite{HS}, Proposition 5.1). 

\begin{proposition}
If $u$ is a Gaussian distribution, then $\lambda_u$ is a
complex Lagrangian plane such that $V:=\lbrace \xi \in \mathbb{C}^n; \; (0,\xi) \in \lambda_u \rbrace$ 
is invariant under complex conjugation, hence generated by its intersection with $\mathbb{R}^n$. 

\noindent Conversely, let $\lambda$ be a Lagrangian plane
and let $u \in\mathscr{D}'(\mathbb{R}^n)$ such that ${\ell \,}^{\mathrm{w}}(x,D)u=0$, for every $\ell=\ell(x,\xi)$ which vanishes on $\lambda$. Then, $u$ is a Gaussian distribution.
Moreover $u=ce^qd$ where $d$ is a $\delta$-function on a linear subspace, $q$ is a quadratic form there (both determined by $\lambda$) 
and $c \in \mathbb{C}\setminus \lbrace 0\rbrace$. 
\end{proposition}

Let now $u\in \mathscr{S}'(\mathbb{R}^n)$ that can be written as the oscillatory integral 
\begin{equation}
\label{ulagrangian}
u=\int_{\mathbb{R}^N}e^{ip(x,\theta)}d\theta, \quad x \in \mathbb{R}^n,
\end{equation}
where $p$ is a complex-valued quadratic form on $\mathbb{R}^{n+N}$, 
\[
p(x,\theta)=\langle (x,\theta), P(x,\theta) \rangle, \quad x \in \mathbb{R}^n, \quad \theta \in \mathbb{R}^N
\]
and $P \in \mathbb{C}^{(n+N)\times (n+N)}$ is symmetric. Writing  
\[
P=
\begin{pmatrix}
P_{x x} & P_{x \theta} \\
P_{\theta x} & P_{\theta \theta},
\end{pmatrix}
\]
we require that $\mathrm{Im} \, P\geq 0$ and that the rows of the submatrix 
\[
\begin{pmatrix}
P_{\theta x} & P_{\theta \theta} 
\end{pmatrix}
\in \mathbb{C}^{N \times (n+N)}
\]
be linearly independent over $\mathbb{C}$. 
By Proposition 5.5 in \cite{HS}, if $u\in \mathscr{S}'(\mathbb{R}^n)$ has the form \eqref{ulagrangian}, then it is a Gaussian distribution
and the corresponding Lagrangian is 
\[
\lambda=\lbrace (x,p_x'(x,\theta)) \in {\mathbb{C}}^{2n}; \ p_\theta'(x,\theta)=0, \ (x,\theta) \in \mathbb{C}^{n+N} \rbrace \subseteq \mathbb{C}^{2n}.
\]
In this setting one has the following result
for the isotropic wave front set of such distributions (cf. \cite{PRW}, Theorem 4.6).
\begin{lemma}
\label{singonlagrangian}
Let $u \in \mathscr{S}'(\mathbb{R}^n)$ be a tempered distribution
of the form \eqref{ulagrangian} and let $\lambda$ be the Lagrangian associated. Then 
\[
WF_{\mathrm{iso}}(u) \subseteq (\lambda \cap \mathbb{R}^{2n})\setminus \lbrace 0 \rbrace.
\]
\end{lemma}
\begin{remark}\label{remWFdense}
Note that since $WF_{\mathrm{iso}}(u)=\overline{\bigcup_{s \in \mathbb{R}}WF_{\mathrm{iso}}^s(u)}$ one has that for every $s \in \mathbb{R}$
\[
WF_{\mathrm{iso}}^s(u) \subseteq (\lambda \cap  \mathbb{R}^{2n})\setminus \lbrace 0 \rbrace.
\]
\end{remark}
At this point, before examining the general case we study
the propagation of singularities for the evolution 
equation
\[
\begin{cases}
\partial_t u+Au=0, \quad \text{in} \ \ \mathbb{R}^+ \times \mathbb{R}^n, \\
u(0,\cdot)=u_0 \in \mathscr{S}'(\mathbb{R}^n),
\end{cases}
\]
where $A \in \Psi_{\mathrm{iso}}^2(\mathbb{R}^n)$ is defined by
\[
A=-\Delta_{x'}+\abs{x'}^2+i(-\Delta_{x''}+\abs{x''}^2),
\]
that we may rewrite as 
\begin{equation}\label{harA1A2}
A=A_R+iA_I= \mathrm{Op}^\mathrm{w}(a_R)+i\mathrm{Op}^{\mathrm{w}}(a_I),
\end{equation}
with
\[
a_R(X)= \abs*{X'}^2, \quad a_I(X)=\abs*{X''}^2,
\]
where $X'=(x',\xi')$ for $x'\in \mathbb{R}^{n_1}$ and $\xi'\in \mathbb{R}^{n_1}$, $X''=(x'',\xi'')$ for $x''\in \mathbb{R}^{n_2}$ 
and $\xi''\in \mathbb{R}^{n_2}$ and $n_1,n_2 \in \mathbb{N}$ satisfy $n_1+n_2=n$.

We wish to study how the operator $A$ propagates singularities, in term of its isotropic wave front set. 
Roughly speaking, we expect that in the first $n_1$ variables it destroys the singularities (that is, it is regularizing), since in those variables 
it behaves like a heat evolution operator in $n_1$ variables.
More precisely, we have the following result. 

\begin{theorem}\label{theoHO}
Let $u_0 \in \mathscr{S}'(\mathbb{R}^n)$, $A$ be as in \eqref{harA1A2}
and let $t>0$ be fixed.
Then
\[
WF_{\mathrm{iso}}(e^{-tA}u_{0})\subseteq\Bigl(\lbrace 0_{n_1}\rbrace \times\mathbb{R}^{n_2}\Bigr)\times\Bigl(\lbrace0_{n_1}\rbrace \times\mathbb{R}^{n_2}\Bigr).
\]
\end{theorem}
\begin{proof}
The idea is to use \textit{the Mehler formulas} for the Harmonic Oscillator
(see \cite{HS}, pp. 424-427). First of all, we note that
if we are able to prove the result for $t<\pi/2$
(in this case the Mehler formulas are easier to handle) 
then the result follows for any $t>0$. 
Indeed, if
$t\geq \pi/2$, by fixing $t_0<\pi/2$ (since $e^{-(t-t_0)A}u_{0}\in \mathscr{S}'(\mathbb{R}^n)$),
we have that 
\[
WF_{\mathrm{iso}}(e^{-tA}u_{0})= WF_{\mathrm{iso}}(e^{-t_0A}e^{-(t-t_0)A}u_{0})\subseteq \Bigl(\lbrace 0_{n_1}\rbrace \times\mathbb{R}^{n_2}\Bigr)\times\Bigl(\lbrace0_{n_1}\rbrace \times\mathbb{R}^{n_2}\Bigr). 
\]
Therefore, we may reduce matters to proving that if $t<\pi/2$ we have
\[
(x_0',\xi_0')\neq (0,0) \implies (x_0,\xi_0) \notin WF_{\mathrm{iso}}(e^{-tA}u_{0}).
\]
Suppose first $x'_0\not=0$ (we will deal with the case $\xi_0'\neq 0$ later).
In this case, roughly speaking, the idea is to work on a cone for which the variable $x'$
controls the other variables. For such purpose we define, for some suitable constants $c>c_1>0$, 
the cones $\Gamma$ and $\Gamma_1$ containing $X_0$ (with $\Gamma\subseteq \Gamma_1$) as
\[
\begin{split}\Gamma & =\lbrace(x,\xi)\in\mathbb{R}^{2n};\;\abs*{x}\geq c\abs*{\xi},\;\abs*{x'}\geq c\abs*{x''}\rbrace,\\
\Gamma_1 & =\lbrace(x,\xi)\in\mathbb{R}^{2n};\;\abs*{x}\geq c_1\abs*{\xi},\;\abs*{x'}\geq c_1\abs*{x''}\rbrace,
\end{split}
\]
and let $b\in S_{\mathrm{iso}}^0(\mathbb{R}^n)$ be such that
$b\equiv 1$ on $\Gamma$ outside a compact neighborhood of the origin, $\mathrm{supp}\,b\subseteq\Gamma_1$.
We define  
\[
B:=\mathrm{Op}^{\mathrm{w}}(b) \in \Psi_{\mathrm{iso}}^0(\mathbb{R}^n).
\]
Since the Hamilton maps $F_R$ and $F_I$ of $a_R$ and $a_I$ commute, we may write (see Proposition 5.9 in \cite{HS} and 
also Example \ref{exaHO} below)
\begin{align*}
e^{-tA}&=e^{-t(A_R+iA_I)}\\
&=e^{-tA_R}e^{-itA_I}\\
&=\mathrm{Op}^{\mathrm{w}}(e^{-(\abs*{x'}^{2}+\abs*{\xi'}^2)/\tanh t}/\cosh t)\mathrm{Op}^{\mathrm{w}}(e^{-i(|x''|^{2}+|\xi''|^2)/\tan t}/\cos t)
\end{align*}
and then 
\[
B\,e^{-tA}=\Bigl(B\,\mathrm{Op}^{\mathrm{w}}(e^{-(\abs*{x'}^{2} + \abs*{\xi'}^{2})/\tanh t}/\cosh t)\Bigr)\mathrm{Op}^{\mathrm{w}}(e^{-i(|x''|^{2}+|\xi''|^2)/\tan t}/\cos t).
\]
Thus, since $\mathrm{Op}^{\mathrm{w}}(e^{-(\abs*{x'}^{2} + \abs*{\xi'}^{2})/\tanh t}/\cosh t)$ is smoothing in $\mathrm{supp}\,b,$ we have
\[
WF'(B)\cap WF'(\mathrm{Op}^{\mathrm{w}}(e^{-(\abs*{x'}^{2} + \abs*{\xi'}^{2})/\tanh t}/\cosh t))=\emptyset,
\]
and hence
\[
Be^{-tA}=S\,\mathrm{Op}^{\mathrm{w}}(e^{-i(|x''|^{2}+|\xi''|^2)/\tan t}/\cos t),
\]
where $S$ is a smoothing operator, and thus $B\,e^{-tA}:\mathscr{S}'(\mathbb{R}^{n})\rightarrow\mathscr{S}(\mathbb{R}^{n})$.
Therefore, since
\[
B\,e^{-tA}u_{0}\in \mathscr{S}(\mathbb{R}^{n})
\]
and $\mathrm{Char}_{\mathrm{iso}}(B)\,\cap\,\Gamma=\emptyset$ ($X_0 \in \Gamma)$,
we get
\[
X_{0}\notin WF_{\mathrm{iso}}(e^{-tA}u_{0}).
\]
The case $\xi'_0\not=0$ is analogous by choosing
\[
\begin{split}\Gamma & =\lbrace(x,\xi)\in\mathbb{R}^{2n};\;\abs*{\xi}\geq c\abs*{x},\;\abs*{\xi'}\geq c\abs*{\xi''}\rbrace;\\
\Gamma_1 & =\lbrace(x,\xi)\in\mathbb{R}^{2n};\; \abs*{\xi}\geq c'\abs*{x},\;\abs*{\xi'}\geq c'\abs*{\xi''}\rbrace.
\end{split}
\]
\end{proof}

\begin{remark}
    If we consider
    \[
    A=-\Delta_{x'}+i(-\Delta_{x''}),
    \]
   by reasoning as above, we get that for all $t>0$ 
    \[
WF_{\mathrm{iso}}(e^{-tA}u_{0})\subseteq\mathbb{R}^{n}\times\Bigl(\lbrace0_{n_1}\rbrace \times\mathbb{R}^{n_2}\Bigr).
    \]
\end{remark}

Now we want to go into details about the propagation of singularities in the general case.
According to Theorem 5.12 in \cite{HS} we have that the propagator of our problem \eqref{Problem} for $t \geq0$ can be written as 
\[
e^{-tA}=\mathscr{K}_{e^{-2itF}},
\]
where $\mathscr{K}_{e^{-2itF}}:\mathscr{S}(\mathbb{R}^n)\rightarrow\mathscr{S}'(\mathbb{R}^n)$ is the
linear operator with Schwartz kernel 
\begin{equation}\label{GaussKernel}
K_{e^{-2itF}}=(2\pi)^{-(n+N)/2}\sqrt{\mathrm{det}
\begin{pmatrix}
p^{''}_{\theta \theta}/i & p^{''}_{\theta y} 
\vspace{2mm} \\
p^{''}_{x \theta} & ip^{''}_{x y}
\end{pmatrix}
}\int e^{ip(x,y,\theta)}d\theta \in \mathscr{S}'(\mathbb{R}^{2n}).
\end{equation}
In this kernel the quadratic form $p$ is a form that defines the Lagrangian $\lambda$ associated with its twisted Lagrangian $\lambda':=\mathrm{graph}(e^{-2itF})$
(see \cite{HS} and \cite{PRW} for more details).
Using Lemma \ref{singonlagrangian} we have the following property for the kernel.

\begin{theorem}
\label{WFlambda}
Let $a(x,\xi)$ be a quadratic form defined by a symmetric matrix $Q$ with $\mathrm{Re} \, Q\geq 0$, and $F$ its Hamilton map. Then (with the notation above) for all $t>0$ we have
\[
WF_{\mathrm{iso}}(K_{e^{-2itF}})\subseteq \lbrace (x,y,\xi,-\eta) \in \dot{\mathbb{R}}^{4n}; \ X=e^{-2itF}Y, \; \mathrm{Im} \ e^{-2itF}Y=0\rbrace.
\]
\end{theorem}

We now use Theorem \ref{WFmainprope}
and Theorem \ref{WFlambda}
to prove the following fundamental lemma. 

\begin{lemma}\label{mainlem}
Let $t>0$ be fixed and let $K_{e^{-2itF}}$ be the kernel of the operator $\mathscr{K}_{e^{-2itF}}$ satisfying $K_{e^{-2itF}} \in B^{-r}(\mathbb{R}^{2n})$ 
and $u_0 \in B^{-r_{0}}(\mathbb{R}^n)$, for some $r_0,r \geq 0$.
Let also $s_{0}, s \in \mathbb{R}$ be such that
$s_\ast:= \min\lbrace s_{0}-r, s-r_{0} \rbrace\leq s_0+s$. Then, for all fixed $\mu>n$, we have
\begin{equation}
\label{WFmainlem}
\begin{split}
WF^{s_\ast-\mu}_{\mathrm{iso}}(e^{-tA}u_0) & \subseteq WF^{s}_{\mathrm{iso}}(K_{e^{-2itF}})' \circ WF^{s_{0}}_{\mathrm{iso}}(u_0) \\
& \subseteq e^{-2itF}(WF_{\mathrm{iso}}^{s_0}(u_0) \cap \mathrm{Ker}_{\mathbb{R}}(\mathrm{Im} \, e^{-2itF})),
\end{split}
\end{equation}
where $\mathrm{Ker}_\mathbb{R}(\mathrm{Im} \, e^{-2itF})=\mathrm{Ker}(\mathrm{Im} \, e^{-2itF})\cap \mathbb{R}^{2n}$.
\end{lemma} 
\begin{proof}
First of all by Theorem \ref{WFlambda} one has that $WF_{\mathrm{iso}}(K_{e^{-2itF}})$ contains neither points of the form $(0,y,0,-\eta)$ for $(y,\eta) \in  \dot{\mathbb{R}}^{2n}$, nor points of the form $(x,0,\xi,0)$ for $(x,\xi) \in \dot{\mathbb{R}}^{2n}$, i.e. using the previous notations (see Theorem \ref{WFmainprope})
\[
WF_{\mathrm{iso},X}(K_{e^{-2itF}})=WF_{\mathrm{iso},Y}(K_{e^{-2itF}})=\emptyset.
\] 
Thus, by Theorem \ref{WFmainprope}, for all fixed $\mu>n$
\[
WF_{\mathrm{iso}}^{s_\ast-\mu}(e^{-tA}u_0)\subseteq WF_{\mathrm{iso}}^{s}(K_{e^{-2itF}})' \circ WF_{\mathrm{iso}}^{s_0}(u_0).
\]
Hence
\begin{equation}
\label{usepreviousresults}
\begin{split}
WF_{\mathrm{iso}}^{s_\ast-\mu}(e^{-tA}u_0)\subseteq \lbrace X \in \mathbb{R}^{2n}; \ & \exists \; (y,\eta) \in WF_{\mathrm{iso}}^{s_0}(u_0); \\
& (x,y,\xi,-\eta) \in WF_{\mathrm{iso}}^{s}(K_{e^{-2itF}}) \rbrace 
\end{split}
\end{equation}
and so, again by Theorem \ref{WFlambda}, 
\begin{equation}
\label{Wfuselambda}
WF_{\mathrm{iso}}^{s_\ast-\mu}(e^{-tA}u_0)\subseteq e^{-2itF}(WF_{\mathrm{iso}}^{s_0}(u_0) \cap \mathrm{Ker}_{\mathbb{R}}(\mathrm{Im} \, e^{-2itF})).
\end{equation}
\end{proof}
We finally give the next result through which we obtain a lower bound for the regularity of the kernel that will lead to a lower bound for the loss of derivatives in Theorem \ref{fundtheocor} and then in Theorem \ref{maintheo}. 

\begin{proposition}\label{prop.Kregglob}
    Let $t>0$ be fixed. Then for $\varepsilon>0$ we have 
    \[
    K_{e^{-2itF}} \in B^{-n-\varepsilon}(\mathbb{R}^{2n}).
    \]
\end{proposition}
\begin{proof}
Recall that (see \eqref{GaussKernel}) $K_{e^{-2itF}}$ is a Gaussian distribution of the form 
\[
K_{e^{-2itF}}=(2\pi)^{-(n+N)/2}\sqrt{\mathrm{det}
\begin{pmatrix}
p^{''}_{\theta \theta}/i & p^{''}_{\theta y} 
\vspace{2mm} \\
p^{''}_{x \theta} & ip^{''}_{x y}
\end{pmatrix}
}\int e^{ip(x,y,\theta)}d\theta \in \mathscr{S}'(\mathbb{R}_{x,y}^{2n}).
\]
Therefore, by Proposition 5.7 in \cite{HS},
there exists an injective map $L:\mathbb{R}_\theta^N\rightarrow \mathbb{R}^{2n}_{x,y}$ and a quadratic form $q\colon\mathbb{R}_{x,y}^{2n}\rightarrow\mathbb{C}$, 
with $\mathrm{Im}\, q \geq 0$, such that 
\[
p(x,y,\theta)=\langle L\theta, (x,y)\rangle + q(x,y)
\]
and hence (again by Proposition 5.7 in \cite{HS})
\[
K_{e^{-2itF}}=c \, e^{iq(x,y)}(\, \transp{L} \,)^\ast \delta_0 \in \mathscr{S}'(\mathbb{R}_{x,y}^{2n}),
\]
for some constant $c\neq0$.
We then use this expression to prove $K_{e^{-2itF}} \in B^{-n-\varepsilon}(\mathbb{R}^{2n})$, that is equivalent (see Proposition \ref{propregu}) to 
\begin{equation*}
WF^{-n-\varepsilon}_{\mathrm{iso}}(K_{e^{-2itF}})=\emptyset.
\end{equation*}
Therefore, writing $q(x,y)=q_R(x,y)+iq_I(x,y)$, we have to show that 
\begin{equation}\label{eq.WF-n-varepsilonK}
WF_{\mathrm{iso}}^{-n-\varepsilon}(e^{-q_I}e^{iq_R}(\, \transp{L} \,)^\ast \delta_0)=\emptyset.
\end{equation}
We begin with noting that, setting $v:=e^{iq_R}(\, \transp{L} \,)^\ast \delta_0$,
one has that
\[
 WF^{s}_{\mathrm{iso}}(v)=\emptyset \Rightarrow  WF^{s}_{\mathrm{iso}}(e^{-q_I}v)=\emptyset, \quad \forall s \in \mathbb{R}.
\]
Indeed, if $v \in B^{s}(\mathbb{R}^{2n})$, we have that there exists a constant $C>0$ such that
\[
\norma*{\Lambda^se^{-q_I}\Lambda^{-s}\Lambda^sv}_0 \leq C \norma*{v}_{B^s},
\]
since $\Lambda^se^{-q_I}\Lambda^{-s}$ is a bounded operator on $L^2(\mathbb{R}^{2n})$ (see \cite{Z}, Theorem 4.23), 
where here $\Lambda^{s}=\mathrm{Op}^{\mathrm{w}}(\langle x,y,\xi,\eta \rangle^s)$, for $s \in \mathbb{R}$.
We thus have reduced the problem to showing
\[
WF^{-n-\varepsilon}_{\mathrm{iso}}(e^{iq_R}(\, \transp{L} \,)^\ast \delta_0)=\emptyset.
\]
Now, since $q_R$ is defined through a real symmetric matrix $Q_R$ we have that $e^{iq_R}$ is the metaplectic operator associated with the symplectomorphism 
\[
\chi:(x,y,\xi,\eta)\mapsto ((x,y),(\xi,\eta)-2Q_R(x,y))
\]
so that, by Proposition \ref{WFmetainv},
\[
WF^{-n-\varepsilon}_{\mathrm{iso}}(e^{iq_R}(\, \transp{L} \,)^\ast \delta_0)=\chi WF_{\mathrm{iso}}^{-n-\varepsilon}((\, \transp{L} \,)^\ast \delta_0).
\]
Therefore, we reduce \eqref{eq.WF-n-varepsilonK} to showing that
\[
WF_{\mathrm{iso}}^{-n-\varepsilon}((\, \transp{L} \,)^\ast \delta_0)=\emptyset.
\]
By possibly changing coordinates in $\mathbb{R}_{x,y}^{2n}$ and using Remark \ref{invcor}, we may suppose that 
\[
L:\mathbb{R}_\theta^{N}\rightarrow \mathbb{R}^{2n}=\mathbb{R}_{x'}^{N}\times\mathbb{R}_{y'}^{2n-N},  \quad L:\theta \mapsto (x'=\theta,y'=0),
\]
whence $(\transp{L} \,)^\ast \delta_0=\delta_0(x') \otimes \mathbbm{1}_{y'}$. 
Therefore, the result follows by reasoning as before, by virtue of the fact that $\delta_0(x') \otimes \mathbbm{1}_{y'}
=\mathscr{F}^{-1}_{\eta\rightarrow x'}(\mathbbm{1}_{\eta,y'})$, that $\mathscr{F}_{\eta\rightarrow x'}^{-1}$ is a metaplectic operator and that
$\mathbbm{1}_{\eta,y'} \in B^{-n-\varepsilon}(\mathbb{R}^{2n})$ by Proposition \ref{prop.1in-mu}.
\end{proof}

Next, note that from Remark \ref{remWFdense}, the property of the kernel of the propagator of Theorem \ref{WFlambda} holds for all $s\in\mathbb{R}$
with $WF^s_{\mathrm{iso}}$ replacing $WF_{\mathrm{iso}}.$  We therefore have the following corollary that allows us to have a loss of derivatives that
does not depend on the regularity of the initial datum. 

\begin{theorem}\label{fundtheocor}
Let $t>0$ be fixed and consider the operator $\mathscr{K}_{e^{-2itF}}:\mathscr{S}(\mathbb{R}^n)\rightarrow \mathscr{S}'(\mathbb{R}^n)$.
Then, for any given $\varepsilon>0$ (small), the Schwartz kernel $K_{e^{-2it'F}}$ of $\mathscr{K}_{e^{-2it'F}}$ belongs to 
$B^{-(n+\varepsilon)}(\mathbb{R}^{2n})$, uniformly in $t' \in (0,t]$.\\
Moreover, on defining $r_0=\inf\{r\geq 0;\,\,K_{e^{-2it'F}}\in B^{-r}(\mathbb{R}^{2n}),\,\,\forall t'\in (0,t]\}$, one has that 
if $u_0 \in \mathscr{S}'(\mathbb{R}^n)$ and $\mu>n+r_0$, then for all $s_0\in\mathbb{R}$ 
\begin{equation}\label{roughinclusion}
WF_{\mathrm{iso}}^{s_0-\mu}(e^{-tA}u_0)\subseteq e^{-2itF}(WF_{\mathrm{iso}}^{s_0}(u_0) \cap \mathrm{Ker}(\mathrm{Im} \, e^{-2itF})).
\end{equation}
\end{theorem} 
\begin{proof}
The fact that the kernel belongs to $B^{-(n+\varepsilon)}$ follows directly from Proposition \ref{prop.Kregglob}. We thus focus on proving the inclusion \eqref{roughinclusion}. 

\noindent If 
\[
u_0 \in \mathscr{S}'(\mathbb{R}^n)=\bigcup_{r\in \mathbb{R}} B^r(\mathbb{R}^n),
\]
there exists $r_0 \geq 0$ such that $u_0 \in B^{-r_0}(\mathbb{R}^n)$.
Moreover, by Theorem \ref{WFlambda}, we have
\[
WF_{\mathrm{iso}}^{s}(K_{e^{-2itF}})\subseteq \lbrace (x,y,\xi,-\eta) \in \dot{\mathbb{R}}^{4n} ; \ X=e^{-2itF}Y, \; \mathrm{Im} \ e^{-2itF}Y=0\rbrace,
\]
for all $s \in \mathbb{R}$. Hence, if we choose $s\geq \max \lbrace r_0+ s_0-r,-r \rbrace$, we have that
\[
s_\ast=\min \lbrace s_0-r,s-r_0 \rbrace=s_0-r \leq s_0+s.
\]
Therefore, the inclusion \eqref{roughinclusion} is a consequence of Lemma \ref{mainlem}.
\end{proof}

The goal now is to use the geometry of the \textit{singular space}
introduced in \cite{HPS} (see also \cite{PRW}) 
\[
S= \Bigl(\bigcap_{j=0}^{2n-1} \mathrm{Ker}(\mathrm{Re} F(\mathrm{Im}F)^j)\Bigr)\cap \mathbb{R}^{2n} \subseteq  \mathbb{R}^{2n},
\]
to refine the inclusion \eqref{roughinclusion}.
The singular space enjoys important properties (see \cite{PRW} pp. 147-149), 
among which are:
\begin{equation}\label{propSing1}
S=\Bigl(\bigcap_{0\leq t' \leq t}\mathrm{Ker}(\mathrm{Im} \, e^{-2it'F})\Bigr)\cap \mathbb{R}^{2n}
\end{equation}
and $S$ is stable for $\mathrm{Re} \, F$ and $\mathrm{Im} \, F$, that is 
\begin{equation}\label{propSing2}
(\mathrm{Re} \, F)S=\lbrace 0 \rbrace, \quad (\mathrm{Im} \, F)S \subseteq S.
\end{equation}

We are in a position to prove our main result concerning propagation of singularities, that is Theorem \ref{maintheo}.
\begin{proof}[Proof of Theorem \ref{maintheo}]
First of all, by Theorem \ref{fundtheocor}, for $\mu'>n+r$ (with $r\leq n+\varepsilon$) we have
\[
WF^{s-\mu'}_{\mathrm{iso}}(e^{-tA}u_0) \subseteq e^{-2itF}(WF_{\mathrm{iso}}^{s}(u_0) \cap \mathrm{Ker}(\mathrm{Im} \ e^{-2itF})).
\]
Now we note that, since $WF_{\mathrm{iso}}^{s}(u_0)\subseteq \mathbb{R}^{2n}$ and $e^{-2itF}$ is an invertible linear map (see Lemma 5.2 in \cite{PRW}),
we may write 
\[
\begin{split}
&e^{-2itF}(WF_{\mathrm{iso}}^{s}(u_0)\cap \mathrm{Ker}(\mathrm{Im} \, e^{-2itF})) \\
&=e^{-2itF}(WF_{\mathrm{iso}}^{s}(u_0)\cap \mathrm{Ker}(\mathrm{Im} \, e^{-2itF})\cap \mathbb{R}^{2n})\\
&=(e^{-2itF}WF_{\mathrm{iso}}^{s}(u_0))\cap (e^{-2itF}(\mathrm{Ker}(\mathrm{Im} \, e^{-2itF})\cap \mathbb{R}^{2n})).
\end{split}
\]
Now, setting $Z=e^{-2itF}W$, with $Z,W \in \mathbb{C}^{2n}$,
\[
Z \in \mathrm{Ker}(\mathrm{Im} \ e^{2itF}) \cap \mathbb{R}^{2n} \Longleftrightarrow W \in \mathrm{Ker}(\mathrm{Im} \, e^{-2itF})\cap \mathbb{R}^{2n}.
\]
Therefore, it follows that 
\[
e^{-2itF}(WF_{\mathrm{iso}}^{s}(u_0)\cap \mathrm{Ker}(\mathrm{Im} \, e^{-2itF}))=e^{-2itF}(WF_{\mathrm{iso}}^{s}(u_0))\cap \mathrm{Ker}(\mathrm{Im} \, e^{2itF})\cap  \mathbb{R}^{2n},
\]
and then
\begin{equation}
\label{eqpropagation}
WF^{s-\mu'}_{\mathrm{iso}}(e^{-tA}u_0)\subseteq (e^{-2itF}WF_{\mathrm{iso}}^{s}(u_0))\cap \mathrm{Ker}(\mathrm{Im} \, e^{2itF})\cap \mathbb{R}^{2n}.
\end{equation}
Now by Theorem 5.9 in \cite{HS} for all $t\geq 0$ the linear continuous operator 
\[
\mathscr{K}_{e^{-2itF}}:\mathscr{S}(\mathbb{R}^n)\rightarrow \mathscr{S}'(\mathbb{R}^n),
\]
extends to a continuous map 
\[
\mathscr{K}_{e^{-2itF}}:\mathscr{S}'(\mathbb{R}^n)\rightarrow \mathscr{S}'(\mathbb{R}^n),
\]
satisfying 
\[
\mathscr{K}_{e^{-2i(t_1+t_2)F}}=\mathscr{K}_{e^{-2it_1F}}\mathscr{K}_{e^{-2it_2F}}
\]
and, then, for all $t_1,t_2\geq 0$
\begin{equation}
\label{tsplit}
e^{-(t_1+t_2)A}=e^{-t_1A}e^{-t_2A},
\end{equation}
Now choosing $t_1,t_2>0$ such that $t_1+t_2=t$ and using \eqref{eqpropagation} and \eqref{tsplit}, with $\mu=2\mu'$ and $\mu'>n+r$, we get
\[
\begin{split}
& WF_{\mathrm{iso}}^{s-2\mu'}(e^{-tA}u_0)\\
&=WF_{\mathrm{iso}}^{s-2\mu'}(e^{-t_1A}e^{-t_2A}u_0) \\
& \subseteq (e^{-2it_1F}WF_{\mathrm{iso}}^{s-\mu'}(e^{-t_2A}u_0))\cap \mathrm{Ker}(\mathrm{Im} \, e^{2it_1F})\cap \mathbb{R}^{2n} \\
& \subseteq \Bigl( e^{-2it_1F}((e^{-2it_2F}(WF_{\mathrm{iso}}^{s}(u_0))\cap \mathrm{Ker}(\mathrm{Im}\, e^{2it_2F})\cap  \mathbb{R}^{2n})\Bigr) \cap \mathrm{Ker}(\mathrm{Im} \, e^{2it_1F})\cap \mathbb{R}^{2n} \\
&\subseteq (e^{-2itF}WF_{\mathrm{iso}}^{s}(u_0))\cap \Bigl((e^{-2it_1F}(\mathrm{Ker}(\mathrm{Im} \, e^{2it_2F})\cap \mathbb{R}^{2n})\Bigr)\cap \mathrm{Ker}(\mathrm{Im} \, e^{2it_1F})\cap \mathbb{R}^{2n} \\
&\subseteq (e^{-2itF}WF_{\mathrm{iso}}^{s}(u_0))\cap \mathrm{Ker}(\mathrm{Im} \, e^{-2it_1F})\cap  \mathbb{R}^{2n}.
\end{split}
\]
Hence, 
\[
WF_{\mathrm{iso}}^{s-\mu}(e^{-tA}u_0)\subseteq (e^{-2itF}WF_{\mathrm{iso}}^{s}(u_0))\cap \Bigl(\bigcap_{0\leq t' \leq t}\mathrm{Ker}(\mathrm{Im} \, e^{-2it'F})\Bigr)\cap \mathbb{R}^{2n}
\]
and, then, by \eqref{propSing1},
\begin{equation}
\label{Usesing1}
WF_{\mathrm{iso}}^{s-\mu}(e^{-tA}u_0)\subseteq (e^{-2itF}WF_{\mathrm{iso}}^{s}(u_0))\cap S.
\end{equation}
If we now take $X \in (e^{-2itF}WF_{\mathrm{iso}}^{s}(u_0))\cap S$,
we find $Y \in WF_{\mathrm{iso}}^{s}(u_0) \subseteq  \mathbb{R}^{2n}$ such that $X=e^{-2itF}Y \in S\subseteq \mathbb{R}^{2n}$. Then
\begin{equation}
\label{Usesing2}
Y=e^{2itF}X=\sum_{j=0}^{+\infty}\frac{(2itF)^jX}{j!}=\sum_{j=0}^{+\infty}\frac{(-2t\mathrm{Im} F)^jX}{j!}=e^{-2t\mathrm{Im}F}X,
\end{equation}
since, by \eqref{propSing2}, $(iF)^jX=(-\mathrm{Im} \, F)^jX$ for all $j \geq 0$. We deduce that 
\[
X=e^{2t\mathrm{Im}F}Y \in (e^{2t\mathrm{Im}F}WF_{\mathrm{iso}}^{s}(u_0))\cap S.
\]
Conversely, from \eqref{Usesing2}, if $X \in (e^{2t\mathrm{Im} F}WF_{\mathrm{iso}}^{s}(u_0))\cap S$, then 
\[
Y=e^{-2t\mathrm{Im}F}X=e^{2itF}X \in WF_{\mathrm{iso}}^{s}(u_0).
\]
On the other hand this implies that 
\begin{equation}
\label{Usesing3}
\begin{split}
(e^{-2itF}WF_{\mathrm{iso}}^{s}(u_0))\cap S&=(e^{2t\mathrm{Im}F}WF_{\mathrm{iso}}^{s}(u_0))\cap S\\
&=e^{2t\mathrm{Im}F}(WF_{\mathrm{iso}}^{s}(u_0)\cap (e^{-2t\mathrm{Im}F}S)).
\end{split}
\end{equation}
Now, since by \eqref{propSing2} the singular space $S$ is invariant
with respect to $\mathrm{Im} \, F$, we get
\[
e^{-2t\mathrm{Im}F}S\subseteq S.
\]
Moreover, as the linear map $e^{-2t\mathrm{Im}F}$ is invertible, we have that  
\[
{e^{-2t\mathrm{Im}F}}_{\bigl|_S}:S\rightarrow S,
\]
is an injective endomorphism and since $S$ is finite dimensional it is also surjective. Therefore 
\[ 
e^{-2t\mathrm{Im}F}S=S.
\]
Hence, putting together \eqref{Usesing1} and \eqref{Usesing3} we obtain 
that 
\[
WF_{\mathrm{iso}}^{s-\mu}(e^{-tA}u_0)\subseteq (e^{2t\mathrm{Im}F}(WF_{\mathrm{iso}}^{s}(u)\cap S))\cap S,
\]
and hence, by \eqref{lin.HVF}
\[
WF_{\mathrm{iso}}^{s-\mu}(e^{-tA}u_0)\subseteq (e^{tH_{a_I}}(WF_{\mathrm{iso}}^{s}(u)\cap S))\cap S.
\]
This concludes the proof of our theorem.
\end{proof}
\begin{remark}
    Let us stress that, in our inclusion \eqref{incmaintheo}, we have a maximum derivative loss of $4n+\varepsilon$, with $\varepsilon>0$ arbitrarily small, since by Proposition \ref{prop.Kregglob} the kernel of the propagator has (at least) Shubin-Sobolev regularity $-(n+\varepsilon)$. 
\end{remark}

Let us now consider the particular case when the operator $\mathrm{Op}^{\mathrm{w}}(a)$ commutes with its adjoint $\mathrm{Op}^{\mathrm{w}}(\bar{a})$, i.e.
\[
\mathrm{Op}^{\mathrm{w}}(a)\mathrm{Op}^{\mathrm{w}}(\bar{a})u=\mathrm{Op}^{\mathrm{w}}(\bar{a})\mathrm{Op}^{\mathrm{w}}(a)u, \quad u \in \mathscr{S}(\mathbb{R}^n).
\] 
This condition is equivalent to the condition on the Poisson bracket
\[
\lbrace a,\bar{a} \rbrace= \langle \nabla_\xi a,\nabla_x\bar{a} \rangle-\langle \nabla_x a,\nabla_\xi\bar{a} \rangle=-2i \lbrace a_R,a_I \rbrace=0,
\]
and in this case (cf. \cite{PRW}) 
\[
S=\mathrm{Ker}(\mathrm{Re} \, F) \cap \mathbb{R}^{2n}=\mathrm{Ker}_{\mathbb{R}}(\mathrm{Re} \, F) \cap \mathbb{R}^{2n}.
\]
\begin{corollary}
Let $a$ be a quadratic form on $ \mathbb{R}^{2n}$ such that $\lbrace a,\bar{a} \rbrace=0$. Then, under the hypotheses of Theorem \ref{maintheo}, for $\mu>2(n+r)$ (with $r\leq n+\varepsilon$) and for all $s \in \mathbb{R}$
\[
WF_{\mathrm{iso}}^{s-\mu}(e^{-tA}u_0) \subseteq (e^{tH_{a_I}}(WF_{\mathrm{iso}}^{s}(u_0) \cap \mathrm{Ker}_\mathbb{R}(\mathrm{Re} \, F))\cap \mathrm{Ker}_{\mathbb{R}}(\mathrm{Re} \, F).
\]
\end{corollary}
This corollary allows us to treat the case in which the quadratic form $a=a(X)$
is \textit{real-valued}. In this context we obtain a more precise information about the propagation of singularities as we see in
the following two examples (cf. \cite{PRW}).
\begin{example}
Let $u_0 \in \mathscr{S}'(\mathbb{R}^n)$ and $a(X)=\abs{\xi}^2$. In this case, as we already pointed out, we have that \eqref{Problem} corresponds to the Cauchy problem for \textit{the Heat equation} 
\begin{equation}
\label{Problem1}
\begin{cases}
\partial_t u-\Delta u=0, \quad \text{in} \ \ \mathbb{R}^+\times \mathbb{R}^n,\\
u(0,\cdot)=u_0 \in \mathscr{S}'(\mathbb{R}^{n}),
\end{cases}
\end{equation}
and since in this case 
$S=\mathrm{Ker}(\mathrm{Re} \, F) \cap \mathbb{R}^{2n}$, we have 
\[
WF_{\mathrm{iso}}^{s-\mu}(e^{-tA}u_0) \subseteq WF_{\mathrm{iso}}^{s}(u_0) \cap (\mathbb{R}^n \times \lbrace 0 \rbrace)
\]
and 
\[
WF_{\mathrm{iso}}(e^{-tA}u_0) \subseteq WF_{\mathrm{iso}}(u_0) \cap (\mathbb{R}^n \times \lbrace 0 \rbrace).
\]
\end{example}
\begin{example}
Let $u_0 \in \mathscr{S}'(\mathbb{R}^n)$ and let
\[
a=a(x)=\langle x, M x \rangle,  
\]
where $0 \leq M \in \mathsf{M}_n(\mathbb{R})$ is symmetric (and where $\langle \cdot, \cdot \rangle$ is the euclidean scalar product on $\mathbb{R}^n$). Then
\[
WF_{\mathrm{iso}}^{s-\mu}(e^{-tA}u_0) \subseteq WF_{\mathrm{iso}}^{s}(u_0) \cap (\mathrm{Ker}(M) \times \lbrace 0 \rbrace)
\]
and 
\[
WF_{\mathrm{iso}}(e^{-tA}u_0) \subseteq WF_{\mathrm{iso}}(u_0) \cap
 (\mathrm{Ker}(M) \times \lbrace 0 \rbrace).
\]
\end{example}
We finally use the results that we proved to refine the inclusion showed in Theorem \ref{theoHO} and obtain a more precise information about the propagation of a quantum \textit{Harmonic Oscillator}. 
To do that recall the notation $X'=(x',\xi')$ with $x'\in \mathbb{R}^{n_1}$, $\xi'\in \mathbb{R}^{n_1}$ and $X''=(x'',\xi'')$ with $x''\in \mathbb{R}^{n_2}$,  
$\xi''\in \mathbb{R}^{n_2}$, where $n_1,n_2 \in \mathbb{N}_0$ are such that $n_1+n_2=n$.

\begin{example}\label{exaHO}
Let us consider the propagation of 
\[
\begin{cases}
\partial_{t}u+Au=0, \quad \text{in} \ \mathbb{R}^+ \times \mathbb{R}^n, \\
u(0,\cdot)=u_{0}\in \mathscr{S}'(\mathbb{R}^{n}),
\end{cases}
\]
where $A=\mathrm{Op}^{\mathrm{w}}(a)\in \Psi^2_\mathrm{iso}(\mathbb{R}^n)$ is defined as 
\[
A=-\Delta_{x'}+\abs{x'}^2+i(-\Delta_{x''}+\abs{x''}^2),
\]
that is $A$ is the Weyl-quantization of the symbol
\[
a(X)=|X'|^2+i|X''|^2.
\]
Let now $0<t<\pi/2$ be fixed. 
In this case, for $t' \in (0,t]$ we have (see \cite{HS} pp. 425-427)
\begin{align*}
&K_{e^{-2it'F}}(x,y)\\
&= \exp\Bigl(-\frac{1}{2}((\abs{x'}^{2}+\abs{y'}^{2})\cosh(2t')-2\langle x',y' \rangle)/\sinh(2t')\Bigr)/\sqrt{2\pi\sinh(2t')}\\
 & \quad \times \exp\Bigl(\frac{i}{2}((\abs{x''}^{2}+\abs{y''}^{2})\cos(2t')-2\langle x'',y''\rangle)/\sin(2t')\Bigr)/\sqrt{2\pi i\sin(2t')},
\end{align*}
where here $\sqrt{2\pi i\sin(2t')}$ is a fixed complex root of $z=2\pi i\sin(2t')$.

\noindent Then, for $\varepsilon>0$, since $K_{e^{-2it'F}}$ is Schwartz in the $(x',y')$ variables, by Proposition \ref{prop.1in-mu} and for all $t' \in (0,t]$ we have
\[
K_{e^{-2it'F}} \in H^{-n_2-\varepsilon,0}(\mathbb{R}^{2n})\subseteq B^{-n_2-\varepsilon}(\mathbb{R}^{2n}).
\] 
Moreover,
\[
Q=\left(\begin{array}{c|c}
\begin{array}{cc}
I_{n_1} & 0_{n_1\times n_2}\\
0_{n_2\times n_1} & iI_{n_2}
\end{array} & 0_{n\times n}\\
\hline 0_{n\times n} & \begin{array}{cc}
I_{n_1} & 0_{n_1\times n_2}\\
0_{n_2\times n_1} & iI_{n_2}
\end{array}
\end{array}\right),
\]
and so
\[
F=JQ=\left(\begin{array}{c|c}
0_{n\times n} & \begin{array}{cc}
I_{n_1} & 0_{n_1\times n_2}\\
0_{n_2\times n_1} & iI_{n_2}
\end{array}\\
\hline \begin{array}{cc}
-I_{n_1} & 0_{n_1\times n_2}\\
0_{n_2\times n_1} & -iI_{n_2}
\end{array} & 0_{n\times n}
\end{array}\right).
\]
Therefore
\[
\mathrm{Im}F=\frac{F-\bar{F}}{2i}=\left(\begin{array}{c|c}
0_{n\times n} & \begin{array}{cc}
0_{n_1} & 0_{n_1\times n_2}\\
0_{n_2\times n_1} & I_{n_2}
\end{array}\\
\hline \begin{array}{cc}
0_{n_1\times n_1} & 0_{n_1\times n_2}\\
0_{n_2\times n_1} & -I_{n_2}
\end{array} & 0_{n\times n}
\end{array}\right)
\]
and 
\[
\mathrm{Re}F=\frac{F+\bar{F}}{2}=\left(\begin{array}{c|c}
0_{n\times n} & \begin{array}{cc}
I_{n_1} & 0_{n_1\times n_2}\\
0_{n_2\times n_1} & 0_{n_2}
\end{array}\\
\hline \begin{array}{cc}
-I_{n_1} & 0_{n_1\times n_2}\\
0_{n_2\times n_1} & 0_{n_2}
\end{array} & 0_{n \times n}
\end{array}\right).
\]
Hence, since $(\mathrm{Re} \, F\mathrm{Im} \, F)=0$, in this case
\[
S=\mathrm{Ker}_\mathbb{R}(\mathrm{Re} \, {F})=(\{0_{n_1}\}\times\mathbb{R}^{n_2})\times(\{0_{n_1}\}\times\mathbb{R}^{n_2})
\]
and as 
\[
e^{2t\mathrm{Im}\,F}=\left(\begin{array}{c|c}
\begin{array}{cc}
I_{n_1} & 0_{n_1\times n_2}\\
0_{n_2\times n_1} & \cos(2t)I_{n_2}
\end{array} & \begin{array}{cc}
0_{n_1} & 0_{n_1\times n_2}\\
0_{n_2\times n_1} & \sin(2t)I_{n_2}
\end{array}\\
\hline \begin{array}{cc}
0_{n_1} & 0_{n_1\times n_2}\\
0_{n_2\times n_1} & -\sin(2t)I_{n_2}
\end{array} & \begin{array}{cc}
I_{n_1} & 0_{n_1\times n_2}\\
0_{n_2\times n_1} & \cos(2t)I_{n_2}
\end{array}
\end{array}\right)
\]
we may conclude that, denoting $\pi_{n_2}:(x,\xi)\mapsto (x'',\xi'')$
the projection onto the $n_2$-group of variables, for $0<t<\pi/2$, for $s \in \mathbb{R}$ and for $\varepsilon>0$,
\begin{equation}\label{eqHarn1n2}
\begin{split}
&WF_{\mathrm{iso}}^{s-2n-2n_2-\varepsilon}(e^{-tA}u_{0}) \\
&\subseteq\Bigl(e^{2t\mathrm{Im}F}(WF_{\mathrm{iso}}^{s}(u_{0})\cap S)\cap S\Bigr)\\
 & \begin{split}= & \Bigl(\{0_{n_1}\}\times\bigl\lbrace \cos(2t)x'' +\sin(2t)\xi''; \; (x'',\xi'') \in \pi_{n_2} WF_{\mathrm{iso}}^{s}(u_0) \bigr \rbrace\Bigr) \\
 & \times \Bigl( \{0_{n_1}\}\times \bigl\lbrace -\sin(2t)x'' +\cos(2t)\xi''; \; (x'',\xi'') \in \pi_{n_2} WF_{\mathrm{iso}}^{s}(u_0) \bigr \rbrace\Bigr).
\end{split}
\end{split}
\end{equation}
\end{example}
\begin{remark}
Note that, in the previous example, when $n_1=n$ we have 
\[
A=-\Delta+\abs{x}^2
\]
and for all $s \in \mathbb{R}$ we get
\[
WF_{\mathrm{iso}}^{s-2n-\varepsilon}(e^{-tA}u_{0})=\emptyset.
\]
When $n_2=n$ we have
\[
A=i(-\Delta+\abs{x}^2)
\]
for which
\[
\begin{split}
WF_{\mathrm{iso}}^{s-4n-\varepsilon}(e^{-tA}u_{0})\subseteq &\Bigl(\bigl\lbrace \cos(2t)x +\sin(2t)\xi; \; (x,\xi) \in  WF_{\mathrm{iso}}^{s}(u_0) \bigr \rbrace\Bigr)\\
&\times \Bigl( \bigl\lbrace -\sin(2t)x +\cos(2t)\xi; \; (x,\xi) \in WF_{\mathrm{iso}}^{s}(u_0) \bigr\rbrace \Bigr),
\end{split}
\]
for all $s \in \mathbb{R}$.
\end{remark}
\begin{remark}
In the latter example we assumed  $t<\pi/2$, with the aim of making an explicit calculation through Mehler's formulas. Nevertheless, once more, using Proposition \ref{prop.Kregglob}, 
we get the inclusion \eqref{eqHarn1n2} for all $t>0$, with a maximum derivative loss of $4n+\varepsilon$, for $\varepsilon>0$ arbitrarily small. 
\end{remark}


\appendix
\section{Auxiliary results}\label{sec.app}
In this section we collect some helpful results. We begin by showing the following useful properties on the isotropic context, and in order to do so, it is convenient to recall 
that, for $m \in \mathbb{R}$, $S_{\mathrm{iso}}^m(\mathbb{R}^n)$ is a Fr\'echet space with semi-norms given by \eqref{semi-norms} and, moreover, for 
$s \in \mathbb{R}$, $\Lambda^s:=\mathrm{Op}^{\mathrm{w}}(\langle X \rangle^s)$.
\begin{proposition}\label{propmoll}
    Let $a \in S^0_\mathrm{iso}(\mathbb{R}^n)$ and set $a_\varepsilon(X):=a(\varepsilon X)$, with $\varepsilon \in [0,1]$. Then $(a_\varepsilon)_\varepsilon$ is 
bounded in $S_{\mathrm{iso}}^0(\mathbb{R}^n)$ and for all $m>0$,
    \[
    a_\varepsilon \rightarrow a_0 \quad \text{in} \ S^m_{\mathrm{iso}}(\mathbb{R}^n) \ \ \text{as} \ \varepsilon \rightarrow 0.
    \]
\end{proposition}
\begin{proof}
The proof of this fact is the analogous (in the isotropic framework) of the proof of Proposition 18.1.2 in \cite{HoV3}. 
We first note that, since $a_0=a(0)$ is constant, it is sufficient to prove that for a fixed $m \in [0,1]$, we have that for all $\alpha \in \mathbb{N}_0^n$ there exists a constant $C_\alpha>0$ such that 
\begin{equation}\label{convarepsilon}
\sup_{X \in \mathbb{R}^{2n}}\abs{(1+\abs{X})^{\abs{\alpha}-m}\partial_X^\alpha(a_\varepsilon(X)-a_0)}\leq C_\alpha\varepsilon^m, \quad \varepsilon \in (0,1].
\end{equation}
If $\alpha=0$ we have that \eqref{convarepsilon} follows from the fact that, by Taylor's formula, 
\[
\abs{a_\varepsilon(X)-a_0}\leq C \abs{\varepsilon X}^m, 
\]
for some constant $C>0$. 

\noindent Finally, if $\alpha \neq 0$, we have that \eqref{convarepsilon} follows from the inequalities
\[
(1+\abs{X})^{\abs{\alpha}-m}(1+\abs{\varepsilon X})^{-\abs{\alpha}}\varepsilon^{\abs{\alpha}-m}\leq \Bigl((1+\abs{X})/(1/\varepsilon +\abs{X})\Bigr)^{\abs{\alpha}-m}\leq 1,
\]
that hold for all $X \in \mathbb{R}^{2n}$ and $\varepsilon \in (0,1]$.
\end{proof}
\begin{proposition}\label{propBs-delta}
Let $A=\mathrm{Op}^\mathrm{w}(a) \in \Psi^m_{\mathrm{iso}}(\mathbb{R}^n)$, with $m \in \mathbb{R}$. Then, for all $s \in \mathbb{R}$ there exist $k=k_s \in \mathbb{N}_0$ and a constant $C=C_s>0$ such that
\[
\norma*{Au}_{B^s} \leq C \abs{a}_k^{(m)} \norma*{u}_{B^{s+m}}, \quad u \in B^{s+m}(\mathbb{R}^n).
\]
\end{proposition}
\begin{proof}
Let $E_{s+m} \in \Psi_{\mathrm{iso}}^{-s-m}(\mathbb{R}^n)$ be a global parametrix of $\Lambda^{s+m}$, i.e. $E_{s+m}\Lambda^{s+m}=I+R$, with $R \in \Psi_{\mathrm{iso}}^{-\infty}(\mathbb{R}^n)$.
 
\noindent Define
\[
A_0:=\Lambda^{s}AE\in \Psi_{\mathrm{iso}}^{0}(\mathbb{R}^n), \quad v:=\Lambda^{s+m}u.
\]
In this way (with $p \in \mathbb{N}, p \geq s$)
\[
\begin{split}
\norma*{Au}_{B^s}&=\norma*{\Lambda^sAu}_0+\sum_{\abs{\alpha}+\abs{\beta}\leq p}\norma*{x^\alpha D_x^\beta RAu}_0\\
&\leq \norma*{A_0 v}_0+\norma*{\Lambda^sARu}_0+\sum_{\abs{\alpha}+\abs{\beta}\leq p}\norma*{x^\alpha D_x^\beta RAu}_0.
\end{split}
\]
Now, writing $A_0=\mathrm{Op}^{\mathrm{w}}(a_0)$, by Theorem \ref{compAB}, we have that the map 
\[
S_{\mathrm{iso}}^m(\mathbb{R}^n) \ni a \mapsto a_0 \in S_\mathrm{iso}^0(\mathbb{R}^n),
\]
is continuous.
Therefore, by Theorem 18.6.3 in  \cite{HoV3} (see also \cite{P1}, Theorem 3.1.13) we get
\[
\norma*{A_0v}_0 \leq C_0 \abs{a_0}_{k_0}^{(0)}\norma*{v}_0
\]
and then (since $R \in \Psi_{\mathrm{iso}}^{-\infty}(\mathbb{R}^n)$)
\begin{equation}\label{apriori.Bsm}
\norma*{Au}_{B^s} \leq C \abs{a}_k^{(m)}\norma*{u}_{B^{s+m}},
\end{equation}
for some $k_0,k=k(s) \in \mathbb{N}_0$ and some constants $C_0,C=C(s)>0$.

\noindent Now, if $ u \in B^{s+m}(\mathbb{R}^n)$ we have that there exists a sequence $(u_j)_{j \geq 1}\subseteq \mathscr{S}(\mathbb{R}^n)$ such that 
\[
u_j \rightarrow u \quad \text{in} \ B^{s+m}(\mathbb{R}^n) \  \text{as} \ j \rightarrow +\infty.
\]
Thus, by the apriori estimate \eqref{apriori.Bsm} we have that $(Au_j)_{j\geq 1}$ is a Cauchy sequence in $B^s(\mathbb{R}^n)$ and then there exists $w \in B^{s}(\mathbb{R}^n)$ such that
\[
Au_j \rightarrow w \quad \text{in} \ B^{s}(\mathbb{R}^n) \  \text{as} \ j \rightarrow +\infty.
\]
Hence, since $Au_j \rightarrow Au$ in $\mathscr{S}'(\mathbb{R}^n)$ we get $w=Au$ and 
\[
\norma*{Au}_{B^s}=\lim_{j \rightarrow +\infty}\norma*{Au_j}_{B^s}\leq \lim_{j\rightarrow +\infty} C \abs{a}_k^{(m)}\norma*{u_j}_{B^{s+m}}=C \abs{a}_k^{(m)}\norma*{u}_{B^{s+m}}.
\]
\end{proof}

We next examine the relation between the Shubin-Sobolev spaces and the standard Sobolev spaces.

\begin{proposition}\label{prop.BsHs}
If $s\geq 0$ then 
\[
B^s(\mathbb{R}^n) \subseteq H^s(\mathbb{R}^n).
\]
\end{proposition} 
\begin{proof}
Let $u \in \mathscr{S}(\mathbb{R}^n)$ and write $E_s\Lambda^{s}=I+R$, where $E_s \in \Psi_{\mathrm{iso}}^{-s}(\mathbb{R}^n)$ and $R \in \Psi_{\mathrm{iso}}^{-\infty}(\mathbb{R}^n)$.
Denoting $\langle D \rangle^s=\mathrm{Op}^{\mathrm{w}}(\langle \xi \rangle^s)$ and using that $\Lambda^{-s} \in \Psi^{-s}_{\mathrm{iso}}(\mathbb{R}^n)\subseteq \Psi^{-s}(\mathbb{R}^n)$ (where the latter denotes the set of pseudodifferential operators in the Kohn-Nirenberg calculus of order $-s$) we get
\[
\norma*{u}_{H^s}=\norma*{\langle D \rangle^s u}_0
\leq \norma*{(\langle D \rangle^s E_s)\Lambda^s u}_0+\norma*{\langle D \rangle^sRu}_0 \leq C\norma*{u}_{B^s},
\]
for some constant $C>0$.
\end{proof}
We now examine the relation between the isotropic metric on $\mathbb{R}^{n+m}$ with respect to the \textit{iso,iso} metric and then the relation of the isotropic metric on $\mathbb{R}^n$ with respect to the $SG$-metric. 
We start with the relation between a symbol on the isotropic metric
$\mathbb{R}^{2(n+m)}_{X,Y}$ given by
\[
g_0=(\abs*{dX}^2+\abs*{dY}^2)/(1+\abs*{X}^2+\abs*{Y}^2),
\]
and a symbol for the split \textit{iso,iso}-metric in $\mathbb{R}^{2n}_X\times\mathbb{R}^{2m}_Y$
\[
g_2=\abs*{dX}^2/(1+\abs*{X}^2)+\abs*{dY}^2/(1+\abs*{Y}^2),
\]
defined as follows.
\begin{definition}\label{defSG_isoiso}
Let $a \in C^\infty(\mathbb{R}^n \times \mathbb{R}^m \times\mathbb{R}^n \times \mathbb{R}^m)$ and $(m_1,m_2)\in \mathbb{R}^2$.
 We say that \textit{$a$ is a iso,iso-symbol of order $(m_1,m_2)$}, denoted $a \in S_{\mathrm{iso,iso}}^{m_1,m_2}(\mathbb{R}^{n}\times \mathbb{R}^{m})$, 
if for all $\alpha_1,\beta_1,\alpha_2,\beta_2 \in \mathbb{N}_0^n$ 
there exists a constant $C=C(\alpha_1,\beta_1,\alpha_2,\beta_2)>0$ such that 
\begin{equation}
\abs*{\partial_x^{\alpha_1}\partial_{\xi}^{\beta_1}\partial_{y}^{\alpha_2}\partial_\eta^{\beta_2}a(x,y,\xi,\eta)}\leq C \braket{X}^{m_1-\abs{\alpha_1}-\abs{\beta_1}}\braket{Y}^{m_2-\abs{\alpha_2}-\abs{\beta_2}}, 
\end{equation}
for $X=(x,\xi) \in \mathbb{R}^n \times \mathbb{R}^n$, $Y=(y,\eta) \in \mathbb{R}^m \times \mathbb{R}^m$. 
We denote the corresponding set of \textit{iso,iso}-pseudodifferential operator
of order $(m_1,m_2)$ as 
\[
\Psi_{\mathrm{iso,iso}}^{m_1,m_2}(\mathbb{R}^n\times\mathbb{R}^m)=\lbrace A:\mathscr{S}'(\mathbb{R}^{n+m}) \rightarrow \mathscr{S}'(\mathbb{R}^{n+m}); \ \exists \, 
a \in S_{\mathrm{iso,iso}}^{m_1,m_2}(\mathbb{R}^{n}\times \mathbb{R}^{m}), \; A=\mathrm{Op}^{\mathrm{w}}(a) \rbrace.
\]
\end{definition}

We briefly describe the composition law in $S^{m_1,m_2}_{\mathrm{iso,iso}}(\mathbb{R}^n\times\mathbb{R}^m).$ Let $\sigma_n$, resp. $\sigma_m$, be the standard symplectic
form in $\mathbb{R}^{2n}$, resp. $\mathbb{R}^{2m}$. For $X,X'\in\mathbb{R}^{2n}$ and $Y,Y'\in\mathbb{R}^{2m}$, and 
$$\sigma(\left[\begin{array}{c}X\\ Y\end{array}\right],
\left[\begin{array}{c}X'\\ Y'\end{array}\right]):=\sigma_n(X,X')+\sigma_m(Y,Y'),$$ 
given symbols $a\in S^{m_1,m_2}_{\mathrm{iso,iso}}(\mathbb{R}^n\times\mathbb{R}^m)$ and 
$b\in S^{m'_1,m'_2}_{\mathrm{iso,iso}}(\mathbb{R}^n\times\mathbb{R}^m)$ one has (from H\"ormander's Weyl-calculus)
$$a\sharp b(X,Y)=e^{i\sigma(D_{X,Y},D_{X',Y'})/2}a(X,Y)b(X',Y')\bigl|_{(X,Y)=(X',Y')}\in S^{m_1+m'_1,m_2+m'_2}_{\mathrm{iso,iso}}(\mathbb{R}^n\times\mathbb{R}^m).$$

\begin{proposition}\label{propSGisoiso} Let $r\geq 0.$
One has that 
\[
S_{\mathrm{iso}}^{-r}(\mathbb{R}^{n}\times\mathbb{R}^m) \subseteq S_{\mathrm{iso,iso}}^{0,-r}(\mathbb{R}^{n}\times \mathbb{R}^{m}).
\]
\end{proposition}
\begin{proof}
Let $a=a(x,y,\xi,\eta) \in S_{\mathrm{iso}}^{-r}(\mathbb{R}^{n}\times\mathbb{R}^m)$. For all $\alpha=(\alpha_1,\alpha_2),\beta=(\beta_1,\beta_2) \in \mathbb{N}_0^{2n}$ one has 
\[
\begin{split}
\abs*{\partial_x^{\alpha_1}\partial_{\xi}^{\beta_1}\partial_{y}^{\alpha_2}\partial_\eta^{\beta_2}a(x,y,\xi,\eta)} & 
=\abs*{\partial_{(x,y)}^\alpha \partial_{(\xi,\eta)}^\beta a(x,y,\xi,\eta)}\leq C\braket {(x,y,\xi,\eta)}^{-r-\abs*{\alpha}-\abs*{\beta}} \\
&=C \braket{(x,y,\xi,\eta)}^{-\abs*{\alpha_1}-\abs*{\beta_1}}\braket{(x,y,\xi,\eta)}^{-r-\abs*{\alpha_2}-\abs*{\beta_2}} \\
& \leq C \braket{X}^{-\abs*{\alpha_1}-\abs*{\beta_1}}\braket{Y}^{-r-\abs*{\alpha_2}-\abs*{\beta_2}},
\end{split}
\]
for some $C>0$.
\end{proof}

Moreover, we recall the notion of $SG$-Symbol 
on $\mathbb{R}^n$, i.e. a symbol with respect to the metric
\[
g_{SG}=\frac{\abs{dx}^2}{1+\abs{x}^2}+\frac{\abs{d\xi}^2}{1+\abs{\xi}^2},
\]
and of $SG$-pseudodifferential 
operator (see \cite{CP} or \cite{SC}).

\begin{definition}\label{defSG}
Let $a \in C^\infty(\mathbb{R}^{2n})$ and $(m_1,m_2)\in \mathbb{R}^2$. 
We say that \textit{$a$ is a SG-symbol of order $(m_1,m_2)$},
denoted $a \in S^{m_1,m_2}(\mathbb{R}^{n})$,
if for all $\alpha,\beta \in \mathbb{N}_0^n$ there exists a constant $C=C(\alpha,\beta)>0$
such that 
\begin{equation}
\abs*{\partial_x^{\beta}\partial_{\xi}^{\alpha}a(x,\xi)}\leq C \braket{x}^{m_1-\abs{\beta}}\braket{\xi}^{m_2-\abs{\alpha}}.
\end{equation}
In addition, we define the set 
of \textit{$SG$-pseudodifferential operator
of order $(m_1,m_2)$} as 
\[
\Psi_{\mathrm{SG}}^{m_1,m_2}(\mathbb{R}^n)=\lbrace A:\mathscr{S}'(\mathbb{R}^n) \rightarrow \mathscr{S}'(\mathbb{R}^n); \ \exists \, a \in S^{m_1,m_2}(\mathbb{R}^n), \; A=\mathrm{Op}(a) \rbrace.
\]
\end{definition}

\begin{proposition}\label{propSG}
If $\mu,\mu_1,\mu_2 \geq 0$ are such that $\mu=\mu_1+\mu_2$ one has
\[
S_{\mathrm{iso}}^{-\mu}(\mathbb{R}^{n}) \subseteq S_{\mathrm{SG}}^{-\mu_1,-\mu_2}(\mathbb{R}^{n}).
\]
\end{proposition}
\begin{proof}
The result follows similarly to Proposition \ref{propSGisoiso} noting that, for all $\alpha,\beta \in \mathbb{N}_0^n$, there exists $C>0$ such that 
\[
\begin{split}
\abs*{\partial_x^{\beta}\partial_{\xi}^{\alpha}a(x,\xi)}&\leq C \braket{(x,\xi)}^{-\mu-\abs*{\alpha}-\abs*{\beta}} \\
&= C \braket{(x,\xi)}^{-\mu_1-\abs*{\beta}} \braket{(x,\xi)}^{-\mu_2-\abs*{\alpha}}  \\
& \leq C \braket{x}^{-\mu_1-\abs*{\beta}}\braket{\xi}^{-\mu_2-\abs*{\alpha}}.
\end{split}
\]
\end{proof}

Finally, the $SG$-Sobolev space of order $(m_1,m_2) \in \mathbb{R}^2$ 
is defined as
\[
H^{m_1,m_2}(\mathbb{R}^n)= \lbrace u \in \mathscr{S}'(\mathbb{R}^n); \ \mathrm{Op}(\langle x \rangle^{m_1}\langle \xi \rangle^{m_2})u \in L^2(\mathbb{R}^n) \rbrace
\]
and the following boundedness property holds. 

\begin{proposition}\label{boundSG}
Let $A \in \Psi_{\mathrm{SG}}^{m_1,m_2}(\mathbb{R}^n)$. Then, for all $(s_1,s_2) \in \mathbb{R}^2$,  
\[
A:H^{s_1,s_2}(\mathbb{R}^n)\rightarrow H^{s_1-m_1,s_2-m_2}(\mathbb{R}^n)
\]
is bounded.
\end{proposition}

By using these spaces, through the Proposition \ref{propSG}, we reach the following result.
 
\begin{proposition}\label{prop.1in-mu}
If $\mu_1,\mu_2 \geq 0$ and $\mu=\mu_1+\mu_2$, then
\begin{equation}\label{eq.BSGneg}
H^{-\mu_1,-\mu_2}(\mathbb{R}^n) \subseteq B^{-\mu}(\mathbb{R}^n).
\end{equation}
In particular, if $\mu>n/2$,
\begin{equation}\label{eq.1in-mu}
\mathbbm{1}_x \in B^{-\mu}(\mathbb{R}^n).
\end{equation}
\end{proposition}
\begin{proof}
Let $u \in H^{-\mu_1,-\mu_2}(\mathbb{R}^n)$.
By definition of the space $B^{-\mu}(\mathbb{R}^n)$ we have to prove that
\[
\Lambda^{-\mu}u\in L^2(\mathbb{R}^n),
\]
and since $\Lambda^{-\mu} \in \Psi_{\mathrm{SG}}^{-\mu_1,-\mu_2}(\mathbb{R}^n)$
(see Proposition \ref{propSG}), 
the result follows
by the boundedness properties on the $SG$-Sobolev spaces 
(see Proposition \ref{boundSG}). 
Consequently, to obtain \eqref{eq.1in-mu} it is sufficient to notice that for $\mu>n/2$ one has 
\[
\mathbbm{1}_x \in H^{-\mu,0}(\mathbb{R}^n).
\]
\end{proof}
Finally, it is possible to give the following characterization of the space 
\[
B^{-k}(\mathbb{R}^n)=(B^k(\mathbb{R}^n))^\ast, \quad  \text{for} \ k \in \mathbb{N}.
\]
\begin{remark}\label{rmksum}
If $u \in \mathscr{S}'(\mathbb{R}^n)$, there exists $k \in \mathbb{N}_0$, such that $u \in B^{-k}(\mathbb{R}^n)=B^k(\mathbb{R}^n)^\ast$. Thus, by the Riesz Representation Theorem, we may find $u_k \in B^{k}(\mathbb{R}^n)$ that satisfies
\[
u(\varphi)=(\varphi,u_k)_{B^k}, \quad \forall \varphi \in B^k(\mathbb{R}^n).
\]
Thus, if $\varphi \in \mathscr{S}(\mathbb{R}^n)$ we have 
\[
\begin{split}
u(\varphi) & =\sum_{\abs*{\alpha}+\abs*{\beta}\leq k}(x^\alpha \partial_x^\beta \varphi,x^\alpha \partial_x^\beta u_k)_0 \\
&=\sum_{\abs*{\alpha}+\abs*{\beta}\leq k} \langle (-1)^{\abs*{\beta}}\partial_x^\beta x^{2\alpha}\partial_x^\beta \bar{u}_k|\varphi \rangle.
\end{split}
\]
Hence, since $B^{-s}(\mathbb{R}^n)\subseteq B^{-s'}(\mathbb{R}^n)$, if $s'>s$, we may write, for some $N \in \mathbb{N}_0$,
\[
u=\sum_{\abs{\alpha}+\abs{\beta}\leq N}D^\beta ( x^\alpha u_{\alpha\beta}),
\]
where $u_{\alpha\beta}\in B^{s}(\mathbb{R}^n)$ for $s>0$ sufficiently large. 
\end{remark}

\vspace{.5cm}
\noindent\textbf{Acknoledgments:} The authors wish to thank Jeffrey Galkowski for useful discussions. The authors wish also to thank the referees for their
thorough reading of the paper, precise observations and suggestions that led to a substantial improvement in clarity and precision.

\vspace{.5cm}
\noindent\textbf{Data availability statement.} There is no associated data to the manuscript.

\vspace{.5cm}
\noindent\textbf{Declarations}\\
\textbf{Conflict of interest.} There is no conflict of interest.



\begin{thebibliography}{10}


\bibitem {CRW} M.~Cappiello, L.~G. Rodino and P.~Wahlberg. 
\newblock Propagation of anisotropic Gabor singularities for Schr\"odinger type equations.
\newblock J. Evol. Equ. {\bf 24}(2024), Paper No. 36, 46 pp. https://doi.org/10.1007/s00028-024-00963-w


\bibitem {CW} E.~Carypis and P.~Wahlberg.
\newblock Propagation of exponential phase space singularities for Schr\"odinger equations with quadratic Hamiltonians,
\newblock J. Fourier Anal. Appl. {\bf 23}(2017), 530--571. https://doi.org/10.1007/s00041-016-9478-6


\bibitem{SC}S.~Coriasco.
\newblock Fourier integral operators in SG classes. I: Composition theorems and action on SG Sobolev spaces.
\newblock Rendiconti del Seminario Matematico {\bf 57}(1999), 249--302. 

\bibitem{CM} S.~Coriasco and L.~Maniccia.
\newblock Wave front set at infinity and hyperbolic linear operators with multiple characteristics.
\newblock Ann. Global Anal. Geom. {\bf 24}(2003), 375--400. https://doi.org/10.1023/A:1026241614722


\bibitem{DGW} M.~Doll, O.~Gannot and J.~Wunsch. 
\newblock Refined Weyl Law for Homogeneous Perturbations of the Harmonic Oscillator.
\newblock Communications in Mathematical Physics {\bf 362}(2017), 269--294. https://doi.org/10.1007/s00220-018-3100-5

\bibitem{G} K.~Gröchenig.
\newblock \textit{Foundations of Time-Frequency Analysis}.
\newblock Birkhäuser, Boston (2001)

\bibitem{H} B.~Helffer. 
\newblock \textit{Th\'eorie spectrale pour des op\'erateurs globalement elliptiques}. 
\newblock Astérisque 112 (1984)

\bibitem{HPS} M.~Hitrik, K.~Pravda-Starov.
\newblock Spectra and semigroup smoothing for non-elliptic quadratic operators.
\newblock Math. Ann. {\bf 344}(2009), 801--846. https://doi.org/10.1007/s00208-008-0328-y

\bibitem{HQ} L.~H\"ormander. 
\newblock Quadratic Hyperbolic Operators. 
\newblock Microlocal Analysis and Applications. Lecture Notes in Mathematics {\bf 1495}, 118--160. Springer, Berlin, Heidelberg (1991). https://doi.org/10.1007/BFb0085123

\bibitem{HS} L.~H\"ormander.
\newblock Symplectic classification of quadratic forms, and general Mehler formulas. 
\newblock Math Z {\bf 219}(1995), 413--449. https://doi.org/10.1007/BF02572374


\bibitem{HoV3} L.~H\"ormander. 
\newblock \textit{The Analysis of Linear Partial Differential Operators III}. 
\newblock Classics in Mathematics. Springer, Berlin (2007)

\bibitem{Le} N.~Lerner.
\newblock\textit{Metrics on the phase space and non-selfadjoint pseudo-differential operators}.
\newblock
Pseudo-Differential Operators. Theory and Applications, {\bf 3}, Birkh\"auser Verlag, Basel, 2010

\bibitem{NR} F.~Nicola, L.~Rodino.
\newblock \textit{Global Pseudo-Differential Calculus on Euclidean Spaces}.
\newblock Birkhäuser, Basel (2010)

\bibitem{CP} C.~Parenti.
\newblock Operatori pseudo-differenziali in Rn e applicazioni.
\newblock Annali di Matematica {\bf 93}(1972), 359--389. https://doi.org/10.1007/BF02412028

\bibitem{P1} A.~Parmeggiani. 
\newblock \textit{Spectral theory of Non-Commutative Harmonic Oscillators: An Introduction}. 
\newblock Lecture Notes in Mathematics {\bf 1992}. Springer-Verlag, Berlin, (2010).
xii+254 pp. doi: 10.1007/978-3-642-11922-4

\bibitem{PRW} K.~Pravda-Starov, L.~Rodino and P.~Wahlberg.
\newblock Propagation of Gabor singularities for Schrödinger equations with quadratic Hamiltonians.
\newblock Math. Nachr. {\bf 291}(2018), 128--159. https://doi.org/10.1002/mana.201600410

\bibitem{RW} L.~Rodino, P.~Wahlberg. 
\newblock The Gabor wave front set.
\newblock Monatshefte für Mathematik {\bf 173}(2012), 625--655. https://doi.org/10.1007/s00605-013-0592-0

\bibitem{SH} M.~A.~Shubin.
\newblock \textit{Pseudodifferential Operators and Spectral Theory}, 2nd Edition.
\newblock Springer-Verlag, Berlin (2001)

\bibitem{SW} R.~Schulz, P.~Wahlberg. 
\newblock Microlocal properties of Shubin pseudodifferential and localization operators.
\newblock J. Pseudo-Differ. Oper. Appl. {\bf 7}(2015), 91--111. https://doi.org/10.1007/s11868-015-0143-7

\bibitem {W} P.~Wahlberg. 
\newblock Propagation of polynomial phase space singularities for Schr\"odinger equations with quadratic Hamiltonians.
\newblock Math. Scand.  {\bf 122}(2018), 107--140. https://doi.org/10.7146/math.scand.a-97187

\bibitem{W1} P.~Wahlberg. 
\newblock Propagation of anisotropic Gabor wave front sets.
\newblock Proc. Edinb. Math. Soc. (2) {\bf 67}(2024), 674--698. https://doi.org/10.1017/S0013091524000269 

\bibitem{Wh} F.~White.
\newblock Propagation of global analytic singularities for Schr\"odinger equations with quadratic Hamiltonians.
\newblock J. Funct. Anal. \textbf{283}(6) (2022), Paper. No. 109569, 45pp. https://doi.org/10.1016/j.jfa.2022.109569  

\bibitem{Z} M. Zworski.
\newblock \textit{Semiclassical analysis}.
\newblock Graduate Studies in Mathematics, {\bf 138}, Amer. Math. Soc., Providence, RI, 2012

\end{thebibliography}
\end{document}